\documentclass[11pt]{article}

\usepackage{hyperref}
\usepackage[utf8]{inputenc}
\usepackage[cyr]{aeguill}
\usepackage[all]{xy}
\usepackage[french, english]{babel}
\usepackage[T1]{fontenc}
\usepackage{amscd}  
\usepackage{amsmath}
\usepackage{amsthm} 
\usepackage{amssymb}
\usepackage{mathabx}
\usepackage{xcolor}
\usepackage{titlesec}
\usepackage{stmaryrd}
\usepackage{enumerate}

\newtheorem{observation}{Remark}[section]
\newtheorem{lemma}[observation]{Lemma}  
\newtheorem{theorem}[observation]{Theorem}
\newtheorem{definition}[observation]{Definition}
\newtheorem{example}[observation]{Example}

\newtheorem{proposition}[observation]{Proposition} 
\newtheorem{corollary}[observation]{Corollary}

\titleformat{\section}
     {\center\Large\scshape}{\thesection.}{1em}{}
 \titleformat{\subsection}[runin]
     {\scshape}{\thesubsection. }{0em}{}
     [. ]
     
\def\-{\operatorname{-}}

\def\Z{\mathbb{Z}}
\def\X{\mathbb{X}}
\def\L{\mathbb{L}}

\def\com{\mathrm{CMon}}
\def\Lie{\mathrm{Lie}}

\def\GRD{\mathsf{GRD}}

\def\P{\mathcal P}

\def\CMON{\mathsf{CMON}}
\def\PMAG{\mathsf{PMAG}}
\def\MON{\mathsf{MON}}
\def\AbG{\mathsf{AbG}}

\newcommand{\Ab}{\mathsf{AB}}
\def\ab{\mathrm{Ab}}

\def\Grp{\mathsf{GRP}}

\def\∑{\Sigma}
\def\diag{%\shorthandoff{;:!?}
\xymatrix}
\def\0{\mathbf 0}
\def\_{\underline}

\author{Sacha Ikonicoff\footnote{Institut de recherche mathématique avancée, Université de Strasbourg, Strasbourg, France}, Jean-Simon Pacaud Lemay\footnote{School of Mathematical and Physical Sciences, Macquarie University, Sydney, New South Wales, Australia} and Tim~Van~der~Linden\footnote{Institut de Recherche en Mathématique et Physique, Université catholique de Louvain, Louvain-la-Neuve, and Mathematics and Data Science, Vrije Universiteit Brussel, Brussel, Belgium}}

\title{From Abelianization to Tangent Categories}

\begin{document}
\allowdisplaybreaks

\maketitle

\begin{abstract}
    A tangent category is a category with an endofunctor, called the tangent bundle functor, which is equipped with various natural transformations that capture essential properties of the classical tangent bundle of smooth manifolds. In this paper, we show that, surprisingly, the category of groups is a tangent category whose tangent bundle functor is induced by abelianization and whose differential bundles correspond to abelian groups. We generalize this construction by introducing the concept of linear assignments, which are endofunctors assigning to every object a commutative monoid in a natural and idempotent manner. We then show that a linear assignment induces a tangent bundle functor, whose differential bundles correspond to a notion of linear algebras. We show that any finitely cocomplete regular unital category is a tangent category whose tangent bundle functor is induced by the canonical abelianization functor, which is a monadic linear assignment. This allows us to provide multiple new examples of tangent categories including monoids, pointed magmas, loops, non-unital rings, Jónsson--Tarski varieties, and pointed Mal'tsev varieties.
\end{abstract}
\small{Keywords: Abelianization, Linear Assignment, Linear Projector, Tangent Category, Unital Category}
\tableofcontents

\section{Introduction}

Tangent categories provide a categorical framework for the foundations of differential calculus over smooth manifolds by abstracting the notion of the tangent bundle. Tangent categories were originally introduced by Rosický in \cite{RosickyTangentCats}, then later rediscovered and further developed by Cockett and Cruttwell in \cite{cockett2014differential}. Briefly, a tangent category is a category equipped with an endofunctor \(\mathcal{T}\), called the tangent bundle functor, where we interpret \(\mathcal{T}(X)\) as an abstract tangent bundle over an object \(X\), and equipped with various natural transformations that capture essential properties of the classical tangent bundle of smooth manifolds. The theory of tangent categories is now a well-established field, having been able to formalize various important concepts in differential geometry, and has applications in various mathematical domains of study, including algebra, algebraic geometry,  and operad theory. There are many interesting examples of tangent categories such as the category of smooth manifolds, the category of (affine) schemes, and the category of commutative algebras.

The observation at the origin of this paper is that, surprisingly, the category of groups is a tangent category whose tangent bundle functor is induced by abelianization. Indeed, as we will see in Example \ref{ex tanL-group}, the tangent bundle functor sends a group \(G\) to the product of \(G\) with its abelianization \(\ab(G)\):
\[ \mathcal{T}(G) = G \times \ab(G) \]
This tangent structure on groups does not fall within the geometric flavoured style of tangent categories, such as smooth manifolds and schemes. Indeed, this tangent structure does not bear much geometric information: it is `everywhere flat', in the sense that for any group \(G\), the fibre over a point \(x\in G\) is always \(\ab(G)\), so it does not depend on \(x\). Moreover, the derivative of a group morphism can be seen as being `simply' its abelianization. That said, while this structure is rather simple from a geometric point of view, it still allows us to interpret the abelianization of a group as a generic fiber space for this group, or in other words, a linear space of directions over this group. As such, this tangent structure on groups instead falls more in line with the more algebraic style examples, such as commutative algebras (which was one of Rosický's original tangent category examples) or categories with finite biproducts, which are of indisputable importance and interest.

The objective of this paper is to extract the properties of the category of groups and of the abelianization functor which allow us to build such a tangent structure. A key part in the definition of a tangent category is that the fibres of the tangent bundle are commutative monoids. In the category of groups, commutative monoid objects correspond precisely to abelian groups. As such, abelianization sends objects to commutative monoids. Moreover, abelianization preserves finite products and is idempotent, that is, abelianizing twice is equivalent to abelianizing only once. These three facts about abelianization are the essential properties which induce a tangent structure. In order to generalize this, we introduce the notion of a \textbf{linear assignment} (Definition \ref{def linear-ass}) on a category with finite products. A linear assignment is an endofunctor which preserves finite products,  and which sends each object to an object equipped with a commutative monoid structure, in a natural and idempotent manner. Alternatively (Theorem \ref{thm ass=proj}), a linear assignment can also be described as a finite product preserving functor from a category with finite products to its category of commutative monoids, which is idempotent in a suitable sense. These are called \textbf{linear projectors} (Definition \ref{def linear-proj}). Of course, one may observe that abelian groups actually correspond to internal abelian groups in the category of groups, and thus, abelianization sends objects to abelian group objects. In fact, the category of groups is a Rosický tangent category, by which we mean a tangent category such that the fibres of the tangent bundles are abelian groups. We then say that an \textbf{additive assignment (resp.\ projector)} is a linear assignment (resp.\ projector) which, furthermore, sends objects to abelian group objects.

The main result of this paper (Theorem \ref{theo:linearprojtotangent}) states that a category with finite products and a linear (resp.\ additive) assignment \(\mathcal{L}\) is a cartesian (Rosický) tangent category with tangent bundle functor given by:
\[ \mathcal{T}(X) = X \times \mathcal{L}(X). \]
We then investigate differential bundles \cite{cockettCruttwellDiffBundles} in such a tangent category. Intuitively, a differential bundle over an object formalizes the notion of a smooth vector bundle over a smooth manifold. Differential bundles over the terminal objects are called \textbf{differential objects}, which capture the notion of Euclidean spaces in a tangent category. For a linear assignment \(\mathcal{L}\), differential bundles and differential objects are very closely related to what we call \textbf{linear algebras} (Definition \ref{def L-alg}), or simply \(\mathcal{L}\)-algebras, which are objects that are isomorphic to their associated commutative monoid. Indeed, we show that for a linear assignment \(\mathcal{L}\), differential objects correspond precisely to \(\mathcal{L}\)-algebras (Theorem \ref{thm diff-obj=l-alg}). We also show that, for an object \(X\), the product of \(X\) with an \(\mathcal{L}\)-algebra is a differential bundle over \(X\) (Proposition \ref{prop Lalg-to-dbun}). When the base category admits zero morphisms and kernels, we also show that every differential bundle over \(X\) is the product of \(X\) with an \(\mathcal{L}\)-algebra (Theorem \ref{thm ker-zero-morph}). In our motivating example, for the abelianization of groups, the linear algebras are the abelian groups (Example \ref{ex L-alg-group}). Therefore, in the category of groups (which admits zero morphisms and kernels), the differential objects correspond precisely to abelian groups, while a differential bundle over a group \(G\) is a product of \(G\) with an abelian group (Example \ref{ex dbun-group}).

The abelianization functor for groups has a very rich structure, beyond being an additive projector. For instance, abelianization on groups has the structure of a monad, and moreover, of an idempotent monad. In general, a linear assignment is not a monad, but it is an idempotent non-unital monad (also sometimes called a semimonad), as there does not need to be a natural transformation \(X \to \mathcal{L}(X)\) playing the role of the unit. A linear assignment with such a unit natural transformation is called a \textbf{monadic linear assignment} (Definition \ref{def monadic-L}). Abelianization on groups is an example of monadic linear assignment. Furthermore, for a monadic linear assignment, the linear algebras correspond precisely to the usual notation of algebras over a monad (Lemma \ref{lemma monad-L-alg}), justifying the terminology of linear algebra and the notation \(\mathcal{L}\)-algebra. In fact, since a monadic linear assignment is an idempotent monad, being a linear algebra is a property of an object rather than an additional structure, and so, the category of linear algebras is in fact a reflective subcategory. Thus, in the same way that idempotent monads correspond to reflective subcategory and reflectors, monadic linear assignments correspond to \textbf{linear reflective subcategories} and \textbf{linear reflectors} (Definition \ref{def L-reflector}), which are reflective subcategories where the product becomes a biproduct and the reflector preserves finite products.

Abelianization is a special kind of monadic linear assignment, since its linear algebras correspond to the internal commutative monoids. This implies that the category of commutative monoids is a linear reflective subcategory. In fact, this observation can be generalized in the context of (regular) unital categories. Unital categories were introduced by Bourn in \cite{Bourn1996} in their study of categorical algebra. Unital categories, along with Mal'tsev categories, are well-studied and have a rich literature. It turns out that in a unital category, being a commutative monoid is a property of an object rather than an additional structure. Thus, for a unital category, we may view its category of commutative monoids as a full subcategory. Moreover, for a finitely cocomplete regular unital category, commutative monoids form a reflective subcategory \cite{borceuxbourn04}, and we show that the reflector induces a monadic linear assignment (Proposition \ref{prop unital-linear-reflector}). Therefore, it follows that every finitely cocomplete regular unital category is a cartesian tangent category whose differential bundles correspond precisely to commutative monoid objects (Theorem \ref{thm unital-tan}). This allows us to provide new examples of tangent categories, including monoids (Example \ref{ex monoids}), pointed magmas (Example \ref{ex pmag}), and more generally any Jónsson--Tarski variety.

Going even deeper, the category of groups is strongly unital \cite{borceuxbourn04}. For strongly unital categories, it turns out that every commutative monoid is in fact an abelian group. Thus, for a finitely cocomplete, regular and strongly unital category, abelian groups are a reflective subcategory, which induces a monadic additive assignment, which we may refer to as abelianization. Therefore, every finitely cocomplete, regular, and strongly unital category is a cartesian Rosický tangent category whose tangent bundle functor is induced by abelianization. This leads us to many new interesting examples of Rosický tangent categories coming from algebra, including groups, non-unital rings (Example \ref{ex rings}), crossed modules (Example \ref{ex crossed}), loops (Example \ref{ex loop}), and more generally, any semiabelian category, and also any pointed Mal'tsev variety (Example \ref{ex Mal}).

We conclude this introduction with a brief discussion of potential future work and ideas to investigate, especially around the concept of `parallelizable objects'. A key feature of the tangent structure induced by a linear assignment is that every object is in fact `parallelizable', meaning that the tangent bundle of each object is the product of the base object with a differential object. However, at the time of writing this paper, the theory of parallelizable objects in a tangent category is part of the folklore, and has yet to be properly developed. Moreover, the tangent bundle of a differential object is simply the product of two copies of the given differential object, thus differential objects are special kinds of parallelizable objects. A tangent category where every object is a differential object is precisely a cartesian differential category, as introduced by Blute, Cockett, and Seely in \cite{blute2009cartesian}. It should be the case that any tangent category where every object is a parallelizable object is a generalized cartesian differential category, as introduced by Cruttwell in \cite{cruttwell17}. Thus, categories with linear assignments should be generalized cartesian differential categories. However, Cruttwell's definition asks for certain strict equalities, while we can only provide isomorphisms. These isomorphisms are quite subtle and appear when showing that we produce a tangent structure. Thus, some work would also be needed to give the proper definition of `non-strict' generalized cartesian differential categories. Lastly, the category of abelian groups is equivalent to the subcategory of differential objects in the tangent category of affine schemes \cite{cruttwell2023differential}. It would be interesting to see if it is possible to build a larger tangent category whose subcategory of parallelizable objects is equivalent to the category of groups.

\section*{Acknowledgement}
The authors would like to thank Steve Lack for providing us with Example \ref{example Lack}, and for fruitful conversations. The second named author is funded by an ARC DECRA award (\# DE230100303) and this material is based upon work supported by the AFOSR under award number FA9550-24-1-0008. The third named author is a Senior Research Associate of the Fonds de la Recherche Scientifique--FNRS.

\section{Linear Assignments and Projectors}
In this section, we introduce the notions of linear (resp.\ additive) assignments and projectors. Essentially, a linear (resp.\ additive) assignment is a way of associating to every object a commutative monoid (resp.\ abelian group), such that this association is natural and idempotent. Our motivating example is given by abelianization, sending groups to abelian groups. Linear (resp.\ additive) assignments are the key ingredient for our construction of (Rosický) tangent structures. Linear (resp.\ additive) projectors are an equivalent way of describing linear (resp.\ additive) assignments, by factoring through the category of commutative monoid objects (resp.\ abelian group objects).

If only to introduce notation, we begin by recalling the definitions of commutative monoid and abelian group objects in a category with finite products. For a category with finite products, we denote the product by \(\times\), projections by \(\pi_j\colon X_1 \times X_2 \to X_j\), pairing by \(\langle -, - \rangle\), the terminal object by \(\ast\), and the unique morphism to the terminal object by \(t_X \colon X \to \ast\). Recall that a \textbf{commutative monoid} in a category with finite products \(\mathbb{X}\) is a triple \(\mathsf{M} = (M, \bullet, e)\) consisting of an object \(M\) equipped with morphisms \(\bullet\colon M \times M \to M\) and \(e\colon \ast \to M\) such that the following diagrams commute:
\begin{equation*}\begin{gathered}\label{diag monoid}
        \xymatrixcolsep{1.9pc}\xymatrix{ M \times (M \times M)  \ar[d]_-{\alpha_{M,M,M}} \ar[r]^-{1_M \times \bullet} & M \times M \ar[dd]^-{\bullet}   & M  \ar@{=}[ddr]_-{} \ar[d]_-{\rho_M} \ar[r]^-{\lambda_M} & \ast \times M \ar[d]^-{e \times 1_M}   & M \times M  \ar[dr]_-{\bullet} \ar[r]^-{\sigma_{M,M}}& M \times M \ar[d]^-{\bullet}   \\
        (M \times M) \times M  \ar[d]_-{\bullet \times 1_M} & & M \times \ast  \ar[d]_-{1_M \times e} & M \times M \ar[d]^-{\bullet} &    & M  \\
        M \times M \ar[r]_-{\bullet} & M &  M \times M \ar[r]_-{\bullet} & M   }
    \end{gathered}\end{equation*}
where \(\alpha\), \(\lambda\) and \(\rho\), and \(\sigma\) are respectively the canonical associativity, unit, and symmetry natural isomorphism for the product. For commutative monoids \(\mathsf{M} = (M, \bullet, e)\) and \(\mathsf{M}^\prime = (M^\prime, \bullet^\prime, e^\prime)\), a \textbf{monoid morphism} \(f\colon \mathsf{M} \to \mathsf{M}^\prime\) is a morphism \(f\colon M \to M^\prime\) such that the following diagrams commute:
\begin{equation*}\begin{gathered}\label{diag monoid-morph}
        \xymatrixcolsep{5pc}\xymatrix{ M \times M  \ar[d]_-{f \times f} \ar[r]^-{\bullet} & M \ar[d]^-{f} & \ast \ar[r]^-{e} \ar[dr]_-{e^\prime} & M  \ar[d]^-{f} \\
        M^\prime \times M^\prime \ar[r]_-{\bullet^\prime}                                 & M^\prime      &                                      & M^\prime }
    \end{gathered}\end{equation*}
We denote the category of commutative monoids of \(\X\) and monoid morphisms between them by \(\CMON[\X]\), and let \(\mathcal{U}\colon \CMON[\X]\to \X\) be the forgetful functor, which sends a commutative monoid \(\mathsf{M}=(M,\bullet,e)\) to its underlying object, \(\mathcal{U}(\mathsf{M}) = M\), and a monoid morphism to itself, \(\mathcal{U}(f) = f\). It is easy to see that this functor is conservative, which means that \(f\colon \mathsf{M} \to \mathsf{M}^\prime\) is a monoid isomorphism if and only if the underlying morphism \(f\colon M \to M^\prime\) is an isomorphism. Recall that the forgetful functor creates products in \(\CMON[\X]\). As such, the terminal object is a commutative monoid \(\star = (\ast, t_{\ast\times \ast}, 1_\ast)\), and for two commutative monoids \(\mathsf{M} = (M, \bullet, e)\) and \(\mathsf{M}^\prime = (M^\prime, \bullet^\prime, e^\prime)\), their product is defined as:
\begin{equation*}\begin{gathered}\label{monoidproduct}
        \mathsf{M} \times \mathsf{M}^\prime = \left( M \times M^\prime, (\bullet \times \bullet^\prime) \circ \tau_{M,M^\prime,M,M^\prime}, \langle e, e^\prime \rangle \right)
    \end{gathered}\end{equation*}
where the morphism \(\tau_{X,Y,Z,W}\colon (X \times Y) \times (Z \times W) \to (X \times Z) \times (Y \times W) \) is the canonical natural interchange isomorphism of the product, that is, the isomorphism which swaps the middle two objects. In fact, \(\times\) is a biproduct and \(\star\) is a zero object in \(\CMON[\X]\).

An \textbf{abelian group} in a category with finite products \(\mathbb{X}\) is a quadruple \(\mathsf{G} = (G, \bullet, e, \iota)\) consisting of a commutative monoid \((G, \bullet, e)\) with a morphism \(\iota\colon G \to G\) such that the following diagram commutes:
\begin{equation*}\begin{gathered}\label{diag group}
        \xymatrixcolsep{5pc}\xymatrix{ & G \times G \ar[dr]^-{\bullet} \\
        G  \ar[ur]^-{\langle 1_G, \iota \rangle} \ar[dr]_-{\langle \iota, 1_G \rangle} \ar[r]^-{t_G} & \ast \ar[r]^-{e} &  G. \\
        & G \times G \ar[ur]_-{\bullet} }
    \end{gathered}\end{equation*}
If a commutative monoid admits such a morphism \(\iota\), this morphism is unique. As such, being an abelian group is a property of a commutative monoid rather than an additional structure. For abelian groups \(\mathsf{G} = (G, \bullet, e, \iota)\) and \(\mathsf{G}^\prime = (G^\prime, \bullet^\prime, e^\prime, \iota^\prime)\), a \textbf{group (iso)morphism} \(f\colon \mathsf{G} \to \mathsf{G}^\prime\) is simply a monoid (iso)morphism \(f\colon (G, \bullet, e) \to (G^\prime, \bullet^\prime, e^\prime)\), since the following diagram automatically commutes:
\begin{equation*}\begin{gathered}\label{diag group-morph}
        \xymatrixcolsep{5pc}\xymatrix{ G  \ar[d]_-{f} \ar[r]^-{\iota} & G  \ar[d]^-{f} \\
        G^\prime  \ar[r]_-{\iota^\prime} & G^\prime }
    \end{gathered}\end{equation*}
We denote the category of abelian group of \(\X\) and group morphisms between them by \(\mathsf{AbG}[\mathbb{X}]\), and, abusing notation, we also write \(\mathcal{U}\colon \mathsf{AbG}[\mathbb{X}]\to \X\) for the forgetful functor, which sends an abelian group to its underlying object and a group morphism to itself. As before, the forgetful functor creates products in \(\mathsf{AbG}[\mathbb{X}]\), which are in fact biproducts. Thus, the terminal object is an abelian group, \(\star = (\ast, t_{\ast\times \ast}, 1_\ast, 1_\ast)\), and for two abelian groups \(\mathsf{G} = (G, \bullet, e, \iota)\) and \(\mathsf{G}^\prime = (G^\prime, \bullet^\prime, e^\prime, \iota^\prime)\), their product is defined by:
\begin{equation*}\begin{gathered}\label{groupproduct}
        \mathsf{G} \times \mathsf{G}^\prime = \left( G \times G^\prime, (\bullet \times \bullet^\prime) \circ \tau_{G,G^\prime, G, G^\prime}, \langle e, e^\prime \rangle, \iota \times \iota^\prime \right).
    \end{gathered}\end{equation*}

Recall that a functor \(\mathcal{F}\colon \mathbb{X} \to \mathbb{Y}\) between categories with finite products is said to \textbf{preserve finite products} if \(t_{F(\ast)}\colon F(\ast) \to \ast\) is an isomorphism, and if the canonical natural transformation \(\omega_{X,Y}\colon \mathcal{F}(X \times Y) \to \mathcal{F}(X) \times \mathcal{F}(Y)\), defined by \(\omega_{X,Y} = \langle \mathcal{F}(\pi_1), \mathcal{F}(\pi_2) \rangle\), is a natural isomorphism.

\begin{definition}\label{def linear-ass} A \textbf{linear} (resp.\ \textbf{additive}) \textbf{assignment} on a category \(\mathbb{X}\) with finite products is a quadruple \((\mathcal{L}, +, 0, \nu)\) consisting of a finite product preserving endofunctor \(\mathcal{L}\colon \mathbb{X} \to \mathbb{X}\), natural transformations \(+_X\colon \mathcal{L}(X) \times \mathcal{L}(X) \to \mathcal{L}(X)\) and \(0_X\colon \ast \to \mathcal{L}(X)\) (and \(-_X\colon \mathcal{L}(X) \to \mathcal{L}(X)\)), and a natural isomorphism \(\nu_X\colon \mathcal{L}\mathcal{L}(X) \to \mathcal{L}(X)\), such that:
    \begin{enumerate}[{\em (i)}]
        \item For each object \(X\), the triple \(\mathsf{L}(X) = (\mathcal{L}(X), +_X, 0_X)\) is a commutative monoid (resp.\ the quadruple \(\mathsf{L}(X) = (\mathcal{L}(X), +_X, 0_X, -_X)\) is an abelian group);
        \item For each object \(X\), \(\nu_X\colon \mathsf{L}(\mathcal{L}(X)) \to \mathsf{L}(X)\) is a monoid (resp.\ group) isomorphism;
        \item For each object \(X\), \(\nu_{\mathcal{L}(X)} = \mathcal{L}(\nu_X)\).
    \end{enumerate}
    As a shorthand, when there is no confusion, we will denote linear and additive assignments simply by their underlying endofunctor \(\mathcal{L}\).
\end{definition}

By definition, every additive assignment is a linear assignment. On the other hand, since inverses, if they exist, are unique for commutative monoids, there exists at most one natural transformation which makes a linear assignment into an additive assignment. Thus, being an additive assignment is a property of a linear assignment, rather than an additional structure.

We first observe that naturality of the monoid structure implies that linear (resp.\ additive) assignments send morphisms to monoid (resp.\ group) morphisms. It follows that the product preserving isomorphisms are in fact monoid (resp.\ group) isomorphisms.

\begin{lemma}\label{lemma LF-w-t} Let \(\mathcal{L}\) be a linear (resp.\ additive) assignment on a category \(\mathbb{X}\) with finite products. Then:
    \begin{enumerate}[{\em (i)}]
        \item \label{Lf-morph} For every morphism \(f\colon X \to Y\), \(\mathcal{L}(f)\colon \mathsf{L}(X) \to \mathsf{L}(Y)\) is a monoid (resp.\ group) morphism;
        \item \label{w-morph} For every pair of objects \(X\) and \(Y\), \(\omega_{X,Y}\colon \mathsf{L}(X \times Y) \to \mathsf{L}(X) \times \mathsf{L}(Y)\) is a monoid (resp.\ group) isomorphism;
        \item \label{t-morph} For every object \(X\), \(t_{\mathcal{L}(\ast)}\colon \mathsf{L}(X) \to \star\) is a monoid (resp.\ group) morphism. Moreover, for the terminal object \(\ast\), \(t_{\mathcal{L}(\ast)}\colon \mathsf{L}(\ast) \to \star\) is a monoid isomorphism and the following equalities hold (where the third equality holds in the case of additive assignment):
              \begin{align}\label{ref t-1-ast}
                  +_\ast = t^{-1}_{\mathcal{L}(\ast)} \circ t_{\mathcal{L}(\ast) \times \mathcal{L}(\ast)} &  & 0_\ast = t^{-1}_{\mathcal{L}(\ast)} &  & -_\ast = 1_{\mathcal{L}(\ast)}
              \end{align}
    \end{enumerate}
\end{lemma}
\begin{proof} For (\ref{Lf-morph}), writing out the natural equalities for \(+\) and \(0\) (and \(-\)) explicitly gives us:
    \begin{align*}\label{eq nat-+-0}
        \mathcal{L}(f) \circ +_X = +_Y \circ \left(\mathcal{L}(f) \times \mathcal{L}(f) \right) &  & \mathcal{L}(f) \circ 0_X = 0_Y &  & \mathcal{L}(f) \circ -_X = -_Y \circ \mathcal{L}(f)
    \end{align*}
    which says precisely that \(\mathcal{L}(f)\) is a monoid morphism (and thus a group morphism for an additive assignment). For (\ref{w-morph}), first observe that it is easy to check that the following equality holds:
    \begin{align*}
        \tau_{\mathcal{L}(X),\mathcal{L}(Y),\mathcal{L}(X),\mathcal{L}(Y)} \circ (\omega_{X,Y} \times \omega_{X,Y}) = \left \langle \mathcal{L}(\pi_1) \times \mathcal{L}(\pi_1), \mathcal{L}(\pi_2) \times \mathcal{L}(\pi_2) \right \rangle
    \end{align*}
    Then using this identity and (\ref{Lf-morph}), we compute that
    \begin{align*}
        \omega_{X,Y} \circ +_{X \times Y} & =~ \langle \mathcal{L}(\pi_1), \mathcal{L}(\pi_2) \rangle \circ +_{X \times Y}                                                                                               \\
                                          & =~ \langle \mathcal{L}(\pi_1) \circ +_{X \times Y}, \mathcal{L}(\pi_2) \circ +_{X \times Y} \rangle                                                                          \\
                                          & =~ \left \langle +_X \circ \left( \mathcal{L}(\pi_1) \times \mathcal{L}(\pi_1) \right), +_Y \circ \left( \mathcal{L}(\pi_1) \times \mathcal{L}(\pi_1) \right) \right \rangle \\
                                          & =~ (+_X \times +_Y) \circ \left \langle \mathcal{L}(\pi_1) \times \mathcal{L}(\pi_1), \mathcal{L}(\pi_2) \times \mathcal{L}(\pi_2) \right \rangle                            \\
                                          & =~ (+_X \times +_Y) \circ \tau_{\mathcal{L}(X),\mathcal{L}(Y),\mathcal{L}(X),\mathcal{L}(Y)}  \circ (\omega_{X,Y} \times \omega_{X,Y})
    \end{align*}
    and
    \begin{gather*}
        \omega_{X,Y} \circ 0_{X \times Y} = \langle \mathcal{L}(\pi_1), \mathcal{L}(\pi_2) \rangle \circ 0_{X \times Y} = \langle \mathcal{L}(\pi_1) \circ 0_{X \times Y}, \mathcal{L}(\pi_2) \circ 0_{X \times Y} \rangle = \langle 0_X, 0_Y \rangle
    \end{gather*}
    So we conclude that \(\omega_{X,Y}\) is a monoid morphism, and since it is also an isomorphism, then it is a monoid isomorphism (and thus, a group isomorphism for an additive assignment). For (\ref{t-morph}), it is automatic from the universal property of the terminal object that for each object \(X\), \(t_{\mathcal{L}(X)}\) is a monoid morphism. In particular, \(t_{\mathcal{L}(\ast)}\) is a monoid isomorphism (and thus, a group isomorphism for an additive assignment). Now, from the universal property of the terminal object, we have that \(t_{\mathcal{L}(\ast)} \circ +_\ast =  t_{\mathcal{L}(\ast) \times \mathcal{L}(\ast)}\) and \(t_{\mathcal{L}(\ast)} \circ 0_\ast = 1_\ast\) (and \(t_{\mathcal{L}(\ast)} \circ -_\ast = t_{\mathcal{L}(\ast)}\)). Post-composing each side of these equalities by \(t^{-1}_{\mathcal{L}(\ast)}\) gives us the desired equalities.
\end{proof}

Recall that any finite product preserving functor sends commutative monoids (resp.\ abelian groups) to commutative monoids (resp.\ abelian groups). Explicitly, if \(\mathcal{F}\colon \mathbb{X} \to \mathbb{Y}\) is a finite product preserving functor and \(\mathsf{M} = (M, \bullet, e)\) is a commutative monoid in \(\X\), then we obtain a commutative monoid \(\mathcal{F}(\mathsf{M})\) in \(\mathbb{Y}\) defined as follows:
\begin{align*}
    \mathcal{F}(\mathsf{M}) = \left( \mathcal{F}(M), \mathcal{F}(\bullet) \circ \omega^{-1}_{\mathcal{F}(M),\mathcal{F}(M)}, \mathcal{F}(e) \circ t^{-1}_{\mathcal{F}(\ast)} \right)
\end{align*}
Similarly, if  \(\mathsf{G} = (G, \bullet, e, \iota)\) is an abelian group in \(\X\), then we obtain an abelian group \(\mathcal{F}(\mathsf{G})\) in \(\mathbb{Y}\) defined as follows:
\begin{align*}
    \mathcal{F}(\mathsf{G}) = \left( \mathcal{F}(M), \mathcal{F}(\bullet) \circ \omega^{-1}_{\mathcal{F}(M),\mathcal{F}(M)}, \mathcal{F}(e) \circ t^{-1}_{\mathcal{F}(\ast)}, \mathcal{F}(\iota) \right)
\end{align*}
Given a linear (resp.\ additive) assignment \(\mathcal{L}\colon \mathbb{X} \to \mathbb{X}\), and an object \(X\), we obtain two commutative monoids (resp.\ abelian groups) \(\mathsf{L}(\mathcal{L}(X))\) and \(\mathcal{L}(\mathsf{L}(X))\) with the same underlying object \(\mathcal{L}\mathcal{L}(X)\). However, it turns out that these commutative monoids (resp.\ abelian groups) are in fact equal on the nose:

\begin{lemma}\label{lemma LL} Let \(\mathcal{L}\) be a linear (resp.\ additive) assignment on a category \(\mathbb{X}\) with finite products. Then, for every object \(X\), \(\mathsf{L}(\mathcal{L}(X)) = \mathcal{L}(\mathsf{L}(X))\). Therefore, the following diagrams commute (and for an additive assignment, the equality below also holds):
    \begin{equation}\begin{gathered}\label{diag L+,L0}
            \xymatrixcolsep{1pc}\xymatrix{ \mathcal{L} \left( \mathcal{L}(X) \times \mathcal{L}(X) \right) \ar[dr]_-{\omega_{\mathcal{L}(X), \mathcal{L}(X)}~~~} \ar[rr]^-{\mathcal{L}(+_X)} &                                        & \mathcal{L}\mathcal{L}(X) & \mathcal{L}(\ast) \ar[rr]^-{\mathcal{L}(0_X)} \ar[dr]_-{t_{\mathcal{L}(\ast)}} &                            & \mathcal{L}\mathcal{L}(X) \\
            & \mathcal{L}\mathcal{L}(X) \times \mathcal{L}\mathcal{L}(X) \ar[ur]_-{~~~+_{\mathcal{L}(X)}} &       &                                                  & \ast \ar[ur]_-{0_{\mathcal{L}(X)}} }\\
            -_{\mathcal{L}(X)} = \mathcal{L}(-_X)
        \end{gathered}\end{equation}
\end{lemma}
\begin{proof} Starting with a linear assignment, on the one hand:
    \[\mathsf{L}(\mathcal{L}(X)) = (\mathcal{L}\mathcal{L}(X), +_{\mathcal{L}(X)}, 0_{\mathcal{L}(X)})\]
    while on the other hand:
    \[\mathcal{L}(\mathsf{L}(X))= \left( \mathcal{L}\mathcal{L}(X), \mathcal{L}(+_X) \circ \omega^{-1}_{\mathcal{L}(X),\mathcal{L}(X)}, \mathcal{L}(0_X) \circ t^{-1}_{\mathcal{L}(\ast)} \right)\]
    Now let us begin by showing that \(0_{\mathcal{L}(X)}\) is equal to \(\mathcal{L}(0_X) \circ t^{-1}_{\mathcal{L}(\ast)}\). To do so, we use naturality of \(0\) and (\ref{ref t-1-ast}):
    \begin{gather*}
        \mathcal{L}(0_X) \circ t^{-1}_{\mathcal{L}(\ast)} =  \mathcal{L}(0_X) \circ 0_\ast = 0_{\mathcal{L}(X)}
    \end{gather*}
    So \(0_{\mathcal{L}(X)} = \mathcal{L}(0_X) \circ t^{-1}_{\mathcal{L}(\ast)}\), which also means that \(0_{\mathcal{L}(X)} \circ t_{\mathcal{L}(\ast)} = \mathcal{L}(0_X)\). To show that \(+_{\mathcal{L}(X)}\) and \(\mathcal{L}(+_X) \circ \omega^{-1}_{\mathcal{L}(X),\mathcal{L}(X)}\) are equal, we apply a version of the Eckmann--Hilton argument. We first compute the following, using naturality of \(+\) and the fact that \(\omega^{-1}\) is a monoid morphism:
    \begin{align*}
         & +_{\mathcal{L}(X)} \circ \left( \mathcal{L}(+_X) \times \mathcal{L}(+_X) \right) \circ \left( \omega^{-1}_{\mathcal{L}(X),\mathcal{L}(X)} \times  \omega^{-1}_{\mathcal{L}(X),\mathcal{L}(X)} \right)                                             \\
         & =~ \mathcal{L}(+_X) \circ +_{\mathcal{L}(X) \times \mathcal{L}(X)} \circ \left( \omega^{-1}_{\mathcal{L}(X),\mathcal{L}(X)} \times  \omega^{-1}_{\mathcal{L}(X),\mathcal{L}(X)} \right)                                                           \\
         & =~ \mathcal{L}(+_X) \circ   \omega^{-1}_{\mathcal{L}(X),\mathcal{L}(X)} \circ (+_{\mathcal{L}(X)} \times +_{\mathcal{L}(X)}) \circ \tau_{\mathcal{L}\mathcal{L}(X),\mathcal{L}\mathcal{L}(X),\mathcal{L}\mathcal{L}(X),\mathcal{L}\mathcal{L}(X)}
    \end{align*}
    So we have that:
    \begin{multline*}
        +_{\mathcal{L}(X)} \circ \left( \mathcal{L}(+_X) \times \mathcal{L}(+_X) \right) \circ \left( \omega^{-1}_{\mathcal{L}(X),\mathcal{L}(X)} \times  \omega^{-1}_{\mathcal{L}(X),\mathcal{L}(X)} \right) \\ = \mathcal{L}(+_X) \circ   \omega^{-1}_{\mathcal{L}(X),\mathcal{L}(X)} \circ (+_{\mathcal{L}(X)} \times +_{\mathcal{L}(X)}) \circ \tau_{\mathcal{L}\mathcal{L}(X),\mathcal{L}\mathcal{L}(X),\mathcal{L}\mathcal{L}(X),\mathcal{L}\mathcal{L}(X)}
    \end{multline*}
    Then, by pre-composing both sides of the above equality by:
    \begin{multline*}\left\langle 1_{\mathcal{L}\mathcal{L}(X)}, \mathcal{L}(0_X) \circ t^{-1}_{\mathcal{L}(\ast)} \circ t_{\mathcal{L}\mathcal{L}(X)} \right\rangle \times \left\langle  \mathcal{L}(0_X) \circ t^{-1}_{\mathcal{L}(\ast)} \circ t_{\mathcal{L}\mathcal{L}(X)}, 1_{\mathcal{L}\mathcal{L}(X)} \right\rangle \\
        = \left\langle 1_{\mathcal{L}\mathcal{L}(X)},  0_{\mathcal{L}(X)} \circ t_{\mathcal{L}\mathcal{L}(X)} \right\rangle \times \left\langle  0_{\mathcal{L}(X)} \circ t_{\mathcal{L}\mathcal{L}(X)}, 1_{\mathcal{L}\mathcal{L}(X)} \right\rangle,
    \end{multline*}
    we get \(+_{\mathcal{L}(X)} = \mathcal{L}(+_X) \circ \omega^{-1}_{\mathcal{L}(X),\mathcal{L}(X)}\), and so, \(+_{\mathcal{L}(X)} \circ \omega_{\mathcal{L}(X),\mathcal{L}(X)} = \mathcal{L}(+_X)\), as desired. Lastly, if \(\mathcal{L}\) is in fact an additive assignment, since we have shown that the underlying commutative monoids of \(\mathsf{L}(\mathcal{L}(X))\) and \(\mathcal{L}(\mathsf{L}(X))\) are equal, then by uniqueness of inverses for commutative monoids, it follows that \(-_{\mathcal{L}(X)} = \mathcal{L}(-_X)\), as desired.
\end{proof}

Clearly, every linear (resp.\ additive) assignment induces a functor from the base category to its category of commutative monoids (resp.\ abelian groups). In fact, it is possible to give an equivalent description of linear assignments as such, which we call a \emph{linear projector}:

\begin{definition}\label{def linear-proj} A \textbf{linear} (resp.\ \textbf{additive}) \textbf{projector} on a category \(\mathbb{X}\) with finite products is a pair \((\mathsf{L}, \nu)\) consisting of a finite product preserving functor \(\mathsf{L} \colon \mathbb{X} \to \CMON[\X]\) (resp.\ \(\mathsf{L} \colon \mathbb{X} \to \mathsf{AbG}[\mathbb{X}]\)) and a natural isomorphism \({\nu_X \colon \mathsf{L}\left(\mathcal{U} (\mathsf{L}(X)) \right) \to \mathsf{L}(X)}\) such that \(\nu_{\mathsf{L}\left(\mathcal{U}(X) \right)} = \mathsf{L}\left(\mathcal{U}(\nu_X) \right)\). As a shorthand, when there is no confusion, we will denote linear assignments simply by their underlying functor \(\mathsf{L}\).
\end{definition}

\begin{theorem}\label{thm ass=proj} For a category with finite products \(\mathbb{X}\), there is a bijective correspondence between linear (resp.\ additive) assignments and linear (resp.\ additive) projectors.
\end{theorem}
\begin{proof} Let us first build a linear projector from a linear assignment. Let \(\mathcal{L}\) be a linear assignment on \(\mathbb{X}\). Define a functor \(\mathsf{L} \colon \mathbb{X} \to \CMON[\X]\) on objects by \(\mathsf{L}(X) = (\mathcal{L}(X), +_X, 0_X)\), which is well-defined by definition, and on morphisms by \(\mathsf{L}(f) = \mathcal{L}(f)\), which is well defined by Lemma \ref{lemma LF-w-t}.(\ref{Lf-morph}). Since \(\mathcal{L}\) is functorial, so is \(\mathsf{L}\). Moreover, since \(\mathcal{L}\) preserves finite products and by Lemma \ref{lemma LF-w-t}.(\ref{w-morph}) and (\ref{t-morph}), it follows that \(\mathsf{L}\) also preserves finite products. Next, observe that \(\mathcal{U} \circ \mathsf{L} = \mathcal{L}\), and thus \(\mathsf{L}\left( \mathcal{U} \left( \mathsf{L}(X) \right) \right) = \mathsf{L}(\mathcal{L}(X))\). Then, by definition, we have that \(\nu_X\colon \mathsf{L}\left( \mathcal{U} \left( \mathsf{L}(X) \right) \right) \to \mathsf{L}(X)\) is a monoid isomorphism, and therefore also gives us our desired natural isomorphism. Lastly, we may express \(\nu_{\mathcal{L}(X)} = \mathcal{L}(\nu_X)\) as \(\nu_{\mathsf{L}\left(\mathcal{U}(X) \right)} = \mathsf{L}\left(\mathcal{U}(\nu_X) \right)\). Thus \(\mathsf{L}\) is a linear projector.

    On the other hand, let \(\mathcal{L}\colon \mathbb{X} \to \CMON[\X]\) be a linear projector on \(\mathbb{X}\). Define an endofunctor on \(\mathbb{X}\) by post-composing \(\mathsf{L}\) with the forgetful functor, \(\mathcal{L} \coloneq \mathcal{U} \circ \mathsf{L}\colon \mathbb{X} \to \mathbb{X}\). Since both \(\mathcal{U}\) and \(\mathsf{L}\) preserve finite products, their composite \(\mathcal{L}\) also preserve finite products. By definition of the forgetful functor, for every object \(X\), the underlying object of the commutative monoid \(\mathsf{L}(X)\) is \(\mathcal{L}(X)\). Thus, we can denote this commutative monoid by \(\mathsf{L}(X) = \left( \mathcal{L}(X), +_X, 0_X \right)\). On morphisms, \(\mathsf{L}(f) = \mathcal{L}(f)\), which are monoid morphisms, so \(+_X\colon \mathcal{L}(X) \times \mathcal{L}(X) \to \mathcal{L}(X)\) and \(0_X\colon \ast \to \mathcal{L}(X)\) are natural transformations. Next, observe that \(\mathsf{L}\left( \mathcal{U} \left( \mathsf{L}(X) \right) \right)  = \mathsf{L}(\mathcal{L}(X))\), and, by definition, we have a natural isomorphism \(\nu_X\colon \mathsf{L}\left( \mathcal{L}(X) \right) \to \mathsf{L}(X)\), which is a monoid morphism. Then, abusing notation slightly, our desired natural isomorphism is \(\nu_X\colon \mathcal{L}\mathcal{L}(X) \to \mathcal{L}(X)\). Lastly, we may express \(\nu_{\mathsf{L}\left(\mathcal{U}(X) \right)} = \mathsf{L}\left(\mathcal{U}(\nu_X) \right)\) as \(\nu_{\mathcal{L}(X)} = \mathcal{L}(\nu_X)\). So we conclude that \(\mathcal{L}\) is a linear assignment.

    It is straightforward to check that these constructions are indeed inverses of each other, so linear assignments are in bijective correspondence with linear projectors. Moreover, this correspondence restricts to a bijective correspondence between additive assignments and additive projectors.
\end{proof}

One can define the notion of morphisms between linear (resp.\ additive) assignments or projectors, and Theorem \ref{thm ass=proj} can be upgraded to an equivalence of categories.

We conclude this section with some preliminary examples of linear assignments and additive assignments, including our main motivating example: the abelianization of groups. Many more examples can be found in Section \ref{sec abelianization}.

\begin{example}\label{ex lin-trival} Every category \(\mathbb{X}\) with finite products admits a trivial additive assignment \(\mathcal{L}_\ast\colon \mathbb{X} \to \mathbb{X}\) which sends every object \(X\) to the terminal object, \(\mathcal{L}_\ast(X) = \ast\). We call \(\mathcal{L}_\ast\) the terminal additive assignment.
\end{example}

\begin{example}\label{ex lin-biprod} Recall that in a \emph{semi-additive category}, by which we mean a category \(\mathbb{X}\) with finite biproducts, every object admits a unique commutative monoid structure, which implies that \(\CMON[\X] \cong \X\). The identity functor \(1_\mathbb{X}\colon \X \to \X\) is then a linear assignment. In fact, the identity functor for any category with finite products is a linear assignment if and only if said category is semi-additive. Similarly, in an \emph{additive category}, by which we mean a semi-additive category \(\mathbb{X}\) which is also enriched over abelian groups, every object admits a unique abelian group structure. In this case, we have that \(\mathsf{AbG}[\mathbb{X}] \cong \X\), and so the identity functor \(1_\mathbb{X}\colon \X \to \X\) is an additive assignment. Again, the identity functor for a category with finite products is a linear assignment if and only if said category is additive. It is worth mentioning that, for a semi-additive category that is \emph{not} an additive category (such as the category of commutative monoids in the classical sense), the identity functor is then an example of a linear assignment that is \emph{not} additive.
\end{example}

\begin{example}\label{ex abelianization is a projector}
    Let \(\Grp\) be the category of groups and \(\Ab\) the category of abelian groups, in the classical sense. The classical Eckmann--Hilton argument implies that commutative monoids and abelian groups objects in \(\Grp\) correspond precisely to abelian groups, thus we have isomorphisms of categories \(\CMON[\Grp] \cong \AbG[\Grp] \cong \Ab\). Then, the \emph{abelianization} functor \(\ab\colon\Grp\to\Ab\), which sends a group \(G\) to its abelianization \(\ab(G) \coloneq G/[G,G]\), can be seen as an additive projector on \(\Grp\). In particular, the required natural isomorphism \(\nu_G\colon \ab\ab(G) \to \ab(G)\) is the obvious isomorphism given by the fact that the abelianization of an abelian group is isomorphic to the starting abelian group. Thus, post-composing the abelianization functor with the forgetful functor from abelian groups to groups gives us the additive assignment \(\mathcal{L}_\ab\colon \Grp \to \Grp\), with \(\mathcal{L}_\ab(G) = \ab(G)\).  We generalize this example in Section \ref{sec abelianization} to the setting of regular (strongly) unital categories.
\end{example}

\section{From Linear Assignments to Tangent Structure}

In this section, we show that a linear (resp.\ additive) assignment on a category \(\X\) induces a (Rosický) tangent structure on \(\X\). In order to keep this paper as self-contained as possible, we review the full definition of a tangent category. For an in-depth introduction to tangent categories, we refer the reader to \cite{cockett2014differential,cockettCruttwellDiffBundles,cruttwell2023differential,RosickyTangentCats}.

We first recall the definition of additive (resp.\ abelian group) bundles \cite[Definition 2.1]{cockett2014differential}, a central concept in the definition of a (Rosický) tangent structure. These are essentially commutative monoids (resp.\ abelian groups) in slice categories. For slice categories to admit all finite products, the base category must admit finite pullbacks. The definition of additive bundles, however, does not require all finite pullbacks, but only the pullbacks of copies of the same morphism.

In a category \(\mathbb{X}\), a \textbf{bundle} is a morphism \(q\colon E \to X\) such that, for all \(n \in \mathbb{N}\), the pullback of \(n\) copies of \(q\) exists. We denote this pullback by \(E_n\), and its \(n\) projections by \(\rho_j\colon E_n \to E\), for all \(1\leq j\leq n\), so that \(q \circ \rho_j = q \circ \rho_i\) for all \(1 \leq i,j \leq n\). By convention, \(E_0=X\) and \(E_1= E\). Then, an \textbf{additive bundle}\footnote{We appreciate that there is a bit of an unfortunate clash of terminology. Here, additive bundle means a bundle without negatives, while we use the term `additive assignment' to mean a linear assignment with negatives.} is a quintuple \(\mathsf{E}= (q\colon E \to X, \bullet, e)\) consisting of a bundle \((X, E, q)\), equipped with morphisms \(\bullet\colon E_2 \to E\) and \(e\colon X \to E\), such that the following diagrams commute:
\begin{equation*}\label{diag bullet-e-slice}\begin{gathered}
        \xymatrixcolsep{5pc}\xymatrix{ E_2 \ar[d]_-{\rho_j}  \ar[r]^-{\bullet} & E \ar[d]^-{q} & A \ar[r]^-{e}  \ar@{=}[dr]^-{} & E  \ar[d]^-{q}  \\
        E \ar[r]_-{q} & X & & X   }  \end{gathered}\end{equation*}
\begin{equation*}\label{diag addbun}\begin{gathered}  \xymatrixcolsep{4pc}\xymatrix{ E_3 \ar[r]^-{\left \langle \bullet \circ \langle \rho_1, \rho_2 \rangle , \rho_3 \right \rangle} \ar[d]_-{\left \langle \rho_1, \bullet \circ \langle \rho_2, \rho_3 \rangle  \right \rangle} & E_2 \ar[d]^-{\bullet} & E \ar[r]^-{\langle e \circ q, 1_{E} \rangle} \ar[d]_-{\langle 1_{E}, e \circ q \rangle}  \ar@{=}[dr]^-{} & E_2  \ar[d]^-{\bullet}  & E_2 \ar[dr]_-{\bullet} \ar[r]^-{\langle \rho_2, \rho_1 \rangle} & E_2 \ar[d]^-{\bullet}  \\
        E_2 \ar[r]_-{\bullet} & E & E_2  \ar[r]_-{\bullet} & E &   & E . } \end{gathered}\end{equation*}
Here, \(\langle -, - \rangle\) is the pairing operator coming from the universal property of the pullback. If \(\mathsf{E}= (q\colon E \to X, \bullet, e)\) is an additive bundle, then we shall also say that \(\mathsf{E}\) is an additive bundle over \(X\). If \(\mathbb{X}\) admits all finite pullbacks, then for each object \(X\), the slice category over \(X\) admits finite products (given by the pullbacks in \(\mathbb{X}\)). In this setting, additive bundles over \(X\) correspond precisely to the commutative monoids in the slice category over \(X\). Indeed, a morphism of type \(q\colon E \to X\) is an object in the slice category, the top two diagrams say that \(\bullet\) and \(e\) are morphisms in the slice category, while the bottom three are precisely the commutative monoid axioms in the slice category. That said, as mentioned above, we will not assume in general that \(\mathbb{X}\) admits all finite pullbacks\footnote{In fact, many key examples of tangent categories do not have all finite pullbacks, such as for instance the category of smooth manifolds and other differential geometry related categories.}.

Morphisms between additive bundles \cite[Definition 2.3]{cockett2014differential} correspond to certain monoid morphisms, but where we can also change the base object: if \(q\colon E \to X\) and \(q^\prime\colon E^\prime \to X^\prime\) are bundles, then a \textbf{bundle morphism} \((f,g)\colon q \to q^\prime\) is a pair of morphisms \(f\colon E \to E^\prime\) and \(g\colon X \to X^\prime\) such that the following diagram commutes:
\begin{equation*}\label{diag bunmorph}\begin{gathered}
        \xymatrixcolsep{5pc}\xymatrix{ E \ar[d]_-{q}  \ar[r]^-{f} & E^\prime \ar[d]^-{q^\prime}   \\
        X \ar[r]_-{g} & X^\prime }  \end{gathered}\end{equation*}
If \(\mathsf{E} = (q\colon E \to X, \bullet, e)\) and \(\mathsf{E}^\prime = (q^\prime\colon E^\prime \to X^\prime, \bullet^\prime, e^\prime)\) are additive bundles, then an \textbf{additive bundle morphism} \((f,g)\colon \mathsf{E} \to \mathsf{E}^\prime\) is a bundle morphism \((f,g)\colon q \to q^\prime\) such that the following diagrams commute:
\begin{equation*}\label{diag addbunmorph}\begin{gathered}
        \xymatrixcolsep{5pc}\xymatrix{ E_2  \ar[rr]^-{\langle f \circ \rho_1, f \circ \rho_2 \rangle} \ar[d]_-{\bullet} && E^\prime_2 \ar[d]^-{\bullet^\prime} & X \ar[r]^-{e} \ar[d]_-{g} & E \ar[d]^-{f}  \\
        E \ar[rr]_-{g} && E^\prime & X^\prime \ar[r]_-{e} &  E^\prime .}  \end{gathered}\end{equation*}
We denote by \(\mathsf{ABUN}[\mathbb{X}]\) the category of additive bundles and additive bundle morphisms between them.

We can also consider the abelian group analogue of an additive bundle. This was considered originally by Rosický in \cite[Section 1]{RosickyTangentCats}. An \textbf{abelian group bundle}\footnote{Abelian group bundles are also known as Beck modules.} is a sextuple \(\mathsf{E}= (q\colon E \to X, \bullet, e, \iota)\) consisting of an additive bundle \((q\colon E \to X, \bullet, e)\) with a morphism \(\iota\colon E \to E\) such that the following diagrams commute:
\begin{equation*}
    \xymatrixcolsep{5pc}\vcenter{\xymatrix{ E \ar[r]^-{\iota}  \ar[dr]_-{q} & E  \ar[d]^-{q} \\
    & X }}
    \qquad\qquad
    \vcenter{\xymatrix{ & E_2  \ar[dr]^-{\bullet} \\
    E \ar[ur]^-{\langle 1_{E}, \iota \rangle} \ar[dr]_-{\langle \iota, 1_E \rangle}  \ar[r]^-{q} & X  \ar[r]^-{e} & E. \\
    & E_2  \ar[ur]_-{\bullet} }}
\end{equation*}
Once again, if the base category \(\mathbb{X}\) has all finite pullbacks, then abelian group bundles over an object \(X\) correspond precisely to the abelian groups in the slice category over \(X\), where the top diagram says that \(\iota\) is a morphism in the slice category and the bottom diagram is the additional abelian group axiom in the slice category. If \(\mathsf{E}=(q\colon E \to X, \bullet, e, \iota)\) and \(\mathsf{E}^\prime = (q^\prime\colon E^\prime \to X^\prime, \bullet^\prime, e^\prime, \iota^\prime)\) are abelian group bundles, then an \textbf{abelian group morphism} \((f,g)\colon q \to q^\prime\) is an additive bundle morphism between their underlying additive bundles, and furthermore the following diagram automatically commutes:
\begin{equation*}\label{diag adgroupbunmorph}\begin{gathered}
        \xymatrixcolsep{5pc}\xymatrix{ E \ar[d]_-{f}  \ar[r]^-{\iota} & E \ar[d]^-{f}   \\
        E^\prime \ar[r]_-{\iota^\prime} & E^\prime }  \end{gathered}\end{equation*}
We denote by \(\mathsf{GBUN}[\mathbb{X}]\) the category of abelian group bundles and abelian group bundle morphisms between them.

We may now review the definition of tangent categories. A \textbf{tangent structure} \cite[Definition 2.3]{cockett2014differential} on a category \(\mathbb{X}\) is a sextuple \(\mathbb{T} \coloneq (\mathcal{T}, \mathsf{p}, \mathsf{s}, \mathsf{z}, \ell, \mathsf{c})\) consisting of the following data:
\begin{itemize}
    \item An endofunctor \(\mathcal{T}\colon \mathbb{X} \to  \mathbb{X}\), called the \textbf{tangent bundle functor}, which we think of as a functor associating, to each object \(X\), an abstract tangent bundle \(\mathcal{T}(X)\).
    \item A natural transformation \(\mathsf{p}_X\colon \mathcal{T}(X) \to X\), called the \textbf{projection}, which is an analogue of the natural projection from the tangent bundle down to its base space, such that for each object \(X\), \(\mathsf{p}_X\) is a bundle over \(X\). We will denote the pullback of \(n\) copies of \(\mathsf{p}_X\) by \(\mathcal{T}_n(X)\), and this induces a family of endofunctors \(\mathcal{T}_n\colon \mathbb{X} \to \mathbb{X}\). One should interpret \(\mathcal{T}_n(X)\) as the space of \(n\) tuples of tangent vectors anchored over the same point. We also ask that for all \(m\in \mathbb{N}\), \(\mathcal{T}^m\) preserves these pullbacks, that is, \(\mathcal{T}^m\left(\mathcal{T}_n(X) \right)\) is the pullback of \(n\) copies of \(\mathcal{T}(\mathsf{p}_X)\colon \mathcal{T}(X) \to  \mathcal{T}\mathcal{T}(X)\), or in other words, that \(\mathcal{T}(\mathsf{p}_X)\) is also a bundle over \(\mathcal{T}(X)\).
    \item Natural transformations \(\mathsf{s}_X\colon \mathcal{T}_2(X) \to \mathcal{T}(X)\), called the \textbf{sum}, and \(\mathsf{z}_X\colon A \to \mathcal{T}(X)\), called the \textbf{zero}, such that for each object \(X\), \(\mathsf{T}(X) = (\mathsf{p}_X\colon \mathcal{T}(X) \to X , \mathsf{s}_X, \mathsf{z}_X)\) is an additive bundle. For smooth manifolds, the sum captures the ability of adding two tangent vectors over the same base point, while the zero picks out the zero tangent vector over a given point. This allows us to view \(\mathcal{T}(X)\) as a kind of smooth vector bundle over~\(X\), where each fibre is a commutative monoid. Note that, since \(\mathcal{T}\) preserves the pullbacks of \(\mathsf{p}_X\), we also have that \(\mathcal{T}(\mathsf{T}(X)) = (\mathcal{T}(\mathsf{p}_X)\colon \mathcal{T}\mathcal{T}(X) \to \mathcal{T}(X), \mathcal{T}(\mathsf{s}_X), \mathcal{T}(\mathsf{z}_X))\) is an additive bundle over \(\mathcal T(X)\).
    \item A natural transformation \(\ell_X\colon \mathcal{T}(X) \to \mathcal{T}\mathcal{T}(X)\), called the \textbf{vertical lift}, which essentially encodes linearity of differentiation. We ask for \((\ell_X, \mathsf{z}_X)\colon \mathsf{T}(X) \to \mathcal{T}(\mathsf{T}(X))\) to be an additive bundle morphism and that the following diagram commutes:
          \begin{equation}\label{diag el1}\begin{gathered}
                  \xymatrixcolsep{2.5pc}\xymatrix{ \mathcal{T}(X) \ar[r]^-{\ell_X} \ar[d]_-{\ell_X}  & \mathcal{T}\mathcal{T}(X) \ar[d]^-{\ell_{\mathcal{T}(X)}} \\
                  \mathcal{T}\mathcal{T}(X) \ar[r]_-{\mathcal{T}(\ell_X)} &  \mathcal{T}\mathcal{T}\mathcal{T}(X)
                  }  \end{gathered}\end{equation}
          However the key feature of the vertical lift is that it is universal, in the sense that the following diagram is a pullback diagram:
          \begin{equation}\label{universalift}\begin{gathered}
                  \xymatrixcolsep{5pc}\xymatrix{ \mathcal{T}_2(X) \ar[d]_-{\mathsf{p}_X \circ \rho_j}  \ar[r]^-{\langle \ell \circ \rho_1, \mathsf{z}_{\mathcal{T}(X)} \circ \rho_2 \rangle} & \mathcal{T}\mathcal{T}_2(X) \ar[r]^-{\mathcal{T}(\mathsf{s}_X)} & \mathcal{T}\mathcal{T}(X) \ar[d]^-{\mathcal{T}(\mathsf{p}_X)}  \\
                  X \ar[rr]_-{\mathsf{z}_X} & & \mathcal{T}(X).
                  }  \end{gathered}\end{equation}
          It comes for free that all powers of the tangent bundle functor \(\mathcal{T}^m\) preserve this pullback \cite[Lemma 2.15]{cockett2014differential}. This universal property of vertical lift means that the tangent bundle embeds into the double tangent bundle via the vertical bundle, and is essential for formalizing important properties of the tangent bundle from differential geometry, see \cite[Section 2.5]{cockett2014differential} for more details.
    \item A natural isomorphism \(\mathsf{c}_X\colon \mathcal{T}\mathcal{T}(X) \to \mathcal{T}\mathcal{T}(X)\), called the \textbf{canonical flip}, which is an analogue of the smooth involution of the same name on the double tangent bundle of smooth manifolds. This transformation encodes symmetry of the mixed partial derivatives. In particular, we ask for \(\mathsf{c}_X\) to be its own inverse, and that it satisfies the following Yang--Baxter identity:
          \begin{equation}\label{cidentitties2}
              \xymatrixcolsep{2.5pc}\vcenter{\xymatrix{ \mathcal{T}\mathcal{T}(X) \ar[r]^-{\mathsf{c}_X} \ar@{=}[dr]_-{}  &  \mathcal{T}\mathcal{T}(X) \ar[d]^-{\mathsf{c}_X} & \mathcal{T}\mathcal{T}\mathcal{T}(X) \ar[d]_-{\mathsf{c}_{\mathcal{T}(X)}} \ar[r]^-{\mathcal{T}(\mathsf{c}_X)} &  \mathcal{T}\mathcal{T}\mathcal{T}(X) \ar[r]^-{\mathsf{c}_{\mathcal{T}(X)}} & \mathcal{T}\mathcal{T}\mathcal{T}(X) \ar[d]^-{\mathcal{T}(\mathsf{c}_X)}  \\
              & \mathcal{T}\mathcal{T}(X) & \mathcal{T}\mathcal{T}\mathcal{T}(X) \ar[r]_-{\mathcal{T}(\mathsf{c}_X)} & \mathcal{T}\mathcal{T}\mathcal{T}(X) \ar[r]_-{\mathsf{c}_{\mathcal{T}(X)}} &  \mathcal{T}\mathcal{T}\mathcal{T}(X)}}
          \end{equation}
          We also require for the morphism \((\mathsf{c}_X, 1_{\mathcal{T}(X)})\colon \mathcal{T}(\mathsf{T}(X)) \to \mathsf{T}(\mathcal{T}(X))\) to be an additive bundle (iso)morphism, and to be compatible with the vertical lift, in the sense that the following diagrams commute:
          \begin{equation}\label{ellcidentitties}\begin{gathered}
                  \xymatrixcolsep{2.5pc}\xymatrix{ \mathcal{T}(X) \ar[r]^-{\ell_X} \ar[dr]_-{\ell_X} & \mathcal{T}\mathcal{T}(X)  \ar[d]^-{\mathsf{c}_X} & \mathcal{T}\mathcal{T}(X)\ar[r]^-{\ell_{\mathcal{T}(X)}} \ar[d]_-{\mathsf{c}_X} &  \mathcal{T}\mathcal{T}\mathcal{T}(X)  \ar[r]^-{\mathcal{T}(\mathsf{c}_X)} & \mathcal{T}\mathcal{T}\mathcal{T}(X) \ar[d]^-{\mathsf{c}_{\mathcal{T}(X)}} \\
                  &  \mathcal{T}\mathcal{T}(X) & \mathcal{T}\mathcal{T}(X) \ar[rr]_-{\mathcal{T}(\ell_X)}  && \mathcal{T}\mathcal{T}\mathcal{T}(X).
                  }  \end{gathered}\end{equation}
\end{itemize}

Rosický's original definition \cite[Section 2]{RosickyTangentCats} required for the fibres of the tangent bundle to be abelian groups rather than commutative monoids. So, a \textbf{Rosický tangent structure} (also sometimes called a \textbf{tangent structure with negatives}) \cite[Section 3.3]{cockett2014differential} on a category \(\mathbb{X}\) is a septuple \({\mathbb{T} \coloneq (\mathcal{T}, \mathsf{p}, \mathsf{s}, \mathsf{z}, \ell, \mathsf{c}, \mathsf{n})}\) consisting of a tangent structure \((\mathcal{T}, \mathsf{p}, \mathsf{s}, \mathsf{z}, \ell, \mathsf{c})\) on \(\mathbb{X}\) with an extra natural transformation \({\mathsf{n}_X\colon \mathcal{T}(X) \to \mathcal{T}(X)}\), called the \textbf{negative}, such that \(\mathsf{T}(X) = (\mathsf{p}_X, \mathsf{s}_X, \mathsf{z}_X, \mathsf{n}_X)\) is an abelian group bundle. Intuitively, in this setting, this allows us to take the negation of tangent vectors. The vertical lift and the canonical flip are automatically abelian group bundle morphisms.

A \textbf{(Rosický) tangent category} is then a pair \((\mathbb{X}, \mathbb{T})\) consisting of a category \(\mathbb{X}\) equipped with a (Rosický) tangent structure \(\mathbb{T}\) on \(\X\). If a tangent category \((\mathbb{X}, \mathbb{T})\) has finite products, and if the tangent bundle functor \(\mathcal T\) preserve those finite products, it is called a \textbf{cartesian (Rosický) tangent category} \cite[Definition 2.8]{cockett2014differential}. Here are now some basic examples of tangent categories. More examples can be found in \cite[Example 2.2]{cockettCruttwellDiffBundles}.

\begin{example}\label{ex manifold-tancat} The archetype of a tangent category is the category of smooth manifolds, where the tangent structure is given by the classical tangent bundle. So let \(\mathsf{SMAN}\) be the category whose objects are (finite-dimensional real) smooth manifolds, and whose morphisms are smooth functions between them. Then, \(\mathsf{SMAN}\) is a Cartesian Rosický tangent category, whose tangent bundle functor is the usual tangent bundle functor which sends a smooth manifold \(M\) to its tangent bundle \(\mathcal{T}(M)\). Recall that, in local coordinates, elements of the tangent bundle \(\mathcal{T}(M)\) can be described as pairs \((x,\vec v)\) of a point \(x \in M\) and a tangent vector \(\vec v\) at \(x\). Then the remaining tangent structure is defined in local coordinates as follows:
    \begin{gather*}
        \mathsf{p}_M(x, \vec v) = x \qquad \mathsf{s}_M(x,\vec v, \vec w) = (x, \vec v + \vec w) \qquad \mathsf{z}_M(x) = (x,\vec 0) \qquad \mathsf{n}_M(x,\vec v) = (x, -\vec v) \\
        \ell_M(x,\vec v) = ( x, \vec0, \vec 0, \vec v ) \qquad  \mathsf{c}_M ( x, \vec v, \vec w, \vec u )  = (  x, \vec w, \vec v, \vec u ) .
    \end{gather*}
\end{example}

\begin{example}\label{ex trivial-tancat} Trivially, any category \(\mathbb{X}\) (with finite products) is a (cartesian) Rosický tangent category whose tangent bundle functor is simply the identity functor \(1_{\mathbb{X}}\), with all the structural natural transformations being identity morphisms.
\end{example}

\begin{example}\label{ex biprod-tancat} Every semi-additive (resp.\ additive) category \(\mathbb{X}\) is a cartesian (Rosický) tangent category whose tangent bundle functor is the diagonal functor \(\Delta_\mathbb{X}\colon \X \to \X\), which is defined by \(\Delta(X) = X \times X\) on objects, and similarly, \(\Delta(f)=f\times f\), on morphisms.
\end{example}

We now show that commutative monoid objects in a category can be used to build certain additive bundles. More precisely, we show that the product of a base object with a commutative monoid object is an additive bundle over said base object. This construction will play a key role in the tangent structure induced by a linear assignment.

\begin{lemma}\label{gluing} Let \(\mathbb{X}\) be a category with finite products.
    \begin{enumerate}[{\em (i)}]
        \item\label{gluing.proj} For every pair of objects \(X\) and \(Y\), the projection \(\pi_1\colon X \times Y \to X\) is a bundle over~\(X\), the pullback of \(n\) copies of \(\pi_1\) is \((X \times Y)_n = X \times Y^n\) (where \(Y^n\) is the \(n\)-ary product of \(n\) copies of \(Y\)), and in particular, \((X \times Y)_2 = X \times (Y \times Y)\). Moreover, for every pair of morphisms \(f\colon Y \to Y^\prime\) and \(g\colon X \to X^\prime\), \((g \times f,g)\colon \pi_1 \to \pi_1\) is a bundle morphism.
        \item\label{gluing.mon} For every object \(X\) and commutative monoid \(\mathsf{M} = (M, \bullet, e)\), the tuple
              \begin{align*}
                  X \times \mathsf{M}\coloneq (\pi_1: X \times M \to X, 1_X \times \bullet, \langle 1_X, e \circ t_x\rangle )
              \end{align*}
              is an additive bundle over \(X\). Moreover, for every monoid morphism \(f\colon \mathsf{M} \to \mathsf{M}^\prime\) and morphism \(g\colon X \to X^\prime\), \((g \times f,g)\colon X \times \mathsf{M} \to X^\prime \times \mathsf{M}^\prime\) is an additive bundle morphism. This induces a functor \(\mathbb{X} \times \CMON[\X] \to \mathsf{ABUN}[\mathbb{X}]\).
        \item\label{gluing.group} For every object \(X\) and any abelian group \(\mathsf{G} = (G, \bullet, e, \iota)\), the tuple
              \begin{align*}
                  X \times \mathsf{G} \coloneq (\pi_1: X \times G \to X, 1_X \times \bullet, \langle 1_X, e \circ t_x\rangle, 1_X \times \iota)
              \end{align*}
              is an additive bundle. Moreover, for every group morphism \(f\colon \mathsf{G} \to \mathsf{G}^\prime\) and morphism \(g\colon X \to X^\prime\), \((g \times f,g)\colon X \times \mathsf{G} \to X^\prime \times \mathsf{G}^\prime\) is an abelian group bundle morphism. This induces a functor \(\mathbb{X} \times \Ab[\X] \to \mathsf{GBUN}[\mathbb{X}]\).
    \end{enumerate}
    Furthermore, if \(\mathcal{F}\colon \mathbb{X} \to \mathbb{Y}\) preserves finite products, then:
    \begin{enumerate}[{\em (i)}]
        \setcounter{enumi}{3}
        \item For every pair of objects \(X\) and \(Y\) of \(\mathbb{X}\), \(\mathcal{F}(\pi_1)\colon \mathcal{F}(X \times Y) \to \mathcal{F}(X)\) is a bundle over \(\mathcal{F}(X)\), and \(\mathcal{F}(X \times Y)_n = \mathcal{F}(X \times Y^n)\).
        \item \label{gluing.mon.fun} For every object \(X\) and commutative monoid \(\mathsf{M}\) of \(\mathbb{X}\), the tuple
              \begin{align*}
                  \mathcal{F}( X \times \mathsf{M}) \coloneq \left( \mathcal{F}(\pi_1): \mathcal{F}(X \times M) \to \mathcal{F}(X), \mathcal{F}(1_X \times \bullet), \mathcal{F}(\langle 1_X, e \circ t_x\rangle) \right)
              \end{align*}
              is an additive bundle over \(\mathcal{F}(X)\), and \((\omega_{X,M}, 1_{\mathcal{F}(X)})\colon \mathcal{F}( X \times \mathsf{M}) \to \mathcal{F}(X) \times \mathcal{F}(\mathsf{M})\) is an additive bundle isomorphism.
        \item  For every object \(X\) and abelian group \(\mathsf{G}\) of \(\mathbb{X}\), the tuple
              \begin{align*}
                  \mathcal{F}( X \times \mathsf{G}) \coloneq \left( \mathcal{F}(\pi_1): \mathcal{F}(X \times G) \to \mathcal{F}(X), \mathcal{F}(1_X \times \bullet), \mathcal{F}(\langle 1_X, e \circ t_x\rangle), \mathcal{F}(1_X \times \iota) \right)
              \end{align*}
              is an abelian group bundle over \(\mathcal{F}(X)\), and \((\omega_{X,M}, 1_{\mathcal{F}(X)})\colon \mathcal{F}( X \times \mathsf{G}) \to \mathcal{F}(X) \times \mathcal{F}(\mathsf{G})\) is an abelian group bundle isomorphism.
    \end{enumerate}
\end{lemma}
\begin{proof} These statements are straightforward to check, so we leave this as an exercise for the reader.
\end{proof}

We may now state the main result of this paper, which shows that a linear assignment induces a tangent structure. Let \(\mathcal{L}\) be a linear assignment on a category \(\mathbb{X}\) with finite products. We define the tangent structure \(\mathbb{T}_\mathcal{L} = (\mathcal{T}_\mathcal{L},\mathsf{p},\mathsf{s},\mathsf{z},\ell,\mathsf{c})\) as follows:
\begin{itemize}
    \item Define the tangent bundle functor \(\mathcal{T}_\mathcal{L}\colon \X \to \X\) as the product of the identity functor with the linear assignment:
          \begin{align*}
              \mathcal{T}_\mathcal{L}(-) = - \times \mathcal{L}(-)
          \end{align*}
    \item Define the projection \(\mathsf{p}_X\colon \mathcal{T}_\mathcal{L}(X) \to X\) as simply the first projection of the product:
          \begin{align*}
              \mathsf{p}_X \coloneq \pi_1: X \times \mathcal{L}(X) \to X
          \end{align*}
    \item Note that \({\mathcal{T}_\mathcal{L}}_2(X) = X \times \left( \mathcal{L}(X) \times \mathcal{L}(X) \right)\). Then define the sum \(\mathsf{s}_X\colon {\mathcal{T}_\mathcal{L}}_2(X) \to \mathcal{T}_\mathcal{L}(X)\) and zero \(\mathsf{z}_X\colon A \to \mathcal{T}(X)\) respectively as follows:
          \begin{align*}
              &\mathsf{s}_X \coloneq 1_X \times +_X  \colon X \times \left( \mathcal{L}(X) \times \mathcal{L}(X) \right) \to X \times \mathcal{L}(X) \\
              &\mathsf{z}_X \coloneq \langle 1_X, 0_X \circ t_X \rangle \colon X \to X \times \mathcal{L}(X)
          \end{align*}
    \item The vertical lift \(\ell_X\) is, up to isomorphism, essentially given by inserting zeroes into the middle components. More precisely, it is the following composite:
          \begin{equation*}
              \xymatrixcolsep{1.5pc}\xymatrix{ X \times \mathcal{L}(X) \ar[rd]|-{\langle 1_X, 0_X \circ t_X \rangle \times \langle 0_X \circ t_{\mathcal{L}(X)}, \nu^{-1}_X \rangle} \ar[rr]^-{\ell_X} && (X \times \mathcal{L}(X)) \times \mathcal{L}\left(X \times \mathcal{L}(X)\right)\\
              & (X \times \mathcal{L}(X)) \times (\mathcal{L}(X) \times \mathcal{L}\mathcal{L}(X)) \ar[ru]|-{1_{X \times \mathcal{L}(X)} \times \omega^{-1}_{X, \mathcal{L}(X)}}
              }
          \end{equation*}
    \item The canonical flip \(\mathsf{c}_X\colon \mathcal{T}_\mathcal{L}(X) \to \mathcal{T}_\mathcal{L}\mathcal{T}_\mathcal{L}(X)\) is, up to isomorphism, essentially given by the natural interchange isomorphism of the product (which swaps the middle arguments). More precisely, it is the following composite:
          \begin{equation*}\label{eq cf2}\begin{gathered}
                  \xymatrixcolsep{6pc}\xymatrix{ (X \times \mathcal{L}(X)) \times \mathcal{L}\left(X \times \mathcal{L}(X)\right) \ar[d]_{\mathsf{c}_X} \ar[r]^-{1_{X \times \mathcal{L}(X)} \times \omega_{X, \mathcal{L}(X)}} & (X \times \mathcal{L}(X)) \times (\mathcal{L}(X) \times \mathcal{L}\mathcal{L}(X)) \ar[d]|-{\tau_{X,\mathcal{L}(X),\mathcal{L}(X),\mathcal{L}\mathcal{L}(X)}} \\
                  (X \times \mathcal{L}(X)) \times \mathcal{L}\left(X \times \mathcal{L}(X)\right) & (X \times \mathcal{L}(X)) \times (\mathcal{L}(X) \times \mathcal{L}\mathcal{L}(X)) \ar[l]^-{1_{X \times \mathcal{L}(X)} \times \omega^{-1}_{X, \mathcal{L}(X)}}
                  }
              \end{gathered}\end{equation*}
\end{itemize}

If \(\mathcal{L}\) is also an additive assignment, slightly abusing the notation, we define a Rosický tangent structure \(\mathbb{T}_\mathcal{L} = (\mathcal{T}_\mathcal{L},\mathsf{p},\mathsf{s},\mathsf{z},\ell,\mathsf{c}, \mathsf{n})\), where the first six components are defined as above, and where the negative \({\mathsf{n}_X\colon \mathcal{T}(X) \to \mathcal{T}(X)}\) is defined as follows:
\begin{gather*}
    \mathsf{n}_X \coloneq 1_X \times -_X : X \times \mathcal{L}(X) \to X \times \mathcal{L}(X)
\end{gather*}

\begin{theorem}\label{theo:linearprojtotangent} Let \(\mathcal{L}\) be a linear (resp.\ additive) assignment on a category \(\mathbb{X}\) with finite products. Then, \((\mathbb{X}, \mathbb{T}_\mathcal{L})\) as defined above is a cartesian (Rosický) tangent category.
\end{theorem}
\begin{proof} By construction, \(\mathcal{T}_\mathcal{L}\) is indeed a functor, and \(\mathsf{p}\), \(\mathsf{s}\), \(\mathsf{z}\), \(\ell\), and \(\mathsf{c}\) (and \(\mathsf{n}\)) are natural transformations. Since \(\mathcal{L}\) preserves finite products, it follows that \(\mathcal{T}_\mathcal{L}\) preserves finite products as well. By definition of the projection, and by Lemma \ref{gluing}.(\ref{gluing.proj}), \(\mathsf{p}_X\) is a bundle over \(X\) and all powers of \(\mathcal{T}_\mathcal{L}\) preserve the necessary pullbacks. Note that the tuple \(\mathsf{T}_\mathcal{L}(X) = (\mathsf{p}_X, \mathsf{s}_X, \mathsf{z}_X, \mathsf{n}_X)\) (resp.\ \(\mathsf{T}(X) = (\mathsf{p}_X, \mathsf{s}_X, \mathsf{z}_X, \mathsf{n}_X)\)) is defined as in Lemma \ref{gluing}.(\ref{gluing.mon}) (resp.\ \ref{gluing.group}). In other words:
    \begin{equation*}
        \mathsf{T}_\mathcal{L}(X) = X \times \mathsf{L}(X).
    \end{equation*}
    Thus, by Lemma \ref{gluing}.(\ref{gluing.mon}) (resp.\ (\ref{gluing.group})), \(\mathsf{T}(X)\) is indeed an additive (resp.\ abelian group) bundle.

    Let us now consider the canonical flip \(\mathsf{c}_X\). By construction, it is the composite of isomorphisms, and so is itself an isomorphism. Moreover, since the interchange morphism is its own inverse: \(\tau^{-1}_{X,\mathcal{L}(X),\mathcal{L}(X),\mathcal{L}(X)} =  \tau_{X,\mathcal{L}(X),\mathcal{L}(X),\mathcal{L}(X)}\), it follows that \(\mathsf{c}_X\) is also its own inverse, thus satisfying the first diagram of (\ref{cidentitties2}). For the Yang--Baxter diagram, first note that:
    \[ \mathcal{T}^3_\mathcal{L}(X) \cong \left( (X \times \mathcal{L}(X)) \times (\mathcal{L}(X) \times \mathcal{L}^2(X)) \right) \times \left( (\mathcal{L}(X) \times \mathcal{L}^2(X)) \times (\mathcal{L}^2(X) \times \mathcal{L}^3(X)) \right)  \]
    Thus, slightly abusing the notation, we may view \(\mathcal{T}^3_\mathcal{L}(X)\) as an octonary product of these objects with projections:
    \begin{gather*}
        \pi_1: \mathcal{T}^3_\mathcal{L}(X) \to X \quad \pi_2: \mathcal{T}^3_\mathcal{L}(X) \to \mathcal{L}(X) \quad \pi_3: \mathcal{T}^3_\mathcal{L}(X) \to \mathcal{L}(X) \\
        \pi_4: \mathcal{T}^3_\mathcal{L}(X) \to \mathcal{L}^2(X) \quad \pi_5: \mathcal{T}^3_\mathcal{L}(X) \to \mathcal{L}(X) \quad \pi_6: \mathcal{T}^3_\mathcal{L}(X) \to \mathcal{L}^2(X) \\
        \pi_7: \mathcal{T}^3_\mathcal{L}(X) \to \mathcal{L}^2(X) \quad \pi_8: \mathcal{T}^3_\mathcal{L}(X) \to \mathcal{L}^3(X)
    \end{gather*}
    Since the interchange morphism satisfies the following Yang--Baxter identity (omitting indices for readability): \((\tau \times \tau) \circ \tau \circ (\tau \times \tau) = \tau \circ (\tau \times \tau) \circ \tau\), one can easily show that the following equalities hold:
    \begin{align*}
        \pi_1 \circ  \mathcal{T}(\mathsf{c}_X) \circ \mathsf{c}_{\mathcal{T}_\mathcal{L}(X)} \circ  \mathcal{T}(\mathsf{c}_X) = \pi_1 = \pi_1 \circ \mathsf{c}_{\mathcal{T}_\mathcal{L}(X)} \circ  \mathcal{T}(\mathsf{c}_X)\circ  \mathcal{T}(\mathsf{c}_X), \\
        \pi_2 \circ  \mathcal{T}(\mathsf{c}_X) \circ \mathsf{c}_{\mathcal{T}_\mathcal{L}(X)} \circ  \mathcal{T}(\mathsf{c}_X) = \pi_5 = \pi_2 \circ \mathsf{c}_{\mathcal{T}_\mathcal{L}(X)} \circ  \mathcal{T}(\mathsf{c}_X)\circ  \mathcal{T}(\mathsf{c}_X), \\
        \pi_3 \circ  \mathcal{T}(\mathsf{c}_X) \circ \mathsf{c}_{\mathcal{T}_\mathcal{L}(X)} \circ  \mathcal{T}(\mathsf{c}_X) = \pi_3 = \pi_3 \circ \mathsf{c}_{\mathcal{T}_\mathcal{L}(X)} \circ  \mathcal{T}(\mathsf{c}_X)\circ  \mathcal{T}(\mathsf{c}_X), \\
        \pi_4 \circ  \mathcal{T}(\mathsf{c}_X) \circ \mathsf{c}_{\mathcal{T}_\mathcal{L}(X)} \circ  \mathcal{T}(\mathsf{c}_X) = \pi_7 = \pi_4 \circ \mathsf{c}_{\mathcal{T}_\mathcal{L}(X)} \circ  \mathcal{T}(\mathsf{c}_X)\circ  \mathcal{T}(\mathsf{c}_X), \\
        \pi_5 \circ  \mathcal{T}(\mathsf{c}_X) \circ \mathsf{c}_{\mathcal{T}_\mathcal{L}(X)} \circ  \mathcal{T}(\mathsf{c}_X) = \pi_2 = \pi_5 \circ \mathsf{c}_{\mathcal{T}_\mathcal{L}(X)} \circ  \mathcal{T}(\mathsf{c}_X)\circ  \mathcal{T}(\mathsf{c}_X), \\
        \pi_6 \circ  \mathcal{T}(\mathsf{c}_X) \circ \mathsf{c}_{\mathcal{T}_\mathcal{L}(X)} \circ  \mathcal{T}(\mathsf{c}_X) = \pi_6 = \pi_6 \circ \mathsf{c}_{\mathcal{T}_\mathcal{L}(X)} \circ  \mathcal{T}(\mathsf{c}_X)\circ  \mathcal{T}(\mathsf{c}_X), \\
        \pi_7 \circ  \mathcal{T}(\mathsf{c}_X) \circ \mathsf{c}_{\mathcal{T}_\mathcal{L}(X)} \circ  \mathcal{T}(\mathsf{c}_X) = \pi_4 = \pi_7 \circ \mathsf{c}_{\mathcal{T}_\mathcal{L}(X)} \circ  \mathcal{T}(\mathsf{c}_X)\circ  \mathcal{T}(\mathsf{c}_X), \\
        \pi_8 \circ  \mathcal{T}(\mathsf{c}_X) \circ \mathsf{c}_{\mathcal{T}_\mathcal{L}(X)} \circ  \mathcal{T}(\mathsf{c}_X) = \pi_8 = \pi_8 \circ \mathsf{c}_{\mathcal{T}_\mathcal{L}(X)} \circ  \mathcal{T}(\mathsf{c}_X)\circ  \mathcal{T}(\mathsf{c}_X).
    \end{align*}
    Thus, \(\mathsf{c}_X\) satisfies the desired Yang--Baxter identity, \(\mathcal{T}(\mathsf{c}_X) \circ \mathsf{c}_{\mathcal{T}_\mathcal{L}(X)} \circ  \mathcal{T}(\mathsf{c}_X)= \mathsf{c}_{\mathcal{T}_\mathcal{L}(X)} \circ  \mathcal{T}(\mathsf{c}_X)\circ  \mathcal{T}(\mathsf{c}_X)\), so the second diagram of (\ref{cidentitties2}) holds.

    We show now that the canonical flip is additive bundle morphism. Observe first that \(\mathcal{T}_\mathcal{L}\left( \mathsf{T}_\mathcal{L}(X) \right) = \mathcal{T}_\mathcal{L}\left( X \times \mathsf{L}(X) \right)\), as in Lemma \ref{gluing}.(\ref{gluing.mon.fun}), and that \(\mathsf{T}_\mathcal{L}(\mathcal{T}_\mathcal{L}(X)) = \mathcal{T}_\mathcal{L}(X) \times \mathsf{L}(\mathcal{T}_\mathcal{L}(X))\). Now, consider the isomorphisms which appear in the product preserving property of the functor \(\mathcal{T}_\mathcal{L}\). To distinguish from those of \(\mathcal{L}\), we will denote these as \(\omega^{\mathcal{T}_\mathcal{L}}_{X,Y}\colon \mathcal{T}_\mathcal{L}(X \times Y) \to \mathcal{T}_\mathcal{L}(X) \times \mathcal{T}_\mathcal{L}(Y)\). It is easy to check that the following equality holds:
    \begin{align}
        \omega^{\mathcal{T}_\mathcal{L}}_{X,Y} = \tau_{X, Y, \mathcal{L}(X), \mathcal{L}(Y)} \circ (1_{X \times Y} \times \omega_{X,Y})
    \end{align}
    Moreover, by Lemma \ref{gluing}.(\ref{gluing.mon.fun}), for any object \(X\) and any commutative monoid \(\mathsf{M}=(M, \bullet,e)\), we get an additive bundle isomorphism \((\omega^{\mathcal{T}_\mathcal{L}}_{X,M}, 1_{\mathcal{T}_\mathcal{L}(X)})\colon \mathcal{T}_\mathcal{L}(X \times \mathsf{M}) \to \mathcal{T}_\mathcal{L}(X) \times \mathcal{T}_\mathcal{L}(\mathsf{M})\). In particular, setting \(\mathsf{M} = \mathsf{L}(X)\), gives us an additive bundle isomorphism:
    \[(\omega^{\mathcal{T}_\mathcal{L}}_{X,\mathcal{L}(Y)}, 1_{\mathcal{T}_\mathcal{L}(X)}): \mathcal{T}_\mathcal{L}\left( \mathsf{T}_\mathcal{L}(X) \right) \to \mathcal{T}_\mathcal{L}(X) \times \mathcal{T}_\mathcal{L}(\mathsf{L}(X)).\]
    On the other hand, note that the underlying object of \(\mathcal{T}_\mathcal{L}(\mathsf{L}(X))\) is \(\mathcal{L}(X) \times \mathcal{L}\mathcal{L}(X)\), the underlying object of \(\mathsf{L}(\mathcal{T}_\mathcal{L}(X))\) is \(\mathcal{L}(X \times \mathcal{L}(X))\), and \(\omega^{-1}_{X, \mathcal{L}(X)}\colon \mathcal{T}_\mathcal{L}(\mathsf{L}(X)) \to \mathsf{L}(\mathcal{T}_\mathcal{L}(X))\) is a monoid morphism. Then, by Lemma \ref{gluing}.(\ref{gluing.mon}), we have that
    \[(1_{\mathcal{T}_\mathcal{L}(X)} \times \omega^{-1}_{X, \mathcal{L}(X)}, 1_{\mathcal{T}_\mathcal{L}(X)}): \mathcal{T}_\mathcal{L}(X) \times \mathcal{T}_\mathcal{L}(\mathsf{L}(X)) \to \mathcal{T}_\mathcal{L}(X) \times \mathsf{L}(\mathcal{T}_\mathcal{L}(X)) = \mathsf{T}_\mathcal{L}(\mathcal{T}_\mathcal{L}(X))\]
    is an additive bundle isomorphism. Now observe that by definition we have that:
    \begin{align*}
        \mathsf{c}_X = (1_{\mathcal{T}_\mathcal{L}(X)} \times \omega^{-1}_{X, \mathcal{L}(X)}) \circ \omega^{\mathcal{T}_\mathcal{L}}_{X,Y}.
    \end{align*}
    Composing the additive bundle morphisms \((1_{\mathcal{T}_\mathcal{L}(X)} \times \omega^{-1}_{X, \mathcal{L}(X)}, 1_{\mathcal{T}_\mathcal{L}(X)})\) and \((\omega^{\mathcal{T}_\mathcal{L}}_{X,\mathcal{L}(Y)}, 1_{\mathcal{T}_\mathcal{L}(X)})\) (recall that composition of bundle morphisms is given pointwise), we get that
    \[(\mathsf{c}_X, 1_{\mathcal{T}_\mathcal{L}(X)}) = (1_{\mathcal{T}_\mathcal{L}(X)} \times \omega^{-1}_{X, \mathcal{L}(X)}, 1_{\mathcal{T}_\mathcal{L}(X)}) \circ (\omega^{\mathcal{T}_\mathcal{L}}_{X,\mathcal{L}(Y)}, 1_{\mathcal{T}_\mathcal{L}(X)}): \mathcal{T}_\mathcal{L}\left( \mathsf{T}_\mathcal{L}(X) \right) \to \mathsf{T}_\mathcal{L}(\mathcal{T}_\mathcal{L}(X)) \]
    is an additive bundle morphism, as desired.

    We now check the identities involving the vertical lift \(\ell_X\). Let us first define the natural transformation \(\widehat{\mathsf{z}}_X\colon \mathcal{L}(X) \to \mathcal{L}(X \times \mathcal{L}(X))\) by \(\widehat{\mathsf{z}}_X = \omega^{-1}_{X, \mathcal{L}(X)} \circ \langle 0_X \circ t_{\mathcal{L}(X)}, \nu^{-1}_X \rangle\). Then, by definition, we have:
    \begin{equation*}
        \ell_X = \mathsf{z}_X \times  \widehat{\mathsf{z}}_X .
    \end{equation*}
    To show that the vertical lift is an additive bundle morphism, first note that
    \[\widehat{\mathsf{z}}_X: \mathsf{L}(X) \to \mathsf{L}(X \times \mathcal{L}(X)) = \mathsf{L}(\mathcal{T}_\mathcal{L}(X)) \]
    is a monoid morphism by construction. Thus, by Lemma \ref{gluing}.(\ref{gluing.mon}), we have that
    \[(\ell_X, \mathsf{z}_X)=(\mathsf{z}_X \times  \widehat{\mathsf{z}}_X, \mathsf{z}_X):  \mathsf{T}_\mathcal{L}(X) = X \times \mathsf{L}(X) \to \mathcal{T}_\mathcal{L}(X) \times \mathsf{L}(\mathcal{T}_\mathcal{L}(X)) = \mathsf{T}_\mathcal{L}(\mathcal{T}_\mathcal{L}(X)) \]
    is an additive bundle morphism, as required.

    We now show that diagrams (\ref{diag el1}) and (\ref{ellcidentitties}) commute. To do so, recall that \(\mathcal{T}^3_\mathcal{L}(X)\) is an octonary product. One can easily show that the following equalities hold:
    \begin{align}
        \pi_1 \circ  \mathcal{T}_\mathcal{L}(\ell_X) \circ \ell_X & = \pi_1 = \pi_1 \circ \ell_{\mathcal{T}_\mathcal{L}(X)} \circ \ell_X  ,                                                                                                       \\
        \pi_2 \circ  \mathcal{T}_\mathcal{L}(\ell_X) \circ \ell_X & = 0_X \circ t_{\mathcal{T}^3_\mathcal{L}(X)} = \pi_2 \circ \ell_{\mathcal{T}_\mathcal{L}(X)} \circ \ell_X   \label{eq liftline2},                                             \\
        \pi_3 \circ  \mathcal{T}_\mathcal{L}(\ell_X) \circ \ell_X & = 0_X \circ t_{\mathcal{T}^3_\mathcal{L}(X)} = \pi_3 \circ \ell_{\mathcal{T}_\mathcal{L}(X)} \circ \ell_X  \label{eq liftline3},                                              \\
        \pi_4 \circ  \mathcal{T}_\mathcal{L}(\ell_X) \circ \ell_X & = 0_{\mathcal{L}(X)} \circ t_{\mathcal{T}^3_\mathcal{L}(X)} = \pi_4 \circ \ell_{\mathcal{T}_\mathcal{L}(X)} \circ \ell_X  \label{eq liftline4},                               \\
        \pi_5 \circ  \mathcal{T}_\mathcal{L}(\ell_X) \circ \ell_X & = 0_X \circ t_{\mathcal{T}^3_\mathcal{L}(X)} = \pi_5 \circ \ell_{\mathcal{T}_\mathcal{L}(X)} \circ \ell_X  \label{eq liftline5},                                              \\
        \pi_6 \circ  \mathcal{T}_\mathcal{L}(\ell_X) \circ \ell_X & = 0_{\mathcal{L}(X)} \circ t_{\mathcal{T}^3_\mathcal{L}(X)} = \pi_6 \circ \ell_{\mathcal{T}_\mathcal{L}(X)} \circ \ell_X  \label{eq liftline6},                               \\
        \pi_7 \circ  \mathcal{T}_\mathcal{L}(\ell_X) \circ \ell_X & = 0_{\mathcal{L}(X)} \circ t_{\mathcal{T}^3_\mathcal{L}(X)} = \pi_7 \circ \ell_{\mathcal{T}_\mathcal{L}(X)} \circ \ell_X  \label{eq liftline7},                               \\
        \pi_8 \circ  \mathcal{T}_\mathcal{L}(\ell_X) \circ \ell_X & = \mathcal{L}(\nu^{-1}_X) \circ \nu^{-1}_X = \nu^{-1}_{\mathcal{L}(X)} \circ \nu^{-1}_X = \pi_7 \circ \ell_{\mathcal{T}_\mathcal{L}(X)} \circ \ell_X   . \label{eq liftline8}
    \end{align}
    For lines (\ref{eq liftline2}) to (\ref{eq liftline6}), we use the fact that \(\nu\) and \(\omega^{-1}\) are monoid morphisms. For lines (\ref{eq liftline6}) and (\ref{eq liftline7}), we also use the equality \(0_{\mathcal{L}(X)} = \mathcal{L}(0_X) \circ t^{-1}_{\mathcal{L}(\ast)}\) (\ref{diag L+,L0}). For line (\ref{eq liftline8}), we use the equality \(\mathcal{L}(\nu_X)=\nu_{\mathcal{L}(X)}\). From this, we conclude that \(\mathcal{T}_\mathcal{L}(\ell_X) \circ \ell_X=\ell_{\mathcal{T}_\mathcal{L}(X)} \circ \ell_X\), and so (\ref{diag el1}) commutes.

    We now show the compatibilities between the vertical lift and the canonical flip. It is easy to check that the following equality holds:
    \[\langle 1_X, 0_X \circ t_X \rangle \times \langle 0_X \circ t_{\mathcal{L}(X)}, \nu^{-1}_X \rangle \circ \tau_{X,\mathcal{L}(X),\mathcal{L}(X),\mathcal{L}(X)} = \langle 1_X, 0_X \circ t_X \rangle \times \langle 0_X \circ t_{\mathcal{L}(X)}, \nu^{-1}_X \rangle.\]
    It follows that \(\mathsf{c}_X \circ \ell_X = \ell_X\), so the first diagram of (\ref{ellcidentitties}) commutes. For the remaining diagram, we view \(\mathcal{T}^2_\mathcal{L}(X)\) as a quaternary product with projections:
    \begin{gather*}
        \pi_1: \mathcal{T}^2_\mathcal{L}(X) \to X \quad \pi_2: \mathcal{T}^2_\mathcal{L}(X) \to \mathcal{L}(X) \quad \pi_3: \mathcal{T}^2_\mathcal{L}(X) \to \mathcal{L}(X) \quad \pi_4: \mathcal{T}^2_\mathcal{L}(X) \to \mathcal{L}^2(X).
    \end{gather*}
    Then, one may easily show that the following equalities hold:
    \begin{align*}
        \pi_1 \circ  \mathcal{T}_\mathcal{L}(\ell_X) \circ  \mathsf{c}_X & = \pi_1 = \pi_1 \circ \mathsf{c}_{\mathcal{T}_\mathcal{L}(X)} \circ  \mathcal{T}(\mathsf{c}_X)\circ  \ell_{\mathcal{T}_\mathcal{L}(X)} ,                                                     \\
        \pi_2 \circ \mathcal{T}_\mathcal{L}(\ell_X) \circ  \mathsf{c}_X  & = 0_X \circ t_{\mathcal{T}^3_\mathcal{L}(X)} = \pi_2 \circ \mathsf{c}_{\mathcal{T}_\mathcal{L}(X)} \circ  \mathcal{T}(\mathsf{c}_X)\circ \ell_{\mathcal{T}_\mathcal{L}(X)} ,                 \\
        \pi_3 \circ  \mathcal{T}_\mathcal{L}(\ell_X) \circ  \mathsf{c}_X & = 0_X \circ t_{\mathcal{T}^3_\mathcal{L}(X)} = \pi_3 \circ \mathsf{c}_{\mathcal{T}_\mathcal{L}(X)} \circ  \mathcal{T}(\mathsf{c}_X)\circ \ell_{\mathcal{T}_\mathcal{L}(X)} ,                 \\
        \pi_4 \circ \mathcal{T}_\mathcal{L}(\ell_X) \circ  \mathsf{c}_X  & = \pi_3 = \pi_4 \circ \mathsf{c}_{\mathcal{T}_\mathcal{L}(X)} \circ  \mathcal{T}(\mathsf{c}_X)\circ  \ell_{\mathcal{T}_\mathcal{L}(X)} ,                                                     \\
        \pi_5 \circ  \mathcal{T}_\mathcal{L}(\ell_X) \circ  \mathsf{c}_X & = \pi_2 = \pi_5 \circ \mathsf{c}_{\mathcal{T}_\mathcal{L}(X)} \circ  \mathcal{T}(\mathsf{c}_X)\circ  \ell_{\mathcal{T}_\mathcal{L}(X)} ,                                                     \\
        \pi_6 \circ \mathcal{T}_\mathcal{L}(\ell_X) \circ  \mathsf{c}_X  & = 0_{\mathcal{L}(X)} \circ t_{\mathcal{T}^3_\mathcal{L}(X)} = \pi_6 \circ \mathsf{c}_{\mathcal{T}_\mathcal{L}(X)} \circ  \mathcal{T}(\mathsf{c}_X)\circ \ell_{\mathcal{T}_\mathcal{L}(X)} ,  \\
        \pi_7 \circ \mathcal{T}_\mathcal{L}(\ell_X) \circ  \mathsf{c}_X  & = 0_{\mathcal{L}(X)} \circ t_{\mathcal{T}^3_\mathcal{L}(X)} = \pi_7 \circ \mathsf{c}_{\mathcal{T}_\mathcal{L}(X)} \circ  \mathcal{T}(\mathsf{c}_X)\circ  \ell_{\mathcal{T}_\mathcal{L}(X)} , \\
        \pi_8 \circ \mathcal{T}_\mathcal{L}(\ell_X) \circ  \mathsf{c}_X  & = \mathcal{L}(\nu^{-1}_X)= \nu^{-1}_{\mathcal{L}(X)} = \pi_8 \circ \mathsf{c}_{\mathcal{T}_\mathcal{L}(X)} \circ  \mathcal{T}(\mathsf{c}_X)\circ  \ell_{\mathcal{T}_\mathcal{L}(X)}.
    \end{align*}
    Here, we used the Yang--Baxter identities involving symmetries and interchange, and in the last line, we again use the equality \(\mathcal{L}(\nu_X)=\nu_{\mathcal{L}(X)}\). It follows that \( \mathcal{T}_\mathcal{L}(\ell_X) \circ  \mathsf{c}_X = \mathsf{c}_{\mathcal{T}_\mathcal{L}(X)} \circ  \mathcal{T}(\mathsf{c}_X)\circ  \ell_{\mathcal{T}_\mathcal{L}(X)}\), so the second diagram of (\ref{ellcidentitties}) commutes.

    It remains to show the universal property of the vertical lift. Translating (\ref{universalift}) in our setting, we must show that the following diagram is a pullback:
    \begin{equation}\label{universalift-L}\begin{gathered}
            \xymatrixcolsep{6pc}\xymatrix{ X \times (\mathcal{L}(X) \times \mathcal{L}(X)) \ar[d]_-{\pi_1}  \ar[r]^-{\langle 1_X \times \pi_2, \widehat{\mathsf{z}}_X \circ \pi_1 \circ \pi_2 \rangle} & (X \times \mathcal{L}(X)) \times \mathcal{L}\left(X \times \mathcal{L}(X)\right) \ar[d]^-{\pi_1 \times \mathcal{L}(\pi_1)}  \\
            X \ar[r]_-{\mathsf{z}_X} & X \times \mathcal{L}(X).
            }  \end{gathered}\end{equation}
    Suppose that we have morphisms \(g\colon Y \to X\) and \(f\colon Y \to (X \times \mathcal{L}(X)) \times \mathcal{L}\left(X \times \mathcal{L}(X)\right)\) such that \(\mathsf{z}_X \circ g = \left( \pi_1 \times \mathcal{L}(\pi_1) \right) \circ f\). We then have unique morphisms \(f_1\colon Y \to X\), \(f_2\colon Y \to \mathcal{L}(X)\), \(f_3\colon Y \to \mathcal{L}(X)\), and \(f_4\colon Y \to \mathcal{L}(X)\) such that:
    \[ f = \bigl\langle \langle f_1, f_2 \rangle, \omega^{-1}_{X, \mathcal{L}(X)} \circ \langle f_3, \nu^{-1}_X \circ f_4 \rangle \bigr\rangle .\]
    The assumptions on \(f\) and \(g\) imply that \(f_1 = g\) and \(f_3 = 0_X \circ t_Y\). Consider the morphism \(h\colon Y \to X \times (\mathcal{L}(X) \times \mathcal{L}(X))\) defined by \(h \coloneq \langle f_1, \langle f_4, f_2 \rangle \rangle\). We then have \(\pi_1 \circ h = g\), and we compute:
    \begin{align*}
          \langle 1_X \times \pi_2, \widehat{\mathsf{z}}_X \circ \pi_1 \circ \pi_2 \rangle \circ h  &= \langle (1_X \times \pi_2) \circ h, \widehat{\mathsf{z}}_X \circ \pi_1 \circ \pi_2 \circ h \rangle                                            \\
         & = \langle (1_X \times \pi_2) \circ h, \widehat{\mathsf{z}}_X \circ \pi_1 \circ \pi_2 \circ h \rangle = \langle \langle f_1, f_2 \rangle, \widehat{\mathsf{z}}_X \circ f_4 \rangle                                                         \\
         & =  \langle \langle f_1, f_2 \rangle, \langle 0_X \circ t_{\mathcal{L}(X)}, \nu^{-1}_X \rangle \circ f_4 \rangle = \langle \langle f_1, f_2 \rangle, \langle 0_X \circ t_{\mathcal{L}(X)} \circ f_4, \nu^{-1}_X \circ f_4 \rangle  \rangle \\
         & = \langle \langle f_1, f_2 \rangle, \langle 0_X \circ t_Y, \nu^{-1}_X \circ f_4 \rangle  \rangle = \langle \langle f_1, f_2 \rangle, \langle f_3, \nu^{-1}_X \circ f_4 \rangle  \rangle = f.
    \end{align*}
    Now, suppose that there exists another morphism \(k\colon Y \to X \times (\mathcal{L}(X) \times \mathcal{L}(X))\) such that \(\pi_1 \circ k = g\) and \(\langle 1_X \times \pi_2, \widehat{\mathsf{z}}_X \circ \pi_1 \circ \pi_2 \rangle \circ k = f\). It is easy to check that \(\pi_1 \circ k = g\), \(\pi_1 \circ \pi_2 \circ k = f_4\) and \(\pi_2 \circ \pi_2 \circ k = f_2\), and thus, that \(k = \langle f_1, \langle f_4, f_2 \rangle \rangle = h\). We conclude that diagram (\ref{universalift-L}) is indeed a pullback diagram, as desired.

    In conclusion, \(\mathbb{T}_\mathcal{L}\) is a (Rosický) tangent structure, and \((\mathbb{X}, \mathbb{T}_\mathcal{L})\) is a cartesian (Rosický) tangent category.
\end{proof}

To conclude this section, let us apply this construction to our main examples of linear assignments:

\begin{example}\label{ex tanL-trival} For any category with finite products, applying Theorem \ref{theo:linearprojtotangent} to the terminal linear assignment from Example \ref{ex lin-trival} gives us (up to isomorphism) the trivial tangent structure from Example \ref{ex trivial-tancat}, that is, \(\mathcal{T}_{\mathcal{L}} \cong 1_\mathbb{X}\).
\end{example}

\begin{example}\label{ex tanL-biproduct} For any semi-additive (resp.\ additive) category, applying Theorem \ref{theo:linearprojtotangent} to the identity linear assignment from Example \ref{ex lin-biprod} results precisely in the canonical (Rosický) tangent structure of a semi-additive category given by the diagonal functor from Example \ref{ex biprod-tancat}, that is, \(\mathcal{T}_{1_\mathbb{X}} = \Delta\).
\end{example}

\begin{example}\label{ex tanL-group} Applying Theorem \ref{theo:linearprojtotangent} to the additive assignment given by abelianization of groups yields a cartesian Rosický tangent structure on the category of groups. This is one of the main new observations of this paper. Let us review this structure in details. Abelian group structures will be denoted additively, and we implicitly identify \(\ab(G\times H)\) with \(\ab(G)\times \ab(H)\), for all groups \(G\), \(H\). For a group \(G\), elements of its abelianization \(\ab(G)\) will be denoted by \([g]\) for all \(g \in G\). Then, \(\Grp\) is a cartesian Rosický tangent category whose tangent bundle functor is given by \(\mathcal{T}_{\mathcal{L}_\ab}(G) = G \times \ab(G),\) and the rest of the Rosický tangent structure is given as follows:
    \begin{gather*}
        \mathsf{p}_G(g,[h]) = g \quad \mathsf{s}_G(g,[h_1],[h_2]) = (g, [h_1]+[h_2])\quad \mathsf{z}_G(g) = (g,0) \quad \mathsf{n}_G(g,[h]) = (g, -[h]) \\
        \ell_G(g,[h]) = \left( (g,0),(0,[h])  \right) \qquad  \mathsf{c}_G \left( (g,[h]),([k],[j])  \right)  = \left( (g,[k]),([h],[j]) \right) .
    \end{gather*}
\end{example}

\section{Linear Algebras, Differential Objects and Bundles}

In a tangent category, an important class of objects are the \emph{differential bundles}, which formalize the notion of smooth vector bundles in a tangent category. For an in-depth introduction to differential bundles, we refer the reader to \cite{ching2024characterization,cockettCruttwellDiffBundles,cruttwell17,MacAdamVectorBundles}. In this section, we show that for the tangent structure induced by a linear assignment, differential bundles are closely related to a special class of commutative monoids associated to the linear assignment, which we call \emph{linear algebras}, and which can be characterized as fixed points for the linear assignment.

\begin{definition}\label{def L-alg} Let \(\mathcal{L}\) be a linear (resp.\ additive) assignment on a category \(\mathbb{X}\) with finite products.
    \begin{enumerate}[{\em (i)}]
        \item A \textbf{linear algebra}, or simply an \textbf{\(\mathcal{L}\)-algebra}, is a pair \((A, \mathsf{a})\) consisting of an object \(A\) and an isomorphism \(\mathsf{a}\colon \mathcal{L}(A) \to A\) such that \(\mathcal{L}(\mathsf{a}) = \nu_A\).
        \item If \((A, \mathsf{a})\) and \((A^\prime, \mathsf{a}^\prime)\) are \(\mathcal{L}\)-algebras, an \textbf{\(\mathcal{L}\)-algebra morphism} \(f\colon (A, \mathsf{a}) \to (A^\prime, \mathsf{a}^\prime)\) is a morphism \(f\colon A \to A^\prime\) such that the following diagram commutes:
              \begin{equation*}\label{}\begin{gathered}
                      \xymatrixcolsep{2pc}\xymatrix{  \mathcal{L}(A) \ar[r]^-{\mathcal{L}(f)} \ar[d]_-{\mathsf{a}}  & \mathcal{L}(A') \ar[d]^-{\mathsf{a}^\prime}  \\
                      A \ar[r]_-{f} &  A^\prime }  \end{gathered}\end{equation*}
    \end{enumerate}
    We denote by \(\mathcal{L}\-\mathsf{ALG}\) the category of \(\mathcal{L}\)-algebras and \(\mathcal{L}\)-algebra morphisms between them.
\end{definition}

For a linear (resp.\ additive) assignment \(\mathcal{L}\), every \(\mathcal{L}\)-algebra is canonically a commutative monoid (resp.\ abelian group). To see this, observe that if \(\mathsf{M}=(M, \bullet, e)\) is a commutative monoid and \(f\colon M \to A\) is an isomorphism, then the tuple:
\begin{align*}
    A_f = (A, (f \times f) \circ \bullet \circ f^{-1}, f \circ e)
\end{align*}
is a commutative monoid and \(f\colon \mathsf{M} \to A_f\) is a monoid isomorphism. Similarly, if \(\mathsf{G}=(G, \bullet, e, \iota)\) is an abelian group and \(f\colon G \to A\) is an isomorphism, then the tuple:
\begin{align*}
    A_f = (A, (f \times f) \circ \bullet \circ f^{-1}, f \circ e, f\circ \iota \circ f^{-1})
\end{align*}
is an abelian group and \(f\colon \mathsf{G} \to A_f\) is a group monoid isomorphism. Applying this construction to linear algebras, we immediately get the following:

\begin{lemma} Let \(\mathcal{L}\) be a linear (resp.\ additive) assignment on a category \(\mathbb{X}\) with finite products.
    \begin{enumerate}[{\em (i)}]
        \item If \((A, \mathsf{a})\) is an \(\mathcal{L}\)-algebra, then \(A_\mathsf{a}\) is a commutative monoid (resp.\ abelian group) and \(\mathsf{a}\colon \mathsf{L}(A) \to A_\mathsf{a}\) is a monoid (resp.\ group) isomorphism.
        \item If \(f\colon (A, \mathsf{a}) \to (A^\prime, \mathsf{a}^\prime)\) is an \(\mathcal{L}\)-algebra morphism, then \(f\colon A_\mathsf{a} \to A^\prime_{\mathsf{a}^\prime}\) is a monoid (resp.\ group) morphism.
    \end{enumerate}
    This induces a functor \(\mathcal{V}\colon \mathcal{L}\-\mathsf{ALG} \to \CMON[\X]\) (resp.\ \(\mathcal{V}\colon \mathcal{L}\-\mathsf{ALG} \to \AbG[\X]\)) defined on objects by \(\mathcal{V}(A,\mathsf{a}) = A_\mathsf{a}\) and on morphisms by \(\mathcal{V}(f) = f\).
\end{lemma}

As a shorthand, we will write the monoid (resp.\ group) structure of \(A_\mathsf{a}\) by \(+_{\mathsf{a}}\colon A \times A \to A\) and \(0_\mathsf{a}\colon \ast \to A\) (and \(-_\mathsf{a}\colon A \to A\)). Observe that, for every object \(X\) of \(\mathbb X\), \((\mathcal{L}(X),\nu_X)\) is an \(\mathcal{L}\)-algebra. Furthermore, seeing \(\mathcal L(X)\) as an \(\mathcal L\)-algebra, and considering the associated monoid (resp.\ group) structure on \(\mathcal L(X)\) yields precisely the monoid (resp.\ group) structure given by the linear (resp.\ additive) assignment. To see this, note that for commutative monoids \(\mathsf{M}=(M, \bullet, e)\) and \(\mathsf{M}^\prime=(M^\prime, \bullet^\prime, e^\prime)\), if \(f\colon \mathsf{M} \to \mathsf{M}^\prime\) is a monoid isomorphism, then \(M^\prime_f = \mathsf{M}^\prime\) (and similarly for abelian groups). We deduce the following:

\begin{lemma} \label{lem:Lnualg} Let \(\mathcal{L}\) be a linear (resp.\ additive) assignment on a category \(\mathbb{X}\) with finite products.
    \begin{enumerate}[{\em (i)}]
        \item \label{lem:Lnualg.1} For every object \(X\), \((\mathcal{L}(X), \nu_X)\) is an \(\mathcal{L}\)-algebra and \(\mathcal{L}(X)_{\nu_X} = \mathsf{L}(X)\), that is, \(+_{\nu_X} = +_X\) and \(0_{\nu_X} = 0\) (and \(-_{\nu_X} = -_X)\).
        \item For every morphism \(f\colon X \to X^\prime\), \(\mathcal{L}(f)\colon (\mathcal{L}(X), \nu_X) \to (\mathcal{L}(X^\prime), \nu_{X^\prime})\) is a \(\mathcal{L}\)-algebra morphism.
    \end{enumerate}
    This induces a functor \(\mathcal{L}^\sharp\colon \X \to \mathcal{L}\-\mathsf{ALG}\) defined on objects by \(\mathcal{L}^\sharp(X) = (\mathcal{L}(X), \nu_X)\) and on morphisms by \(\mathcal{L}^\sharp(f) = \mathcal{L}(f)\). Moreover, the induced linear (resp.\ additive) projector \(\mathsf{L}\) factors through \(\mathcal{L}\-\mathsf{ALG}\) in the sense that \(\mathsf{L} = \mathcal{V} \circ \mathcal{L}^\sharp\).
\end{lemma}

Let us now identify the linear algebras of our main examples of linear assignments:

\begin{example} For the terminal additive assignment on a category with finite products, the only \(\mathcal{L}_\ast\)-algebra, up to isomorphism, is the terminal object \(\ast\), and so \(\mathcal{L}_\ast\-\mathsf{ALG}\) is the terminal category (with one object and one morphism).
\end{example}

\begin{example} For the identity linear (resp.\ additive) assignment on a semi-additive (resp.\ additive) category \(\mathbb{X}\), every object is a \(1_{\mathbb{X}}\)-algebra (where the \(1_{\mathbb{X}}\)-algebra structure is simply the identity), and so, \(1_{\mathbb{X}}\-\mathsf{ALG} \cong \mathbb{X}\).
\end{example}

\begin{example}\label{ex L-alg-group} For the abelianization \(\mathcal{L}_\ab\) of groups, seen as an additive assignment, the \(\mathcal{L}_\ab\)-algebras are precisely the abelian groups. For an abelian group \(A\), its \(\mathcal{L}_\ab\)-algebra structure is given by the canonical group isomorphism \(A \cong \ab(A)\). We then have \(\mathcal{L}_\ab\-\mathsf{ALG} \cong \Ab\).
\end{example}

Let us now turn our attention to differential bundles, which are the analogues of smooth vector bundles for tangent categories. There are various ways to define differential bundles. Here, we have chosen the original definition found in \cite{cockettCruttwellDiffBundles}. Equivalent definitions were given by MacAdam in \cite{MacAdamVectorBundles} and by Ching in \cite{ching2024characterization}. In a tangent category \((\mathbb{X}, \mathbb{T})\), a \textbf{differential bundle} \cite[Definition 2.3]{cockettCruttwellDiffBundles} is a pair \((\mathsf{E}, \lambda)\), consisting of an additive bundle \(\mathsf{E}=(X, E, \mathsf{q}, \sigma, \zeta)\) (where we call \(\mathsf{q}\) the \textbf{projection}, \(\sigma\) the \textbf{sum}, and \(\zeta\) the \textbf{zero}) and a morphism \(\lambda\colon E \to \mathcal{T}(E)\), called the \textbf{lift} satisfying the following axioms:
\begin{itemize}
    \item The tangent bundle functor preserves all pullbacks powers of the projection, that is, for all \(n\), \(m \in \mathbb{N}\), \(\mathcal{T}^m(E_n)\) is the pullback of \(n\) copies of \(\mathcal{T}(\mathsf{q})\). In particular, this implies that \(\mathcal{T}(\mathsf{E}) = (\mathcal{T}(X), \mathcal{T}(E), \mathcal{T}(\mathsf{q}), \mathcal{T}(\sigma), \mathcal{T}(\zeta))\) is an additive bundle.
    \item Both \((\lambda, \mathsf{z}_X)\colon \mathsf{E} \to \mathcal{T}(\mathsf{E})\) and \((\lambda, \zeta)\colon \mathsf{E} \to \mathsf{T}(E)\) are additive bundle morphisms, and the following diagram commutes:
          \begin{equation}\label{diag lift-dbun}\begin{gathered}
                  \xymatrixcolsep{2pc}\xymatrix{  E \ar[r]^-{\lambda} \ar[d]_-{\lambda}  & \mathcal{T}(E) \ar[d]^-{\mathcal{T}(\lambda)}  \\
                  \mathcal{T}(E) \ar[r]_-{\ell_E} &  \mathcal{T}\mathcal{T}(E) }  \end{gathered}\end{equation}
    \item The lift is universal in the sense that the following diagram is a pullback:
          \begin{equation}\label{universaliftdbun}\begin{gathered}
                  \xymatrixcolsep{5pc}\xymatrix{ E_2 \ar[d]_-{\mathsf{p}_X \circ \rho_j}  \ar[r]^-{\langle \lambda \circ \rho_1, 0_E \circ \rho_2 \rangle} & \mathcal{T}(E_2) \ar[r]^-{\mathcal{T}(\sigma)} & \mathcal{T}(E) \ar[d]^-{\mathcal{T}(q)}  \\
                  X \ar[rr]_-{\mathsf{z}_X} & & \mathcal{T}(X)
                  }  \end{gathered}\end{equation}
          and all powers of the tangent bundle functor \(\mathcal{T}^m\) preserve these pullbacks.
\end{itemize}
There is of course an obvious notion of a differential bundle with negatives, where we upgrade the definition from an additive bundle to an abelian group bundle. However, it turns out that in a Rosický tangent category, every differential bundle is canonically an abelian group bundle, where the negative for the differential bundle is induced from the negative of the tangent bundle. Thus, in a Rosický tangent category, the notions of differential bundles and differential bundles with negatives coincide \cite[Proposition 2.13]{cruttwell2023differential}, and we only need to consider differential bundles as additive bundles.

Morphisms of differential bundles \cite[Definition 2.3]{cockettCruttwellDiffBundles} are bundle morphisms that also preserve the lift. That is to say, for two differential bundles \((\mathsf{E}, \lambda)\) and \((\mathsf{E}^\prime, \lambda^\prime)\), with underlying bundles \(\mathsf{q}\colon E \to X\) and \(\mathsf{q}^\prime\colon E^\prime \to X^\prime\) respectively, a \textbf{linear bundle morphism} \((f,g)\colon (\mathsf{E}, \lambda) \to (\mathsf{E}^\prime, \lambda^\prime)\) is a bundle morphism \((f,g)\colon \mathsf{q} \to \mathsf{q}^\prime\) such that the following diagram commutes:
\begin{equation}\label{diag dbun-morph}\begin{gathered}
        \xymatrixcolsep{2pc}\xymatrix{  E \ar[r]^-{f} \ar[d]_-{\lambda}  & E^\prime \ar[d]^-{\lambda^\prime}  \\
        \mathcal{T}(E) \ar[r]_-{\mathcal{T}(f)} &  \mathcal{T}(E^\prime) }  \end{gathered}\end{equation}
In this setting, \((f,g)\colon \mathsf{E} \to \mathsf{E}^\prime\) is automatically an additive bundle morphism (and so, an abelian group bundle morphism in the Rosický setting) \cite[Proposition 2.16]{cockettCruttwellDiffBundles}. In the special case where \(X=Y\), so for differential bundles \((\mathsf{E}, \lambda)\) and \((\mathsf{E}^\prime, \lambda^\prime)\) over the same object \(X\), we say that \(f\colon (\mathsf{E}, \lambda) \to (\mathsf{E}^\prime, \lambda^\prime)\) is a \textbf{linear bundle morphism over \(X\)} if \((f,1_X)\colon (\mathsf{E}, \lambda) \to (\mathsf{E}^\prime, \lambda^\prime)\) is a linear bundle morphism. For a (Rosický) tangent category \((\mathbb{X}, \mathbb{T})\), we denote by \(\mathsf{DBUN}[(\mathbb{X}, \mathbb{T})]\) its category of differential bundles and linear bundle morphisms, and for each object \(X\), we denote by \(\mathsf{DBUN}[(\mathbb{X}, \mathbb{T})]_X\) the category of differential bundles over \(X\) and linear bundle morphisms over \(X\).

We now show that linear algebras give rise to differential bundles. Let \(\mathcal{L}\) be a linear assignment on a category \(\mathbb{X}\) with finite products. For any object \(X\) and any \(\mathcal{L}\)-algebra \((A, \mathsf{a})\), define the morphism \(\lambda_{X,\mathsf{a}}\colon X \times A \to \mathcal{T}_\mathcal{L}(X \times A)\) to be the following composite:
\begin{equation*}\label{}\begin{gathered}
        \xymatrixcolsep{5pc}\xymatrix{ X \times A \ar[rr]^-{\lambda_{X,\mathsf{a}}} \ar[rd]|-{\langle 1_X, 0_\mathsf{a} \circ t_X \rangle \times \langle 0_X \circ t_{A}, \mathsf{a}^{-1} \rangle} && (X \times A) \times \mathcal{L}\left(X \times A\right)\\
        & (X \times A) \times (\mathcal{L}(X) \times \mathcal{L}(A)) \ar[ru]|-{1_{X \times A} \times \omega^{-1}_{X, A}}
        }
    \end{gathered}\end{equation*}

\begin{proposition}\label{prop Lalg-to-dbun} Let \(\mathcal{L}\) be a linear assignment on a category \(\mathbb{X}\) with finite products.
    \begin{enumerate}[{\em (i)}]
        \item  For any object \(X\) and any \(\mathcal{L}\)-algebra \((A, \mathsf{a})\), \((X \times A_\mathsf{a}, \lambda_{X,\mathsf{a}})\) is a differential bundle.
        \item For any morphism \(g\colon X \to X^\prime\) and any \(\mathcal{L}\)-algebra morphism \(f\colon (A, \mathsf{a}) \to (A^\prime, \mathsf{a}^\prime)\), \((g \times f, g)\colon (X \times A_\mathsf{a}, \lambda_{X,\mathsf{a}}) \to (X^\prime \times A^\prime_{\mathsf{a}^\prime}, \lambda_{X^\prime,\mathsf{a}^\prime})\) is a linear bundle morphism.
    \end{enumerate}
    As such, this induces a functor \(\mathcal{D}\colon \mathbb{X} \times \mathcal{L}\-\mathsf{ALG} \to \mathsf{DBUN}[(\mathbb{X}, \mathbb{T}_\mathcal{L})]\), defined on objects by $\mathcal{D}(X, (A, \mathsf{a})) = (X \times A_\mathsf{a}, \lambda_{X,\mathsf{a}})$, and on morphisms by \(\mathcal{D}(g,f) = (g \times f, g)\). Similarly, for every object \(X\), we also have a functor \(\mathcal{D}_X\colon \mathcal{L}\-\mathsf{ALG} \to \mathsf{DBUN}[(\mathbb{X}, \mathbb{T}_\mathcal{L})]_X\) defined on objects by \(\mathcal{D}_X(A, \mathsf{a}) = (X \times A_\mathsf{a}, \lambda_{X,\mathsf{a}})\) and on morphisms by \(\mathcal{D}_X(f) = g \times f\).
\end{proposition}

\begin{proof} Consider an object \(X\) and an \(\mathcal{L}\)-algebra \((A, \mathsf{a})\). By Lemma \ref{gluing}.(\ref{gluing.proj}) and (\ref{gluing.mon}), we already know that \(X \times A_\mathsf{a}\) is an additive bundle morphism such that all powers of \(\mathcal{T}_\mathcal{L}\) preserve the pullbacks of copies of \(\pi_1\colon X \times A \to X\). Let us show that the lift induces two additive bundle morphisms \(X \times A_\mathsf{a}\to \mathcal T(X \times A_\mathsf{a})\). First, note that \(\mathsf{T}_\mathcal{L}(X \times A) = (X \times A) \times \mathsf{L}(X \times A)\). Define the morphism \(\zeta\colon X \to X \times A\) by \(\zeta \coloneq \langle 1_X, 0_\mathsf{a} \circ t_X \rangle\), and \(\widehat{\zeta}\colon A \to \mathcal{L}\left(X \times A\right)\) by \(\widehat{\zeta} \coloneq \omega^{-1}_{X, A} \circ \langle 0_X \circ t_{A}, \mathsf{a}^{-1} \rangle\). By definition, we have: \(\lambda_{X,\mathsf{a}} = \zeta \times \widehat{\zeta}\).
    Now, observe that \(\widehat{\zeta}\colon A_\mathsf{a} \to \mathsf{L}\left(X \times A\right)\) is a monoid morphism by construction. Thus, by Lemma \ref{gluing}.(\ref{gluing.mon}), the following is an additive bundle morphism:
    \[(\lambda_{X,\mathsf{a}}, \zeta)=(\zeta \times  \widehat{\zeta}, \zeta):  X \times A_\mathsf{a} \to (X \times A) \times \mathsf{L}(X \times A) =\mathsf{T}_\mathcal{L}(X \times A) .\]
    On the other hand, recall from the proof of Theorem \ref{theo:linearprojtotangent} that, for any object \(X\), and any commutative monoid \(\mathsf{M}=(M, \bullet,e)\), we have that \(({\omega^{\mathcal{T}_\mathcal{L}}}^{-1}_{X,M}, 1_{\mathcal{T}_\mathcal{L}(X)})\colon \mathcal{T}_\mathcal{L}(X \times \mathsf{M}) \to \mathcal{T}_\mathcal{L}(X) \times \mathcal{T}_\mathcal{L}(\mathsf{M})\) is an additive bundle isomorphism. So, setting \(\mathsf{M} = A_\alpha\), and using the inverse of the isomorphism above, we get the additive bundle isomorphism:
    \[({\omega^{\mathcal{T}_\mathcal{L}}}^{-1}_{X,A}, 1_{\mathcal{T}_\mathcal{L}(X)})^{-1}: \mathcal{T}_\mathcal{L}(X) \times \mathcal{T}_\mathcal{L}(A_\alpha) \to \mathcal{T}_\mathcal{L}(X \times A_\alpha) . \]
    Now, define the morphism \(\tilde{\zeta}\colon A \to A \times \mathcal{L}(A)\) by \(\tilde{\zeta} \coloneq \langle 0_\mathsf{a} \circ t_{A}, \mathsf{a}^{-1} \rangle\). Note that \(\tilde{\zeta}\colon A_\mathsf{a} \to A_\alpha \times \mathcal{L}(A_\alpha) = \mathcal{T}_\mathcal{L}(A_\mathsf{a})\) is a monoid morphism by construction. Thus, by Lemma \ref{gluing}.(\ref{gluing.mon}), the following is an additive bundle morphism:
    \[ ( \mathsf{z}_X \times \tilde{\zeta}, \mathsf{z}_X): X \times A_\mathsf{a} \to \mathcal{T}_\mathcal{L}(X) \times \mathcal{T}_\mathcal{L}(A_\mathsf{a}) \]
    Now, it is easy to check that the following equality holds:
    \begin{align*}
        {\omega^{\mathcal{T}_\mathcal{L}}}^{-1}_{X,A} \circ (\mathsf{z}_X \times \tilde{\zeta}) = \lambda_{X,\mathsf{a}} ,
    \end{align*}
    and thus, by composing additive bundle morphisms, we get that:
    \[ (\lambda_{X,\mathsf{a}}, \mathsf{z}_X) = ({\omega^{\mathcal{T}_\mathcal{L}}}^{-1}_{X,A}, 1_{\mathcal{T}_\mathcal{L}(X)}) \circ ( \mathsf{z}_X \times \tilde{\zeta}, \mathsf{z}_X)'\colon X \times A_\mathsf{a} \to \mathcal{T}_\mathcal{L}(X \times A_\alpha) \]
    is a bundle morphism, as desired.

    Let us now show that diagram (\ref{diag lift-dbun}) commutes. This is again analogous to an argument from the proof of Theorem \ref{theo:linearprojtotangent}. We first observe that:
    \[ \mathcal{T}^2_\mathcal{L}(X \times A) \cong \left( (X \times A) \times (\mathcal{L}(X) \times \mathcal{L}(A)) \right) \times \left( (\mathcal{L}(X) \times \mathcal{L}(A)) \times (\mathcal{L}^2(X) \times \mathcal{L}^2(A)) \right) .\textbf{} \]
    Thus, by a slight abuse of notation, we may view \(\mathcal{T}^2_\mathcal{L}(X \times A)\) as an octonary product with projections:
    \begin{gather*}
        \pi_1\colon \mathcal{T}^2_\mathcal{L}(X \times A) \to X \quad \pi_2\colon \mathcal{T}^2_\mathcal{L}(X \times A) \to A \quad \pi_3\colon \mathcal{T}^2_\mathcal{L}(X \times A) \to \mathcal{L}(X) \\
        \pi_4\colon \mathcal{T}^2_\mathcal{L}(X \times A) \to \mathcal{L}(A) \quad \pi_5\colon \mathcal{T}^2_\mathcal{L}(X \times A) \to \mathcal{L}(X) \quad \pi_6: \mathcal{T}^3_\mathcal{L}(X) \to \mathcal{L}(A) \\
        \pi_7\colon \mathcal{T}^2_\mathcal{L}(X \times A) \to \mathcal{L}^2(X) \quad \pi_8\colon \mathcal{T}^2_\mathcal{L}(X \times A) \to \mathcal{L}^2(A) .
    \end{gather*}
    So, to show that (\ref{diag lift-dbun}) commutes one can check that the following equalities hold:
    \begin{align}
        \pi_1 \circ  \mathcal{T}_\mathcal{L}(\lambda_{X,\mathsf{a}}) \circ \lambda_{X,\mathsf{a}} & = \pi_1 = \pi_1 \circ \ell_{X \times A} \circ \lambda_{X,\mathsf{a}}  \label{eq dbunLalg1}                                                                                 \\
        \pi_2 \circ  \mathcal{T}_\mathcal{L}(\lambda_{X,\mathsf{a}}) \circ \lambda_{X,\mathsf{a}} & = 0_X \circ t_{\mathcal{T}^2_\mathcal{L}(X \times A)} = \pi_2 \circ \ell_{X \times A} \circ \lambda_{X,\mathsf{a}}   \label{eq dbunLalg2}                                  \\
        \pi_3 \circ  \mathcal{T}_\mathcal{L}(\lambda_{X,\mathsf{a}}) \circ \lambda_{X,\mathsf{a}} & = 0_X \circ t_{\mathcal{T}^2_\mathcal{L}(X \times A)} = \pi_3 \circ \ell_{X \times A} \circ \lambda_{X,\mathsf{a}}  \label{eq dbunLalg3}                                   \\
        \pi_4 \circ  \mathcal{T}_\mathcal{L}(\lambda_{X,\mathsf{a}}) \circ \lambda_{X,\mathsf{a}} & = 0_{\mathcal{L}(X)} \circ t_{\mathcal{T}^2_\mathcal{L}(X \times A)} = \pi_4 \circ \ell_{X \times A} \circ \lambda_{X,\mathsf{a}}  \label{eq dbunLalg4}                    \\
        \pi_5 \circ  \mathcal{T}_\mathcal{L}(\lambda_{X,\mathsf{a}}) \circ \lambda_{X,\mathsf{a}} & = 0_X \circ t_{\mathcal{T}^2_\mathcal{L}(X \times A)} = \pi_5 \circ \ell_{X \times A} \circ \lambda_{X,\mathsf{a}}  \label{eq dbunLalg5}                                   \\
        \pi_6 \circ  \mathcal{T}_\mathcal{L}(\lambda_{X,\mathsf{a}}) \circ \lambda_{X,\mathsf{a}} & = 0_{\mathcal{L}(X)} \circ t_{\mathcal{T}^2_\mathcal{L}(X \times A)} = \pi_6 \circ \ell_{X \times A} \circ \lambda_{X,\mathsf{a}}  \label{eq dbunLalg6}                    \\
        \pi_7 \circ  \mathcal{T}_\mathcal{L}(\lambda_{X,\mathsf{a}}) \circ \lambda_{X,\mathsf{a}} & = 0_{\mathcal{L}(X)} \circ t_{\mathcal{T}^2_\mathcal{L}(X \times A)} = \pi_7 \circ \ell_{X \times A} \circ \lambda_{X,\mathsf{a}}  \label{eq dbunLalg7}                    \\
        \pi_8 \circ  \mathcal{T}_\mathcal{L}(\lambda_{X,\mathsf{a}}) \circ \lambda_{X,\mathsf{a}} & = \mathcal{L}(\mathsf{a}^{-1}) \circ \mathsf{a}^{-1} = \nu^{-1}_{A} \circ \mathsf{a}^{-1} = \pi_7 \circ \ell_{X \times A} \circ \lambda_{X,\mathsf{a}}\label{eq dbunLalg8}
    \end{align}
    These are easy to check. For lines (\ref{eq dbunLalg2}) to (\ref{eq dbunLalg6}), we use the fact that \(\mathsf{a}\) is a monoid morphism, while for line (\ref{eq dbunLalg8}) we use the equality \(\mathcal{L}(\mathsf{a}) = \nu_A\). We conclude that \(\mathcal{T}_\mathcal{L}(\lambda_{X,\mathsf{a}}) \circ \lambda_{X,\mathsf{a}} =\ell_{X \times A} \circ \lambda_{X,\mathsf{a}}\), and so, (\ref{diag lift-dbun}) commutes.

    It remains to show the universal property of the lift. Translating diagram (\ref{universaliftdbun}) to our setting, we must show that the following diagram is a pullback:
    \begin{equation}\label{universaliftdbun-L}\begin{gathered}
            \xymatrixcolsep{6pc}\xymatrix{ X \times (A \times A) \ar[d]_-{\pi_1}  \ar[r]^-{\langle 1_X \times \pi_2, \widehat{\zeta} \circ \pi_1 \circ \pi_2 \rangle} & (X \times A) \times \mathcal{L}\left(X \times A\right) \ar[d]^-{\pi_1 \times \mathcal{L}(\pi_1)}  \\
            X \ar[r]_-{\zeta} & X \times A
            }  \end{gathered}\end{equation}
    This is done by a very similar argument to the argument used in the proof of Theorem \ref{theo:linearprojtotangent} to show that (\ref{universalift-L}) is a pullback, so we omit this part of the proof. Lastly, it is not difficult to see that all powers of \(\mathcal{T}_\mathcal{L}\) preserve this pullback. So, we conclude that \((X \times A_\mathsf{a}, \lambda_{X,\mathsf{a}})\) is a differential bundle.

    Suppose now that we are given a morphism \(g\colon X \to X^\prime\) and an \(\mathcal{L}\)-algebra morphism \(f\colon (A, \mathsf{a}) \to (A^\prime, \mathsf{a}^\prime)\). By Lemma \ref{gluing}.(\ref{gluing.proj}), we know that \((g \times f, g)\) is a bundle morphism. So, it remains to show that \(g \times f\) commutes with the lifts. Using the fact that \(f\colon A_\mathsf{a} \to A^\prime_{\mathsf{a}^\prime}\) is a monoid morphism, we compute that:
    \begin{align*}
          \mathcal{T}_\mathcal{L}( g\times f) \circ \lambda_{X,\mathsf{a}}                                                                                                                                                                                                                                                 
         & =~ ( ( g\times f) \times \mathcal{L}(g \times f)) \circ (1_{X \times A} \times \omega^{-1}_{X, A}) \circ \left( \langle 1_X, 0_\mathsf{a} \circ t_X \rangle \times \langle 0_X \circ t_{A}, \mathsf{a}^{-1} \rangle \right)                                                                                      \\
         & =~  (1_{X^\prime \times A^\prime} \times \omega^{-1}_{X^\prime, A^\prime}) \circ \left( ( g\times f) \times \left( \mathcal{L}(g) \times \mathcal{L}(f) \right) \right) \circ \left( \langle 1_X, 0_\mathsf{a} \circ t_X \rangle \times \langle 0_X \circ t_{A}, \mathsf{a}^{-1} \rangle \right)                 \\
         & =~ (1_{X^\prime \times A^\prime} \times \omega^{-1}_{X^\prime, A^\prime}) \circ \left( \left\langle g, f \circ 0_\mathsf{a} \circ t_X \right\rangle \times \langle \mathcal{L}(g) \circ 0_X \circ t_{A}, \mathcal{L}(f) \circ \mathsf{a}^{-1} \rangle  \right)                                                   \\
         & =~ (1_{X^\prime \times A^\prime} \times \omega^{-1}_{X^\prime, A^\prime}) \circ \left( \left\langle g, 0_\mathsf{a^\prime} \circ t_X \right\rangle \times \langle 0_{X^\prime} \circ t_{A}, {\mathsf{a}^\prime}^{-1} \circ f \rangle  \right)                                                                    \\
         & =~ (1_{X^\prime \times A^\prime} \times \omega^{-1}_{X^\prime, A^\prime}) \circ \left( \left\langle g, 0_\mathsf{a^\prime} \circ t_{X^\prime} \circ g \right\rangle \times \langle 0_{X^\prime} \circ t_{A^\prime} \circ f, {\mathsf{a}^\prime}^{-1} \circ f \rangle  \right)                                    \\
         & =~  (1_{X^\prime \times A^\prime} \times \omega^{-1}_{X^\prime, A^\prime}) \circ \left( \left(\left\langle 1_{X^\prime}, 0_\mathsf{a^\prime} \circ t_{X^\prime}  \right\rangle \circ g \right) \times \left( \langle 0_{X^\prime} \circ t_{A^\prime} , {\mathsf{a}^\prime}^{-1} \rangle \circ f \right)  \right) \\
         & =~ (1_{X^\prime \times A^\prime} \times \omega^{-1}_{X^\prime, A^\prime}) \circ \left( \left\langle 1_{X^\prime}, 0_\mathsf{a^\prime} \circ t_{X^\prime}  \right\rangle  \times \langle 0_{X^\prime} \circ t_{A^\prime} , {\mathsf{a}^\prime}^{-1} \rangle   \right) \circ (g \times f)                          \\
         & =~ \lambda_{X^\prime,\mathsf{a}^\prime} \circ (g \times f),
    \end{align*}
    and we conclude that \((g \times f, g)\colon (X \times A_\mathsf{a}, \lambda_{X,\mathsf{a}}) \to (X^\prime \times A^\prime_{\mathsf{a}^\prime}, \lambda_{X^\prime,\mathsf{a}^\prime})\) is a linear bundle morphism.
\end{proof}

Alternatively, one can show that \(\mathcal{L}\)-algebras give differential bundles by pulling back along tangent bundles. Indeed, in a tangent category, for every object \(X\), \((\mathsf{T}(X), \ell_X)\) is a differential bundle. Moreover, assuming the existence and preservation of certain limits, the pullback of the projection of a differential bundle along a morphism is again a differential bundle \cite[Lemma 2.7]{cockettCruttwellDiffBundles}. In the case of a linear assignment \(\mathsf{L}\), \((X \times A_\mathsf{a}, \lambda_{X,\mathsf{a}})\) is in fact the pullback of the tangent bundle \((\mathcal{T}_\mathcal{L}(A), \ell_A)\) along the morphism \(0_\mathsf{a} \circ t_X\colon X \to A\).

In general, \(\mathcal{D}\colon \mathbb{X} \times \mathcal{L}\-\mathsf{ALG} \to \mathsf{DBUN}[(\mathbb{X}, \mathbb{T}_\mathcal{L})]\) will not be an equivalence: an arbitrary differential bundle over \(X\) may not necessarily be of the form \(X \times A\) for some \(\mathcal{L}\)-algebra \(A\). However, in many interesting cases, for example in the category of groups with the abelianization additive assignment, this functor will be an equivalence. We will show in Theorem \ref{thm ker-zero-morph} that this functor is an equivalence as soon as \(\X\) admits zero morphisms and kernels.

Let us apply the above construction to our running examples of linear assignments.

\begin{example} For the terminal additive assignment on a category with finite products, applying this construction to an object \(X\) and the terminal object \(\ast\) results in the trivial differential bundle over \(X\), which is itself \cite[Example 2.4.(i)]{cockettCruttwellDiffBundles}.
\end{example}

\begin{example} For the identity linear (resp.\ additive) assignment on a semi-additive (resp.\ additive) category \(\mathbb{X}\), we get that, for every pair of objects \((X,A)\), \(X \times A\) is a differential bundle over \(X\).
\end{example}

\begin{example} For the additive assignment given by the abelianization \(\mathcal{L}_\ab\) of groups, we get that, for any group \(G\) and any abelian group \(A\), the product \(G \times A\) is a differential bundle over~\(G\), where the differential bundle structure is given as follows:
    \begin{align*}
        \pi_1(g,a) = g &  & \sigma(g,(a,b)) = (g, a+b) &  & \zeta(g) = (g,0) &  & \lambda(g,a) = \left( (g,0), (0,a) \right)
    \end{align*}
\end{example}

Let us now turn our attention to differential objects in our tangent categories: in a cartesian tangent category \((\mathbb{X}, \mathbb{T})\), a \textbf{differential object} \cite[Section 3.1]{cockettCruttwellDiffBundles} is a differential bundle over the terminal object \(\ast\). Note that additive (resp.\ abelian group) bundles over the terminal object correspond precisely to commutative monoids (resp.\ abelian groups), where the projection must be the unique morphism to the terminal object. We will not distinguish between commutative monoids and additive bundles over the terminal object. Therefore, we will characterize differential objects as pairs \((\mathsf{A}, \lambda)\) consisting of a commutative monoid \(\mathsf{A}=(A, +, 0)\) with a morphism \(\lambda\colon A \to \mathcal{T}(A)\). We will refer to linear bundle morphisms over the terminal object as linear morphisms. Explicitly, for two differential objects \((\mathsf{A}, \lambda)\) and \((\mathsf{A}^\prime, \lambda^\prime)\), a \textbf{linear morphism} \(f\colon (\mathsf{A}, \lambda) \to (\mathsf{A}^\prime, \lambda^\prime)\) is a morphism \(f\colon A \to A^\prime\) between the underlying objects such that (\ref{diag dbun-morph}) commutes. We denote by \(\mathsf{DO}[(\mathbb{X}, \mathbb{T}_\mathcal{L})] \coloneq \mathsf{DBUN}[(\mathbb{X}, \mathbb{T}_\mathcal{L})]_\ast\) the category of differential objects and linear morphisms.

Applying Proposition \ref{prop Lalg-to-dbun} to the terminal object, we may build a differential object out of any \(\mathcal{L}\)-algebra: if \((A, \mathsf{a})\) is an \(\mathcal{L}\)-algebra, we define the following morphism:
\begin{equation*}
    \lambda_\mathsf{a} \coloneq \langle 0_\mathsf{a}, \mathsf{a}^{-1} \rangle: A \to A \times \mathcal{L}(A).
\end{equation*}

\begin{corollary}\label{coro LalgDO}Let \(\mathcal{L}\) be a linear assignment on a category \(\mathbb{X}\) with finite products.
    \begin{enumerate}[{\em (i)}]
        \item  For any \(\mathcal{L}\)-algebra \((A, \mathsf{a})\), the pair \((A_\mathsf{a}, \lambda_{\mathsf{a}})\) is a differential object.
        \item For any morphism \(f\colon (A, \mathsf{a}) \to (A^\prime, \mathsf{a}^\prime)\) between \(\mathcal{L}\)-algebras, \(f\colon (A_\mathsf{a}, \lambda_{\mathsf{a}}) \to (A_\mathsf{a}, \lambda_{\mathsf{a}})\) is a linear morphism.
    \end{enumerate}
    This induces a functor \(\mathcal{D}_{\ast}\colon \mathcal{L}\-\mathsf{ALG} \to \mathsf{DO}[(\mathbb{X}, \mathbb{T}_\mathcal{L})]\) defined on objects by \(\mathcal{D}_\ast(A, \mathsf{a}) = (A_\mathsf{a}, \lambda_{\mathsf{a}})\) and on morphisms by \(\mathcal{D}_\ast(f) = f\).
\end{corollary}

We now show that \(\mathcal{D}_{\ast}\colon \mathcal{L}\-\mathsf{ALG} \to \mathsf{DO}[(\mathbb{X}, \mathbb{T}_\mathcal{L})]\) is in fact an isomorphism. We will need to use the following property of differential objects: in a cartesian tangent category \((\mathbb{X}, \mathbb{T})\), if \((\mathsf{A}, \lambda)\) is a differential object, then the following diagram is an equalizer \cite[Lemma 2.14]{cockettCruttwellDiffBundles}:
\begin{equation}\label{eqalizer-diff-obj}\begin{gathered}
        \xymatrixcolsep{5pc}\xymatrix{ A \ar[r]^-{\lambda} & \mathcal{T}(A) \ar@<.6ex>[r]^-{\mathsf{p}_A}  \ar@<-.6ex>[r]_-{0 \circ t_{\mathcal{T}(A)}} & A
        .}
    \end{gathered}\end{equation}
Let \(\mathcal{L}\) be a linear assignment on a category \(\mathbb{X}\) with finite products, and let \((\mathsf{A}, \lambda)\) be a differential object in \((\mathbb{X}, \mathbb{T}_\mathcal{L})\). Define the morphism \(\mathsf{a}_{\mathsf{q}, \lambda}\colon \mathcal{L}(A) \to A\) using the universal property of the above equalizer. That is, \(\mathsf{a}_{\mathsf{q}, \lambda}\) is the unique morphism which makes the following diagram commute:
\begin{equation*}\label{diag a-lambda}\begin{gathered}
        \xymatrixcolsep{5pc}\xymatrix{ A \ar[r]^-{\lambda} & \mathcal{T}_\mathcal{L}(A) = A \times \mathcal{L}(A) \ar@<.6ex>[r]^-{\mathsf{p}_A}  \ar@<-.6ex>[r]_-{0 \circ t_{\mathcal{T}_\mathcal{L}(A)}} & A. \\
        \mathcal{L}(A) \ar@{-->}[u]^-{\mathsf{a}_{\mathsf{q}, \lambda}} \ar[ur]_-{~~~\langle 0 \circ t_A, 1_{\mathcal{L}(A)} \rangle}
        }
    \end{gathered}\end{equation*}

\begin{lemma}\label{lemm DOLalg}Let \(\mathcal{L}\) be a linear assignment on a category \(\mathbb{X}\) with finite products.
    \begin{enumerate}[{\em (i)}]
        \item  For any differential object \((\mathsf{A}, \lambda)\), the pair \((A, {\mathsf{a}_{\mathsf{q}, \lambda}})\) is an \(\mathcal{L}\)-algebra.
        \item For any linear morphism \(f\colon (\mathsf{A}, \lambda) \to (\mathsf{A}^\prime, \lambda^\prime)\), \(f\colon (A, {\mathsf{a}_{\mathsf{q}, \lambda}}) \to (A^\prime, {\mathsf{a}_{\lambda^\prime}})\) is an \(\mathcal{L}\)-algebra morphism.
    \end{enumerate}
    This induced a functor \(\mathcal{D}^{-1}_{\ast}\colon \mathsf{DO}[(\mathbb{X}, \mathbb{T}_\mathcal{L})] \to \mathcal{L}\-\mathsf{ALG} \to \) defined on objects by \(\mathcal{D}^{-1}_\ast(\mathsf{A}, \lambda) = (A, {\mathsf{a}_{\mathsf{q}, \lambda}})\) and on morphisms by \(\mathcal{D}^{-1}_\ast(f) = f\).
\end{lemma}
\begin{proof} Let \((\mathsf{A}, \lambda)\) be a differential object. We first show that \(\mathsf{a}_{\mathsf{q}, \lambda}\) is an isomorphism. Define the morphism \(\mathsf{a}^{-1}_\lambda\colon A \to \mathcal{L}(A)\) by \(\mathsf{a}^{-1}_\lambda \coloneq \pi_2 \circ \lambda\). Observe that, since \(\mathsf{p}_A \circ \lambda = \pi_1 \circ \lambda = 0 \circ t_{\mathcal{T}_\mathcal{L}(A)} \circ \lambda = 0 \circ t_{A}\), by the universal property of the product, we have \(\lambda = \langle  0 \circ t_{A}, \mathsf{a}^{-1}_\lambda \rangle\). We then get:
    \begin{equation*}
        \lambda \circ \mathsf{a}_{\mathsf{q}, \lambda} \circ \mathsf{a}^{-1}_\lambda = \langle  0 \circ t_{\mathcal{L}(A)}, 1_{\mathcal{L}(A)} \rangle \circ \mathsf{a}^{-1}_\lambda = \langle  0 \circ t_{A \times \mathcal{L}(A)}  \circ \mathsf{a}^{-1}_\lambda,  \mathsf{a}^{-1}_\lambda\rangle =  \langle  0 \circ t_{A},  \mathsf{a}^{-1}_\lambda\rangle = \lambda
    \end{equation*}
    Since \(\lambda\) is monic, we have: \(\mathsf{a}_{\mathsf{q}, \lambda} \circ \mathsf{a}^{-1}_\lambda = 1_A\). We can also easily compute that:
    \[ \mathsf{a}^{-1}_\lambda \circ \mathsf{a}_{\mathsf{q}, \lambda} = \pi_2 \circ \lambda \circ \mathsf{a}_{\mathsf{q}, \lambda}  = \pi_2 \circ  \langle  0 \circ t_{\mathcal{L}(A)}, 1_{\mathcal{L}(A)} \rangle = 1_{\mathcal{L}(A)}, \]
    so \(\mathsf{a}^{-1}_\lambda \circ \mathsf{a}_{\mathsf{q}, \lambda} = 1_{\mathcal{L}(A)}\) as well, and thus, \(\mathsf{a}\) is an isomorphism.

    Now, observe that \(\mathcal{T}^2(A) \cong (A \times \mathcal{L}(A)) \times ( \mathcal{L}(A) \times \mathcal{L}^2(A))\). By a slight abuse of notation, we consider the fourth projection \(\pi_4\colon \mathcal{T}^2(A) \to \mathcal{L}^2(A)\). It is easy to check that the following equalities holds:
    \begin{equation*}
        \pi_4 \circ \ell_A = \nu^{-1}_A \circ \pi_2,\qquad \pi_4 \circ \mathcal{T}_\mathcal{L}(\lambda) = \mathcal{L}(\mathsf{a}^{-1}_\lambda) \circ \pi_2 .
    \end{equation*}
    Then, using these equalities and the commutativity of (\ref{diag lift-dbun}), we get:
    \begin{gather*}
        \nu^{-1}_A \circ  \mathsf{a}^{-1}_\lambda =  \nu^{-1}_A \circ \pi_2 \circ \lambda =    \pi_4 \circ \ell_A \circ \lambda = \pi_4 \circ \mathcal{T}_\mathcal{L}(\lambda) \circ \lambda  = \mathcal{L}(\mathsf{a}^{-1}_\lambda) \circ \pi_2 \circ \lambda = \mathcal{L}(\mathsf{a}^{-1}_\lambda) \circ \mathsf{a}^{-1}_\lambda.
    \end{gather*}
    Since \(\mathsf{a}^{-1}_\lambda\) is epic, we get \(\nu^{-1}_A  = \mathcal{L}(\mathsf{a}^{-1}_\lambda)\), and it follows that \(\nu_A = \mathcal{L}(\mathsf{a})\). So, we conclude that  \((A, {\mathsf{a}_{\mathsf{q}, \lambda}})\) is an \(\mathcal{L}\)-algebra.

    Let now \(f\colon (\mathsf{A}, \lambda) \to (\mathsf{A}^\prime, \lambda^\prime)\) be a linear morphism. Since \(f\) commutes with the vertical lifts and since \(\pi_2 \circ \mathcal{T}_\mathcal{L}(f) = \mathcal{L}(f) \circ \pi_2\), we first compute that:
    \begin{equation*}
        \mathcal{L}(f) \circ  \mathsf{a}^{-1}_\lambda = \mathcal{L}(f) \circ \pi_2 \circ \lambda = \pi_2 \circ \mathcal{T}_\mathcal{L}(f) \circ \lambda= \pi_2 \circ \lambda^\prime \circ f = \mathsf{a}^{-1}_{\lambda^\prime} \circ f
    \end{equation*}
    So \(\mathcal{L}(f) \circ  \mathsf{a}^{-1}_\lambda = \mathsf{a}^{-1}_{\lambda^\prime} \circ f\), which implies that \(\mathsf{a}_{\lambda^\prime} \circ \mathcal{L}(f)  =  f \circ  \mathsf{a}_{\mathsf{q}, \lambda}\), and so, finally, \(f\colon (A, {\mathsf{a}_{\mathsf{q}, \lambda}}) \to (A^\prime, {\mathsf{a}_{\lambda^\prime}})\) is an \(\mathcal{L}\)-algebra morphism.
\end{proof}

\begin{theorem}\label{thm diff-obj=l-alg} Let \(\mathcal{L}\) be a linear (resp.\ additive) assignment on a category \(\mathbb{X}\) with finite products. Then, \(\mathsf{DO}[(\mathbb{X}, \mathbb{T}_\mathcal{L})] \cong \mathcal{L}\-\mathsf{ALG}\).
\end{theorem}
\begin{proof} We must show that the functors \(\mathcal{D}_\ast\) and \(\mathcal{D}^{-1}_\ast\) defined respectively in Corollary \ref{coro LalgDO} and Lemma \ref{lemm DOLalg} are mutual inverses. Clearly, on morphisms, one has \(\mathcal{D}_\ast\left(\mathcal{D}^{-1}_\ast(f) \right) = f\) and \(\mathcal{D}^{-1}_\ast\left(\mathcal{D}_\ast(f) \right) = f\). It remains to check that they are inverses on objects as well.

    Starting with a \(\mathcal{L}\)-algebra \((A, \mathsf{a})\), we first notice that, by definition, we have \(\mathsf{a}^{-1}_{\lambda_\mathsf{a}} = \pi_2 \circ \lambda_\mathsf{a} = \mathsf{a}^{-1}\). It follows that \(\mathsf{a}_{\lambda_\mathsf{a}} = \mathsf{a}\), and thus \(\mathcal{D}^{-1}_\ast\left(\mathcal{D}_\ast(A, \mathsf{a}) \right) = (A, \mathsf{a})\).

    On the other hand, let \((\mathsf{A}, \lambda)\) be a differential object. We first show that \(\mathsf{A}\) and \(A_{\mathsf{a}_{\mathsf{q}, \lambda}}\) are the same monoid. Indeed, since \((\lambda, 0)\colon \mathsf{A} \to \mathsf{T}_\mathcal{L}(A)\) is an additive bundle morphism, it is straightforward to see that \(\mathsf{a}^{-1}_\lambda = \pi_2 \circ \lambda\colon \mathsf{A} \to \mathsf{L}(A)\) is a monoid morphism, and since \(\mathsf{a}^{-1}_\lambda\) is an isomorphism, \(\mathsf{a}_{\mathsf{q}, \lambda}\colon \mathsf{L}(A) \to \mathsf{A}\) is a monoid isomorphism. It follows that \(A_{\mathsf{a}_{\mathsf{q}, \lambda}} = \mathsf{A}\). Now, recall that \(\lambda = \langle  0 \circ t_{A}, \mathsf{a}^{-1}_\lambda \rangle\), and that, by definition, \(\lambda_{\mathsf{a}_{\mathsf{q}, \lambda}} = \langle 0_{\mathsf{a}_{\mathsf{q}, \lambda}} \circ t_A,  \mathsf{a}^{-1}_\lambda \rangle\). Since \(A_{\mathsf{a}_{\mathsf{q}, \lambda}} = \mathsf{A}\), we have \(0 = 0_{\mathsf{a}_{\mathsf{q}, \lambda}}\), and so, \(\lambda = \lambda_{\mathsf{a}_{\mathsf{q}, \lambda}}\). Thus, \(\mathcal{D}_\ast\left(\mathcal{D}^{-1}_\ast(\mathsf{A}, \lambda) \right) = (\mathsf{A}, \lambda)\). 
    
    We conclude that \(\mathcal{D}_\ast\) and \(\mathcal{D}^{-1}_\ast\) are inverses of each other, and so, \(\mathsf{DO}[(\mathbb{X}, \mathbb{T}_\mathcal{L})] \cong \mathcal{L}\-\mathsf{ALG}\).
\end{proof}

While differential objects always correspond to \(\mathcal{L}\)-algebras, we mentioned previously that differential bundles over a fixed object do not necessarily correspond to \(\mathcal{L}\)-algebras. However, this will be the case as soon as we have access to zero morphisms and kernels in our category. A category \(\mathbb{X}\) is said to admit \textbf{zero morphisms} if, for every pair of objects \(X\) and \(Y\), there is a distinguished morphism \(0_{X,Y}\colon X \to Y\), and this family of morphisms satisfy \(0_{Y, Y^\prime} \circ f = 0_{X^\prime, Y^\prime} = f \circ 0_{X^\prime, X}\) for every morphism \(f\colon X \to Y\). A category \(\mathbb{X}\) with zero morphisms is said to admit \textbf{kernels} if, for every morphism \(f\colon X \to Y\), the equalizer of \(0_{X,Y}\) and \(f\) exists, which we will denote by \(k_f\colon \mathsf{ker}(f) \to X\).

We will need the use the following universal property of differential bundles: in a cartesian tangent category \((\mathbb{X}, \mathbb{T})\), every differential bundle \((\mathsf{E}, \lambda)\) over \(X\) satisfies the \textbf{Rosický's universality diagram} \cite[Proposition 6]{MacAdamVectorBundles}. That is, the following diagram commutes:
\begin{equation}\label{universal-rosicky-db}\begin{gathered}
        \xymatrixcolsep{5pc}\xymatrix{ E \ar[r]^-{\lambda} \ar[d]_-{\mathsf{q}} & \mathcal{T}(E) \ar[d]^-{\langle \mathcal{T}(\mathsf{q}), \mathsf{p}_E \rangle}  \\
        X \ar[r]_-{\langle \mathsf{z}_X, \zeta \rangle} & \mathcal{T}(X) \times E
        }
    \end{gathered}\end{equation}
Observe that, for differential objects, Rosický's universality diagram is, up to isomorphism, the equalizer diagram (\ref{eqalizer-diff-obj}).

Let \(\mathcal{L}\) be a linear assignment on a category \(\mathbb{X}\) with finite products, and suppose that \(\mathbb{X}\) also admits zero morphisms which are preserved by \(\mathcal{L}\), that is, \(\mathcal{L}(0_{X,Y}) = 0_{\mathcal{L}(X),\mathcal{L}(Y)}\). Let \((\mathsf{E}, \lambda)\) be a differential bundle over \(X\) in \((\mathbb{X}, \mathbb{T}_\mathcal{L})\). Then, using Rosický's universality diagram, let \(\tilde{\lambda}\colon \mathcal{L}(E) \to E\) to be the be the unique morphism making the following diagram commute:
\begin{equation*}\label{}\begin{gathered}
        \xymatrixcolsep{5pc}\xymatrix{ \mathcal{L}(E)  \ar[drr]^-{\langle 0_{\mathcal{L}(E), E}, 1_{\mathcal{L}(E)} \rangle} \ar[ddr]_-{0_{E,X}}  \ar@{-->}[dr]_-{\tilde{\lambda}} \\
        & E \ar[r]^-{\lambda} \ar[d]_-{\mathsf{q}} & \mathcal{T}_\mathcal{L}(E) = E \times \mathcal{L}(E) \ar[d]^-{\langle \mathcal{T}_\mathcal{L}(\mathsf{q}), \mathsf{p}_E \rangle}  \\
        & X \ar[r]_-{\langle \mathsf{z}_X, \zeta \rangle} & \mathcal{T}_\mathcal{L}(X) \times E .
        }
    \end{gathered}\end{equation*}
Suppose now that \(\mathbb{X}\) also admits kernels which are preserved by \(\mathcal{L}\), that is, for every morphism \(f\colon X \to Y\), the morphism \(\mathcal{L}(k_f)\colon \mathcal{L}\left( \mathsf{ker}(f) \right) \to \mathcal{L}(X)\) is the kernel of \(\mathcal{L}(f)\colon \mathcal{L}(X) \to  \mathcal{L}(Y)\). Then, consider the kernel \(k_\mathsf{q}\colon \mathsf{ker}(\mathsf{q}) \to E\) of the projection \(\mathsf{q}\colon E \to X\). Let \(\mathsf{a}_{\mathsf{q}, \lambda}\colon \mathcal{L}\left( \mathsf{ker}(\mathsf{q}) \right) \to \mathsf{ker}(\mathsf{q})\) be the unique morphism which makes the following diagram commute:
\begin{equation*}\label{}\begin{gathered}
        \xymatrixcolsep{5pc}\xymatrix{ \mathsf{ker}(\mathsf{q}) \ar[r]^-{k_\mathsf{q}} & E \ar@<1ex>[r]^-{\mathsf{q}}  \ar@<-1ex>[r]_-{0_{E,X}} & X. \\
        \mathcal{L}\left( \mathsf{ker}(\mathsf{q}) \right) \ar@{-->}[u]^-{\mathsf{a}_{\mathsf{q}, \lambda}} \ar[r]_-{\mathcal{L}(k_\mathsf{q})} & \mathcal{L}(E) \ar[u]_-{\tilde{\lambda}}
        }
    \end{gathered}\end{equation*}
For a linear bundle \((f,g)\colon (\mathsf{E}^\prime, \lambda^\prime) \to (\mathsf{E}^\prime, \lambda^\prime)\), let \(\mathsf{ker}_f\colon \mathsf{ker}(\mathsf{q}) \to \mathsf{ker}(\mathsf{q}^\prime)\) be the unique morphism making the following diagram commute:
\begin{equation*}\label{}\begin{gathered}
        \xymatrixcolsep{5pc}\xymatrix{ \mathsf{ker}(\mathsf{q}^\prime) \ar[r]^-{k_{\mathsf{q}^\prime}} & E^\prime \ar@<1ex>[r]^-{\mathsf{q}^\prime}  \ar@<-1ex>[r]_-{0_{E^\prime,X^\prime}} & X^\prime \\
        \mathsf{ker}(\mathsf{q}) \ar@{-->}[u]^-{\mathsf{ker}_f} \ar[r]_-{k_{\mathsf{q}}} & E \ar[u]_-{f}\ar@<1ex>[r]^-{\mathsf{q}}  \ar@<-1ex>[r]_-{0_{E,X}} & X. \ar[u]_-{g}
        }
    \end{gathered}\end{equation*}

\begin{lemma}\label{lemma L-alg-to-d-bun}Let \(\mathcal{L}\) be a linear assignment on a category \(\mathbb{X}\) with finite products, such that \(\mathbb{X}\) also has zero morphisms and kernels which are preserved by \(\mathcal{L}\).
    \begin{enumerate}[{\em (i)}]
        \item  For any differential bundle \((\mathsf{E}, \lambda)\), the pair \((\mathsf{ker}(\mathsf{q}), {\mathsf{a}_{\mathsf{q},\lambda}})\) is a \(\mathcal{L}\)-algebra.
        \item For any linear bundle morphism \((f,g)\colon (\mathsf{E}^\prime, \lambda^\prime) \to (\mathsf{E}^\prime, \lambda^\prime)\), we have that
              \[
                  \mathsf{ker}_f\colon (\mathsf{ker}(\mathsf{q}), {\mathsf{a}_{\mathsf{q},\lambda}}) \to (\mathsf{ker}(\mathsf{q}^\prime), {\mathsf{a}_{\mathsf{q}^\prime,\lambda^\prime}})
              \]
              is an \(\mathcal{L}\)-algebra morphism.
    \end{enumerate}
    This induces a functor \(\mathcal{D}^\flat\colon \mathsf{DBUN}[(\mathbb{X}, \mathbb{T}_\mathcal{L})] \to \mathbb{X} \times \mathcal{L}\-\mathsf{ALG}\) defined on objects by \(\mathcal{D}^\flat(\mathsf{E}, \lambda) = \left(X, (\mathsf{ker}(\mathsf{q}), {\mathsf{a}_{\mathsf{q},\lambda}}) \right)\) and on morphisms by \(\mathcal{D}(f,g) = (g, \mathsf{ker}_f)\). 
    
    Similarly, for every object \(X\), we also have a functor \(\mathcal{D}^{-1}_X\colon \mathsf{DBUN}[(\mathbb{X}, \mathbb{T}_\mathcal{L})]_X \to  \mathcal{L}\-\mathsf{ALG}\) defined on objects by \(\mathcal{D}^\flat_X(\mathsf{E},\lambda) = (\mathsf{ker}(\mathsf{q}), {\mathsf{a}_{\mathsf{q},\lambda}})\) and on morphisms by \(\mathcal{D}^\flat_X(f) = \mathsf{ker}_f\).
\end{lemma}
\begin{proof} Let \(\lambda_2\colon E \to \mathcal{L}(E)\) be the composite \(\lambda_2 \coloneq \pi_2 \circ \lambda\). Then define the morphism \(\mathsf{a}^{-1}_{\mathsf{q}, \lambda}\colon \mathsf{ker}(\mathsf{q}) \to \mathcal{L}\left( \mathsf{ker}(\mathsf{q}) \right)\) using the universal property of \(\mathcal{L}\left( \mathsf{ker}(\mathsf{q}) \right)\): \(\mathsf{a}^{-1}_{\mathsf{q}, \lambda}\) is the unique morphism which makes the following diagram commute:
    \begin{equation*}\label{}\begin{gathered}
            \xymatrixcolsep{5pc}\xymatrix{  \mathcal{L}\left( \mathsf{ker}(\mathsf{q}) \right) \ar[r]^-{\mathcal{L}(k_\mathsf{q})} & \mathcal{L}(E)  \ar@<.6ex>[rr]^-{\mathcal{L}(\mathsf{q})}  \ar@<-.6ex>[rr]_-{\mathcal{L}(0_{X,Y}) = 0_{\mathcal{L}(X),\mathcal{L}(Y)}} && \mathcal{L}(X). \\
            \mathsf{ker}(\mathsf{q})  \ar@{-->}[u]^-{\mathsf{a}_{\mathsf{q}, \lambda}^{-1}} \ar[r]_-{k_\mathsf{q}} & E \ar[u]_-{\lambda_2}
            }
        \end{gathered}\end{equation*}
    We then compute:
    \begin{equation*}
        \lambda_2 \circ \tilde{\lambda} = \pi_2 \circ \lambda \circ \tilde{\lambda} = \pi_2 \circ \langle 0_{\mathcal{L}(E), E}, 1_{\mathcal{L}(E)} \rangle = 1_{\mathcal{L}(E)},
    \end{equation*}
    so \(\lambda_2 \circ \tilde{\lambda}  = 1_{\mathcal{L}(E)}\), and we get:
    \begin{equation*}
        \mathcal{L}(k_\mathsf{q}) \circ \mathsf{a}^{-1}_{\mathsf{q}, \lambda} \circ \mathsf{a}_{\mathsf{q}, \lambda} = \lambda_2 \circ k_\mathsf{q} \circ \mathsf{a}_{\mathsf{q}, \lambda} =   \lambda_2 \circ \tilde{\lambda} \mathcal{L}(k_\mathsf{q}) =  \mathcal{L}(k_\mathsf{q}).
    \end{equation*}
    Since, by the universal property of the equalizer, \(\mathcal{L}(k_\mathsf{q})\) is monic, it follows that \(\mathsf{a}^{-1}_{\mathsf{q}, \lambda} \circ \mathsf{a}_{\mathsf{q}, \lambda} = 1_{\mathcal{L}\left( \mathsf{ker}(\mathsf{q}) \right)}\). On the other hand, since \((\lambda, \zeta)\colon \mathsf{E} \to \mathsf{T}(E)\) is an additive bundle morphism, \(\pi_1 \circ \lambda = \zeta \circ \mathsf{q}\), and therefore, by the universal property of the product, we have  \(\lambda = \langle \zeta \circ \mathsf{q}, \lambda \rangle\). We then get:
    \begin{equation*}
        \mathsf{q} \circ \tilde{\lambda} \circ \lambda_2 \circ k_\mathsf{q} = 0_{E,X} \circ \lambda_2 \circ k_\mathsf{q} = 0_{\mathsf{ker}(\mathsf{q}),X} = \mathsf{q} \circ k_\mathsf{q},
    \end{equation*}
    and
    \begin{align*}
        \lambda \circ \tilde{\lambda} \circ \lambda_2 \circ k_\mathsf{q} & = \langle 0_{\mathcal{L}(E), E}, 1_{\mathcal{L}(E)} \rangle \circ \lambda_2 \circ k_\mathsf{q} = \langle 0_{\mathcal{L}(E), E} \circ \lambda_2 \circ k_\mathsf{q}, \lambda_2 \circ k_\mathsf{q} \rangle = \langle 0_{\mathsf{ker}(\mathsf{q}), E}, \lambda_2 \circ k_\mathsf{q} \rangle         \\
                                                                         & = \langle \zeta \circ 0_{\mathsf{ker}(\mathsf{q}), X}, \lambda_2 \circ k_\mathsf{q} \rangle   = \langle \zeta \circ \mathsf{q} \circ k_\mathsf{q} , \lambda_2 \circ k_\mathsf{q} \rangle =  \langle \zeta \circ \mathsf{q} , \lambda_2 \rangle \circ k_\mathsf{q} = \lambda \circ k_\mathsf{q}.
    \end{align*}
    Since, by the universal property of the pullback, \(\mathsf{q}\) and \(\lambda\) are jointly monic, it follows that \(\tilde{\lambda} \circ \lambda_2 \circ k_\mathsf{q}  = k_\mathsf{q}\). Using this, we get:
    \begin{equation*}
        k_\mathsf{q} \circ \mathsf{a}_{\mathsf{q}, \lambda} \circ \mathsf{a}^{-1}_{\mathsf{q}, \lambda} = \tilde{\lambda} \circ \mathcal{L}(k_\mathsf{q}) \circ \mathsf{a}^{-1}_{\mathsf{q}, \lambda} = \tilde{\lambda} \circ \lambda_2 \circ k_\mathsf{q} = k_\mathsf{q}.
    \end{equation*}
    Since, by the universal property of an equalizer, \(k_\mathsf{q}\) is monic, it follows that \(\mathsf{a}_{\mathsf{q}, \lambda} \circ \mathsf{a}^{-1}_{\mathsf{q}, \lambda} = 1_{\mathsf{ker}(\mathsf{q})}\). So, \(\mathsf{a}_{\mathsf{q}, \lambda}\) is an isomorphism.

    It remains to show that \(\mathcal{L}(\mathsf{a}_{\mathsf{q}, \lambda}) = \nu_{\mathsf{ker}(\mathsf{q})}\). Translating the commutativity of (\ref{diag lift-dbun}), we have \(\ell_E \circ \lambda = \mathcal{T}_\mathcal{L}(\lambda) \circ \lambda\). Post-composing both sides by the projection \(\pi_4\colon \mathcal{T}^2_\mathcal{L}(E) \to \mathcal{L}\mathcal{L}(E)\), we get:
    \begin{equation}\label{nu=L-lambda_2}
        \nu^{-1}_E = \mathcal{L}(\lambda_2).
    \end{equation}
    From this identity we obtain:
    \begin{equation*}
        \mathcal{L}\mathcal{L}(k_\mathsf{q}) \circ \mathcal{L}(\mathsf{a}^{-1}_{\mathsf{q}, \lambda})  = \mathcal{L}(\lambda_2) \circ \mathcal{L}(k_\mathsf{q}) =\nu^{-1}_E \circ \mathcal{L}(k_\mathsf{q}) =   \mathcal{L}\mathcal{L}(k_\mathsf{q}) \circ   \nu^{-1}_{\mathsf{ker}(\mathsf{q})} .
    \end{equation*}
    Since, by the universal property of the equalizer, \(\mathcal{L}\mathcal{L}(k_\mathsf{q})\) is monic, it follows that \(\mathcal{L}(\mathsf{a}^{-1}_{\mathsf{q}, \lambda}) = \nu^{-1}_{\mathsf{ker}(\mathsf{q})} \), and so \(\mathcal{L}(\mathsf{a}_{\mathsf{q}, \lambda}) = \nu_{\mathsf{ker}(\mathsf{q})}\). We then conclude that \((\mathsf{E}, \lambda)\), \((\mathsf{ker}(\mathsf{q}), {\mathsf{a}_{\mathsf{q},\lambda}})\) is a \(\mathcal{L}\)-algebra.

    Suppose now that \((f,g)\colon (\mathsf{E}^\prime, \lambda^\prime) \to (\mathsf{E}^\prime, \lambda^\prime)\) is a linear bundle morphism. Then \(\lambda^\prime_2 \circ f = \mathcal{L}(f) \circ \lambda_2\), and we get:
    \begin{align*}
        \mathcal{L}(k_{\mathsf{q}^\prime}) \circ \mathcal{L}(\mathsf{ker}_f)  \circ  \mathsf{a}^{-1}_{\mathsf{q}, \lambda} & = \mathcal{L}(f) \circ \mathcal{L}(k_\mathsf{q}) \circ  \mathsf{a}^{-1}_{\mathsf{q}, \lambda} = \mathcal{L}(f) \circ \lambda_2 \circ k_\mathsf{q} = \lambda^\prime_2 \circ f \circ k_\mathsf{q} \\
                                                                                                                           & = \lambda^\prime_2 \circ k_{\mathsf{q}^\prime} \circ \mathsf{ker}_f = \mathcal{L}(k_{\mathsf{q}^\prime}) \circ  \mathsf{a}^{-1}_{\mathsf{q}^\prime, \lambda^\prime} \circ \mathsf{ker}_f .
    \end{align*}
    Since, by the universal property of the equalizer, \(\mathcal{L}(k_{\mathsf{q}^\prime})\) is monic, it follows that \(\mathcal{L}(\mathsf{ker}_f)  \circ  \mathsf{a}^{-1}_{\mathsf{q}, \lambda} = \mathsf{a}^{-1}_{\mathsf{q}^\prime, \lambda^\prime} \circ \mathsf{ker}_f \), and so \(\mathsf{a}_{\mathsf{q}^\prime, \lambda^\prime} \circ \mathcal{L}(\mathsf{ker}_f)  =  \mathsf{ker}_f \circ \mathsf{a}_{\mathsf{q}, \lambda}\). We then conclude that \(\mathsf{ker}_f\) is an \(\mathcal{L}\)-algebra morphism, as desired.
\end{proof}

\begin{theorem}\label{thm ker-zero-morph} Let \(\mathcal{L}\) be a linear (resp.\ additive) assignment on a category \(\mathbb{X}\) with finite products. Suppose that \(\mathbb{X}\) admits zero morphisms and kernels. There is an equivalence of categories \(\mathsf{DBUN}[(\mathbb{X}, \mathbb{T}_\mathcal{L})] \simeq \mathbb{X} \times \mathcal{L}\-\mathsf{ALG}\), and for every object \(X\), \(\mathsf{DBUN}[(\mathbb{X}, \mathbb{T}_\mathcal{L})]_X \simeq \mathcal{L}\-\mathsf{ALG}\).
\end{theorem}
\begin{proof} Let \(X\) be an object of \(\mathbb X\), and  \((A, \mathsf{a})\) be an \(\mathcal{L}\)-algebra. We have \(\mathcal{D}^\flat\left(\mathcal{D}\left( X,(A, \mathsf{a}) \right)\right) = \mathcal{D}^\flat(X \times A_\mathsf{a}, \lambda_{X,\mathsf{a}}) = \left( X, (\mathsf{ker}(\pi_1), {\mathsf{a}_{\pi_1,\lambda_{X,\mathsf{a}}}}) \right)\). Observe that \(\langle 0_{X,Y}, 1_Y \rangle\colon Y \to X \times Y\) is a kernel for the morphism \(\pi_1\), so \(A \cong \mathsf{ker}(\pi_1)\). It is easy to check that this extends to an isomorphism of \(\mathcal{L}\)-algebras \((A, \mathsf{a}) \cong (\mathsf{ker}(\pi_1), {\mathsf{a}_{\pi_1,\lambda_{X,\mathsf{a}}}})\), essentially by construction. Thus, \(\mathcal{D}^\flat\left(\mathcal{D}\left( X,(A, \mathsf{a}) \right)\right) \cong \left( X,(A, \mathsf{a}) \right)\). Moreover, it is straightforward to check that this extends to a natural isomorphism \(\mathcal{D}^\flat \circ \mathcal{D} \cong 1_{\mathbb{X} \times \mathcal{L}\-\mathsf{ALG}}\).

    On the other hand, let \((\mathsf{E}, \lambda)\) be a differential bundle over \(X\). Then,
    \[\mathcal{D}\left(\mathcal{D}^\flat\left( E, \lambda\right)\right) = \mathcal{D}^\flat\left(X, (\mathsf{ker}(\mathsf{q}), {\mathsf{a}_{\mathsf{q},\lambda}}) \right) = (X \times \mathsf{ker}(\mathsf{q})_{\mathsf{a}_{\mathsf{q},\lambda}}, \lambda_{X, \mathsf{a}_{\mathsf{q},\lambda}}).\]
    Let \(\lambda^\flat\colon E \to \mathsf{ker}(\mathsf{q})\) be the unique morphism which makes the following diagram commute:
    \begin{equation*}\label{diag lambda-lambda}\begin{gathered}
            \xymatrixcolsep{5pc}\xymatrix{ \mathsf{ker}(\mathsf{q}) \ar[r]^-{k_\mathsf{q}} & E \ar@<.6ex>[r]^-{\mathsf{q}}  \ar@<-.6ex>[r]_-{0_{E,X}} & X. \\
            E \ar@{-->}[u]^-{\lambda^\flat} \ar[r]_-{\lambda_2} & \mathcal{L}(E) \ar[u]_-{\tilde{\lambda}}
            }
        \end{gathered}\end{equation*}
    Define \(\phi\colon E \to X \times \mathsf{ker}(\mathsf{q})\) as the pairing \(\phi \coloneq \langle \mathsf{q}, \lambda^\flat \rangle\). Then, using Rosický's universality diagram, define the morphism \(\phi^{-1}\colon X \times \mathsf{ker}(\mathsf{q}) \to E\) as the unique morphism which makes the following diagram commute:
    \begin{equation*}\label{}\begin{gathered}
            \xymatrixcolsep{5pc}\xymatrix{ X \times \mathsf{ker}(\mathsf{q}) \ar[drr]^-{\zeta \times (\lambda_2 \circ k_\mathsf{q})} \ar[ddr]_-{\pi_1}  \ar@{-->}[dr]_-{\phi^{-1}} \\
            & E \ar[r]^-{\lambda} \ar[d]_-{\mathsf{q}} & \mathcal{T}_\mathcal{L}(E) = E \times \mathcal{L}(E) \ar[d]^-{\langle \mathcal{T}_\mathcal{L}(\mathsf{q}), \mathsf{p}_E \rangle}  \\
            & X \ar[r]_-{\langle \mathsf{z}_X, \zeta \rangle} & \mathcal{T}_\mathcal{L}(X) \times E .
            }
        \end{gathered}\end{equation*}
    Recall from the proof of Lemma \ref{lemma L-alg-to-d-bun} that \(\tilde{\lambda} \circ \lambda_2 \circ \mathsf{k}_\mathsf{q} = \mathsf{k}_\mathsf{q}\). Then, we get:
    \[  \mathsf{k}_\mathsf{q} \circ \lambda^\flat \circ  \mathsf{k}_\mathsf{q} = \tilde{\lambda} \circ \lambda_2 \circ \mathsf{k}_\mathsf{q} = \mathsf{k}_\mathsf{q}.  \]
    Since \(\mathsf{k}_\mathsf{q}\) is monic, we then get \(\lambda^\flat \circ  \mathsf{k}_\mathsf{q}  = 1_{\mathsf{ker}(\mathsf{q})}\), and
    \begin{align*}
        \mathsf{k}_\mathsf{q} \circ \lambda^\flat \circ  \phi^{-1} &= \tilde{\lambda} \circ \lambda_2 \circ \phi^{-1} = \tilde{\lambda} \circ \pi_2 \circ \lambda \circ \phi^{-1}  = \tilde{\lambda} \circ \pi_2 \circ \left( \zeta \times (\lambda_2 \circ k_\mathsf{q}) \right)  \\
        &= \tilde{\lambda} \circ \lambda_2 \circ \mathsf{k}_\mathsf{q} \circ \pi_2 = \mathsf{k}_\mathsf{q} \circ \pi_2.
    \end{align*}
    Again since \(\mathsf{k}_\mathsf{q}\) is monic, we get \(\lambda^\flat \circ  \phi^{-1} = \pi_2\). It follows that:
    \[ \phi \circ \phi^{-1} =  \langle \mathsf{q}, \lambda^\flat \rangle \circ \phi^{-1}  = \langle \mathsf{q}  \circ \phi^{-1} , \lambda^\flat  \circ \phi^{-1}  \rangle = \langle \pi_1, \pi_2 \rangle = 1_{X \times \mathsf{ker}(\mathsf{q})},  \]
    so \(\phi \circ \phi^{-1} =  1_{X \times \mathsf{ker}(\mathsf{q})}\). In the other direction, we first have:
    \[ \mathsf{q} \circ \phi^{-1} \circ \phi = \pi_1 \circ \phi = \mathsf{q}  .\]
    Recall from the proof of Lemma \ref{lemma L-alg-to-d-bun} that \(\lambda = \langle \zeta \circ \mathsf{q}, \lambda_2 \rangle\) and \(\lambda_2 \circ \tilde{\lambda} = 1_{\mathcal{L}(E)}\). Then, we get:
    \begin{align*}
        \lambda \circ \phi^{-1} \circ \phi & = \left( \zeta \times (\lambda_2 \circ k_\mathsf{q}) \right) \circ \langle \mathsf{q}, \lambda^\flat \rangle = \langle \zeta \circ \mathsf{q}, \lambda_2 \circ k_\mathsf{q} \circ \lambda^\flat \rangle \\
                                           & = \langle \zeta \circ \mathsf{q}, \lambda_2 \circ \tilde{\lambda} \circ \lambda_2 \rangle = \langle \zeta \circ \mathsf{q}, \lambda_2 \rangle = \lambda
    \end{align*}
    Since \( \mathsf{q}\) and \(\lambda\) are jointly monic, we get \(\phi^{-1} \circ \phi = 1_E\).

    It remains to show that \(\phi\) is in fact a linear bundle morphism over \(X\). By definition, we have \(\pi_1 \circ \phi = \mathsf{q}\), so we need to check that \(\phi\) also commutes with the lifts. Recall from the proof of Lemma \ref{lemma L-alg-to-d-bun} that \(\nu^{-1}_E = \mathcal{L}(\lambda_2)\). However, since \(\lambda_2 \circ \tilde{\lambda} = 1_{\mathcal{L}(E)}\), we have \(\mathcal{L}(\lambda_2) \circ \mathcal{L}(\tilde{\lambda}) = 1_{\mathcal{L}\mathcal{L}(E)}\), and so \(\nu_E = \mathcal{L}(\tilde{\lambda})\), or again, \(\mathcal{L}(\tilde{\lambda}) \circ \mathcal{L}(\lambda_2)  = 1_{\mathcal{L}(E)}\). Using this, we get:
    \begin{align*}
        \mathcal{L}(k_\mathsf{q}) \circ \mathsf{a}^{-1}_{\mathsf{q}, \lambda} \circ \lambda^\flat & = \lambda_2 \circ k_\mathsf{q} \circ \lambda^\flat = \lambda_2 \circ \tilde{\lambda} \circ \lambda_2 = \lambda_2                                           \\
                                                                                                  & = \mathcal{L}(\tilde{\lambda}) \circ \mathcal{L}(\lambda_2) \circ \lambda_2 = \mathcal{L}(k_\mathsf{q}) \circ \mathcal{L}(\lambda^\flat) \circ \lambda_2 .
    \end{align*}
    Since \(\mathcal{L}(k_\mathsf{q})\) is monic, we get \(\mathsf{a}^{-1}_{\mathsf{q}, \lambda} \circ \lambda^\flat = \mathcal{L}(\lambda^\flat) \circ \lambda_2\). Now, consider \(\mathcal{T}_\mathcal{L}(X \times  \mathsf{ker}(\mathsf{q}))\) as a quaternary product with projections:
    \begin{align*}
        \pi_1 & \colon \mathcal{T}_\mathcal{L}(X \times  \mathsf{ker}(\mathsf{q})) \to X              & \pi_2 & \colon \mathcal{T}_\mathcal{L}(X \times  \mathsf{ker}(\mathsf{q})) \to \mathsf{ker}(\mathsf{q})                            \\
        \pi_3 & \colon \mathcal{T}_\mathcal{L}(X \times  \mathsf{ker}(\mathsf{q})) \to \mathcal{L}(X) & \pi_4 & \colon \mathcal{T}_\mathcal{L}(X \times  \mathsf{ker}(\mathsf{q})) \to \mathcal{L}\left( \mathsf{ker}(\mathsf{q}) \right).
    \end{align*}
    Then, one can easily check that the following equalities hold:
    \begin{align*}
        \pi_1 \circ \lambda_{\mathsf{a}_{\mathsf{q}, \lambda}} \circ \phi & = \mathsf{q} = \pi_1 \circ \mathcal{T}_\mathcal{L}(\phi) \circ \lambda,                                                                                             \\
        \pi_2 \circ \lambda_{\mathsf{a}_{\mathsf{q}, \lambda}} \circ \phi & = 0_{E, \mathsf{ker}(\mathsf{q})} = \pi_2 \circ \mathcal{T}_\mathcal{L}(\phi) \circ \lambda,                                                                        \\
        \pi_3 \circ \lambda_{\mathsf{a}_{\mathsf{q}, \lambda}} \circ \phi & = 0_{E, \mathcal{L}(X)} = \pi_3 \circ \mathcal{T}_\mathcal{L}(\phi) \circ \lambda ,                                                                                 \\
        \pi_4 \circ \lambda_{\mathsf{a}_{\mathsf{q}, \lambda}} \circ \phi & = \mathsf{a}^{-1}_{\mathsf{q}, \lambda} \circ \lambda^\flat = \mathcal{L}(\lambda^\flat) \circ \lambda_2 = \pi_4 \circ \mathcal{T}_\mathcal{L}(\phi) \circ \lambda.
    \end{align*}
    where, for the fourth and  last equality, we used the fact that \(\mathsf{a}^{-1}_{\mathsf{q}, \lambda} \circ \lambda^\flat = \mathcal{L}(\lambda^\flat) \circ \lambda_2\), which we showed above. We then deduce that \(\phi\colon (\mathsf{E}, \lambda) \to (\mathsf{ker}(\mathsf{q}), {\mathsf{a}_{\mathsf{q},\lambda}})\) is a linear bundle isomorphism. Thus, \(\mathcal{D}\left(\mathcal{D}^\flat\left( E, \lambda\right)\right) \cong (E, \lambda)\). It is not difficult to check that this extends to a natural isomorphism \(\mathcal{D} \circ \mathcal{D}^\flat \cong 1_{\mathsf{DBUN}[(\mathbb{X}, \mathbb{T}_\mathcal{L})]}\).

    We conclude that \(\mathsf{DBUN}[(\mathbb{X}, \mathbb{T}_\mathcal{L})] \simeq \mathbb{X} \times \mathcal{L}\-\mathsf{ALG}\), as desired. Fixing an object \(X\), we obtain natural isomorphisms \(\mathcal{D}^\flat_X \circ \mathcal{D}_X \cong 1_{\mathcal{L}\-\mathsf{ALG}}\) and \(\mathcal{D}_X \circ \mathcal{D}^\flat_X \cong 1_{\mathsf{DBUN}[(\mathbb{X}, \mathbb{T}_\mathcal{L})]_X}\), and thus, \(\mathsf{DBUN}[(\mathbb{X}, \mathbb{T}_\mathcal{L})]_X \simeq \mathcal{L}\-\mathsf{ALG}\) as well.
\end{proof}

An alternative way to prove the desired equivalences of categories is to use Ching's equivalent characterization of differential bundles in terms of wide pullbacks \cite[Theorem 6]{ching2024characterization}, and then using similar arguments as in \cite[Example 17]{ching2024characterization}, which also involves kernels.

Let us apply the above theorem to characterize differential bundles for our tangent structure on groups.

\begin{example}\label{ex dbun-group} \(\Grp\) has zero morphisms and kernels, and the abelianization functors preserves both. By Theorem \ref{thm ker-zero-morph}, differential bundles over a group \(G\) correspond precisely to groups of the form \(G \times A\) for some abelian group \(A\). We then get \(\mathsf{DBUN}[(\Grp, \mathbb{T}_{\mathcal{L}_\ab})] \simeq \Grp \times \Ab\), and for every group \(G\), \(\mathsf{DBUN}[(\Grp, \mathbb{T}_{\mathcal{L}_\ab})]_G \simeq \Ab\).
\end{example}

Here is an example where the conditions of the above theorem fail and where we do not have a correspondence between linear algebras and differential bundles:

\begin{example} Let \(k\) be a field and let \(\mathsf{FVEC}^\flat_k\) be the category of finite-dimensional \(k\)-vector spaces that are \emph{not} of dimension \(1\). Now \(\mathsf{FVEC}^\flat_k\) is an additive category, so the identity functor is an additive assignment. Observe that \(k^3\) is equipped with the structure of a differential bundle over \(k^2\), where:
    \begin{gather*}
        \mathsf{q}(x,y,z) = (x,y) \qquad \sigma( (x,y), (z,w) ) = (x, y, z+w) \qquad \zeta(x,y) = (x,y,0) \\
        \lambda(x,y,z) = \left( ( (x,y,0), (0,0,0) ), ( (0,0,0), (0,0,z) ) \right)
    \end{gather*}
    However, since there is no object \(X\) in \(\mathsf{FVEC}^\flat_k\) such that \(k^2 \times X \cong k^3\) and \(\mathsf{q}\) does not have a kernel in \(\mathsf{FVEC}^\flat_k\), the functor \(\mathcal{D}\) is not an equivalence of categories.
\end{example}

Here is an example where the conditions of the above theorem fail, but where we nevertheless still have a correspondence between linear algebras and differential bundles:

\begin{example}
    Let \(\mathbb{X}\) be a category with finite products, and consider the terminal additive assignment. In this setting, the only differential bundle over an object \(X\) is \(X\) itself. Thus, \(\mathsf{DBUN}[(\mathbb{X}, \mathbb{T}_{\mathcal{L}_\ast})] \cong \mathbb{X}\), and for every object \(X\), \(\mathsf{DBUN}[(\mathbb{X}, \mathbb{T}_{\mathcal{L}_\ast})]\) is trivial. We then trivially have \(\mathsf{DBUN}[(\mathbb{X}, \mathbb{T}_{\mathcal{L}_\ast})] \simeq \mathbb{X} \times \mathcal{L}_\ast\-\mathsf{ALG}\) and \(\mathsf{DBUN}[(\mathbb{X}, \mathbb{T}_{\mathcal{L}_\ast})]_X \simeq \mathcal{L}_\ast\-\mathsf{ALG}\).
\end{example}

\section{Monadic Linear Assignments and Linear Reflectors}

In this section, we consider linear assignments equipped with a monad structure, which we call \textit{monadic} linear assignments. In fact, these will always be idempotent monads. We will show that monadic linear assignments are closely related to reflectors.

\begin{definition}\label{def monadic-L} A \textbf{monadic linear} (resp.\ \textbf{additive}) \textbf{assignment} on a category \(\mathbb{X}\) with finite products is a quintuple \((\mathcal{L}, +, 0, \nu, \eta)\) (resp.\ a sextuple \((\mathcal{L}, +, 0, -, \nu, \eta)\)) consisting of a linear (resp.\ additive) assignment \((\mathcal{L}, +, 0, \nu)\) (resp.\ \((\mathcal{L}, +, 0, -, \nu)\)) such that \((\mathcal{L}, \nu, \eta)\) is a monad. As a shorthand, when there is no confusion, we will denote monadic linear and additive assignments simply by their underlying endofunctor \(\mathcal{L}\).
\end{definition}

Recall that a monad \((\mathcal{L}, \nu, \eta)\) is \textbf{idempotent} \cite[Proposition 4.2.3]{borceux1994handbook} if its multiplication \(\nu_X\colon \mathcal{L}\mathcal{L}(X) \to \mathcal{L}(X)\) is an isomorphism. By definition of a linear assignment, every monadic linear assignment is an idempotent monad. Moreover, for an idempotent monad \(\mathcal L\), postcomposing the monad axiom diagrams by \(\nu^{-1}_X\) gives us the following equalities:
\begin{equation*}
    \mathcal{L}(\nu_X) = \nu_{\mathcal{L}(X)}, \qquad \nu^{-1}_X = \eta_{\mathcal{L}(X)} = \mathcal{L}(\eta_X) .
\end{equation*}
This allows us to characterize monadic linear assignments in terms of idempotent monads:
\begin{proposition} For a category \(\mathbb{X}\) with finite products, a quintuple \((\mathcal{L}, +, 0, \nu, \eta)\) (resp.\ a sextuple\\ \((\mathcal{L}, +, 0, -, \nu, \eta)\)) is a monadic linear (resp.\ additive) assignment if and only if
    \begin{enumerate}[{\em (i)}]
        \item \(\mathcal{L}\) preserves finite products;
        \item \((\mathcal{L}, \nu, \eta)\) is an idempotent monad;
        \item For each object \(X\), the triple \(\mathsf{L}(X) = (\mathcal{L}(X), +_X, 0_X)\) is a commutative monoid (resp.\ the quadruple \(\mathsf{L}(X) = (\mathcal{L}(X), +_X, 0_X, -_X)\) is an abelian group);
        \item For each object \(X\), \(\nu_X\colon \mathsf{L}(\mathcal{L}(X)) \to \mathsf{L}(X)\) is a monoid (resp.\ group) isomorphism.
    \end{enumerate}
\end{proposition}

All of the examples of linear assignments discussed so far are equipped with the structure of a monadic linear assignment:

\begin{example} For any category \(\mathbb{X}\) with finite products, the terminal additive assignment \(\mathcal{L}_\ast\) is monadic, with unit \(\eta_X = t_X\colon X \to \mathcal{L}_\ast(X) = \ast\).
\end{example}

\begin{example} For any semi-additive (resp.\ additive) category \(\mathbb{X}\), the identity linear (resp.\ additive) assignment \(1_\mathbb{X}\) is monadic, with unit \(\eta_X = 1_X\).
\end{example}

\begin{example} The abelianization functor of groups, \(\mathcal{L}_\ab\), seen as a linear assignment, is monadic, and its unit, \(\eta_G\colon G \to \mathcal{L}_\ab(G) = \ab(G)\), is the quotient morphism, \(\eta_G(g) = [g]\).
\end{example}

Of course, not every linear assignment is monadic, here is a counter-example.

\begin{example} Let \(\ab^*(G)\) be the endofunctor on \(\Grp\) which sends a group \(G\) to the torsion subgroup of \(\ab(G)\), which we denote as \(\ab^*(G)\). On can equip this endofunctor with the structure of a linear assignment, which we denote by \(\mathcal{L}_{\ab^*}\). The \(\mathcal{L}_{\ab^*}\)-algebras are the torsion abelian groups. However, \(\mathcal{L}_{\ab^*}\) does not preserve colimits, and so, it is not monadic. Indeed, consider the group morphism \(\phi\colon \Z\to\Z\) defined by \(\phi(n) = 2n\). The cokernel is \(\mathrm{coker}(\phi) \cong \Z/2\Z\), and since \(\Z/2\Z\) a torsion group, we have \(\ab^*(\mathrm{coker}(\phi)) \cong \Z/2\Z\). On the other hand, since \(\Z\) is torsion-free, \(\ab^*(\Z) \cong \mathsf{0}\). Thus, \(\ab^*(m)\) is a zero morphism, and so \(\mathrm{coker}\left( \ab^*(m) \right) \cong \mathsf{0}\).
\end{example}

We now study the linear algebras of a monadic linear assignment. It turns out that, unsurprisingly, these correspond precisely to the algebras over the underlying monad:

\begin{lemma}Let \(\mathcal{L}\) be a monadic linear assignment on a category \(\mathbb{X}\) with finite products. Then, the Eilenberg-Moore category of the monad \((\mathcal{L}, \nu, \eta)\) corresponds precisely to the category \(\mathcal{L}\-\mathsf{ALG}\) described in Definition \ref{def L-alg}.
\end{lemma}
\begin{proof} This argument holds for any idempotent monad. Indeed, for an idempotent monad \((\mathcal{L}, \nu, \eta)\), and an \(\mathcal{L}\)-algebra \((A, \mathsf{a}\colon \mathcal{L}(A) \to A)\) in the usual sense, the morphism \(\mathsf{a}\) is in fact an isomorphism \cite[Proposition 4.2.3]{borceux1994handbook}, and the following equalities hold:
    \begin{equation*}
        \mathcal{L}(\mathsf{a}) = \nu_{A}, \qquad \mathsf{a}^{-1} = \eta_{A},
    \end{equation*}
    so \((A, \mathsf{a})\) is an \(\mathcal{L}\)-algebra in the sense of Definition \ref{def L-alg}. On the other hand, let \((A, \mathsf{a})\) be an \(\mathcal{L}\)-algebra in the sense of Definition \ref{def L-alg}. Then, \(\mathsf{a}\) is an isomorphism and \(\mathcal{L}(\mathsf{a}) = \nu_{A}\), so \(\mathcal{L}(\mathsf{a}) \circ \mathsf{a} = \nu_{A} \circ \mathsf{a}\). Furthermore, one can easily check that \(\mathsf{a} \circ \eta_A \circ \mathsf{a} = \mathsf{a}\), and since \(\mathsf{a}\) is an isomorphism, this implies that \(\mathsf{a} \circ \eta_A= 1_A\). Thus, \((A, \mathsf{a})\) is an \(\mathcal{L}\)-algebra in the usual sense.
\end{proof}

It turns out that, for an idempotent monad \((\mathcal{L}, \nu, \eta)\), the \(\mathcal{L}\)-algebras correspond precisely to the objects \(A\) for which \(\eta_A\colon A \to \mathcal{L}(A)\) is an isomorphism \cite[Corollary 4.2.4]{borceux1994handbook}. Thus, for an idempotent monad, objects have at most one \(\mathcal{L}\)-algebra structure. As such, the Eilenberg-Moore category of an idempotent monad can be associated to a full subcategory of the base category. Let us denote by \(\mathcal{L}\-\mathsf{ALG}^\sharp\) the full subcategory of \(\mathcal{L}\-\mathsf{ALG}\) consisting of objects \(A\) such that \(\eta_A\colon A \to \mathcal{L}(A)\) is an isomorphism. Then, we have an isomorphism of categories \(\mathcal{L}\-\mathsf{ALG}^\sharp \cong \mathcal{L}\-\mathsf{ALG}\), and in particular:

\begin{lemma}\label{lemma monad-L-alg} For a monadic linear assignment \(\mathcal{L}\) on a category \(\mathbb{X}\) with finite products, we have an isomorphism of categories \(\mathcal{L}\-\mathsf{ALG}^\sharp \cong \mathcal{L}\-\mathsf{ALG}\).
\end{lemma}

Idempotent monads are closely connected to the notion of reflective subcategories and reflectors. Recall that a \textbf{reflective subcategory} of a category \(\mathbb{X}\) is a full subcategory \(\mathbb{Y}\) of \(\mathbb{X}\) such that the inclusion function \(\mathcal{I}_\mathbb{Y}\colon \mathbb{Y} \to \mathbb{X}\) admits a left adjoint \(\mathfrak{L}\colon \mathbb{X} \to \mathbb{Y}\), which is called a \textbf{reflector}. For every reflector \(\mathfrak{L}\colon \mathbb{X} \to \mathbb{Y}\), with unit \(\eta_X\colon X \to \mathcal{I}_\mathbb{Y}\left( \mathfrak{L}(X) \right)\) and counit \(\epsilon_Y\colon \mathfrak{L}\left( \mathcal{I}_\mathbb{Y}(Y) \right) \to Y\), the induced monad, \((\mathcal{L}_\mathfrak{L} = \mathcal{I}_\mathbb{Y} \circ \mathfrak{L}, \nu = \mathcal{I}_\mathbb{L}\left( \epsilon_{\mathfrak{L}(-)} \right), \eta) \) on \(\mathbb{X}\) is an idempotent monad. Conversely, given any idempotent monad \((\mathcal{L}, \nu, \eta)\) on \(\mathbb{X}\), the full subcategory \(\mathcal{L}\-\mathsf{ALG}^\sharp\) of \(\mathbb{X}\) is reflective, with reflector \(\mathfrak{L}\colon \mathbb{X} \to \mathcal{L}\-\mathsf{ALG}^\sharp\) defined by \(\mathfrak{L} = \mathcal{L}\). This induces an bijective correspondence between idempotent monads and reflective subcategories \cite[Corollary 4.2.4]{borceux1994handbook}.

We now study this correspondence from the point of view of monadic linear assignments. For a category \(\mathbb{X}\) with finite products, a \textbf{linear} (resp.\ \textbf{additive}) subcategory is a subcategory \(\mathbb{L}\) of \(\mathbb{X}\) such that \(\mathbb{L}\) is a semi-additive (resp.\ additive) category and the inclusion function \(\mathcal{I}_\mathbb{L}\colon \mathbb{L} \to \mathbb{X}\) preserves finite products strictly. In other words, \(\mathbb{L}\) is closed under the product structure of \(\mathbb{X}\), and these products are in fact biproducts in \(\mathbb{L}\).

\begin{definition}\label{def L-reflector} For a category \(\mathbb{X}\) with finite products, a \textbf{linear} (resp.\ \textbf{additive}) \textbf{reflective subcategory} of \(\mathbb{X}\) is a linear (resp.\ additive) subcategory \(\mathbb{L}\) of \(\mathbb{X}\) such that the inclusion functor \(\mathcal{I}_\mathbb{L}\colon \mathbb{L} \to \mathbb{X}\) admits a left adjoint \(\mathfrak{L}\colon \mathbb{X} \to \mathbb{L}\) which preserves finite products. Such a left adjoint is called a \textbf{linear} (resp.\ \textbf{additive}) \textbf{reflector}.
\end{definition}

\begin{theorem} For a category \(\mathbb{X}\) with finite products, there is an bijective correspondence between monadic linear (resp.\ additive) assignments on \(\X\) and linear (resp.\ additive) reflective subcategories of \(\X\).
\end{theorem}
\begin{proof}
    If \(\mathcal{L}\) is a monadic linear (resp.\ additive) assignment on a category \(\mathbb{X}\) with finite products, then \(\mathcal{L}\-\mathsf{ALG}^\sharp\) is a linear (resp.\ additive) reflective subcategory of \(\mathbb{X}\) and \(\mathfrak{L}\colon \mathbb{X} \to \mathcal{L}\-\mathsf{ALG}^\sharp\) is a linear reflector. In the other direction, let \(\mathbb{L}\) be a linear (resp.\ additive) subcategory of a category \(\mathbb{X}\) with finite products, and let \(\mathfrak{L}\colon \mathbb{X} \to \mathbb{L}\) be a linear reflector. As discussed above, this gives us an idempotent monad \((\mathcal{L}_\mathfrak{L}, \nu, \eta)\) on \(\mathbb{X}\). Since \(\mathbb{L}\) is semi-additive (resp.\ additive), every object \(A \in \mathbb{L}\) is canonically a commutative monoid (resp.\ abelian group), with structure morphisms \(\nabla_A\colon A \times A \to A\) and \(b_A\colon \ast \to A\) (and \(i_A\colon A \to A\)), and every morphism in \(\mathbb{L}\) is a monoid (resp.\ group) morphism. In particular, for every object \(X \in \mathbb{X}\), \(\mathfrak{L}(X)\) is a commutative monoid (resp.\ abelian group) and for every morphism \(f\) in \(\mathbb{X}\), \(\mathfrak{L}(f)\) is a monoid (resp.\ group) morphism. As such, we obtain natural transformations \(+_X\colon \mathcal{L}_{\mathfrak{L}}(X) \times \mathcal{L}_{\mathfrak{L}}(X) \to \mathcal{L}_{\mathfrak{L}}(X)\) and \(0_X\colon \ast \to \mathcal{L}_{\mathfrak{L}}(X)\) (and \(-_X\colon \mathcal{L}_{\mathfrak{L}}(X) \to \mathcal{L}_{\mathfrak{L}}(X)\)) defined by:
    \begin{align*}
        +_X \coloneq \nabla_{\mathfrak{L}(X)} &  & 0_X \coloneq b_{\mathfrak{L}(X)} &  & -_X \coloneq i_{\mathfrak{L}(X)},
    \end{align*}
    and which equip \(\mathcal{L}_{\mathfrak{L}}(X)\) with a commutative monoid (resp.\ an abelian group) structure. Moreover, since \(\epsilon_{\mathfrak{L}(X)}\) is a morphism in \(\mathbb{L}\), it is a monoid (resp.\ group) isomorphism, and thus, \(\nu_X\) is also a monoid (resp.\ group) isomorphism. So, \(\mathcal{L}_{\mathfrak{L}}\) is a linear (resp.\ additive) assignment.

    One can easily check that the two processes above, sending monadic linear (resp.\ additive) assignments on \(\X\) to linear reflectors on \(\X\), and vice versa, are mutually inverse.
\end{proof}
It is important to stress that, even if a linear subcategory is reflective, a reflector does not necessarily preserve finite products. Here is an example of a semi-additive (and even additive), reflective subcategory, whose reflector is not linear\footnote{We thank Steve Lack for suggesting this example.}

\begin{example}\label{example Lack} Let \(\GRD\) the category of (small) groupoids. By a slight abuse of notation, we view \(\Ab\) as the full subcategory of \(\GRD\), whose objects are the one-object, abelian groupoids. The inclusion functor \(\Ab\to\GRD\) admits a left adjoint \(\ab^\flat\colon \GRD\to \Ab\). For a groupoid \(\mathbb{G}\), one obtains the abelian group \(\ab^\flat(\mathbb{G})\) by identifying all the objects of \(\mathbb{G}\), adding formal iterations of all morphisms which are not automorphisms, and abelianizing the resulting group. Thus, \(\Ab\) is a reflective subcategory of \(\GRD\). Furthermore, \(\Ab\) is also an additive subcategory of \(\GRD\). However, \(\ab^\flat\) does not preserve finite products. Indeed, let \(\mathbb{I}\) be the groupoid with two objects and morphisms generated by a single isomorphism between these two objects. Then, the product of \(\mathbb{I}\) with itself is a commuting complete diagram with 4 objects:
    \[
        I\colon \diag{0\ar@/^/[r]^{a}&1\ar@/^/[l]},\qquad \mathbb{I} \times \mathbb{I}\colon \vcenter{\diag@R+4pc@C+4pc{(1,0)\ar@/^/[r]^{c}\ar@/^/[d]\ar@/^/[rd]^(.3){f}&(1,1)\ar@/^/[l]\ar@/^/[d]\ar@/^/[ld]\\
        (0,0)\ar@/^/[r]^{a}\ar@/^/[u]^{b}\ar@/^/[ru]^(.3){e}&(0,1)\ar@/^/[l]\ar@/^/[u]^{d}\ar@/^/[lu]}}
    \]
    In the diagrams above, we gave names to one of each pairs of inverse isomorphisms, but we did not label their inverse. For example, in \(\mathbb{I}\), we named \(a\) the isomorphism from \(0\) to~\(1\), and the arrow going the opposite way is \(a^{-1}\). In \(\mathbb{I} \times \mathbb{I}\), we named \(f\) the arrow going from \((1,0)\) to \((0,1)\), and \(e\) the arrow going from \((0,0)\) to \((1,1)\). The abelian group \(\ab^\flat(\mathbb{I})\) is obtained by identifying the objects \(0\) and \(1\), then adding all formal iterations of \(a\) and its inverse. In other words, \(\ab^\flat(\mathbb{I})\) is isomorphic to \(\Z\), freely generated by \(a\). In the abelian group \(\ab(\mathbb{I} \times \mathbb{I})\), all objects of \(\mathbb{I} \times \mathbb{I}\) are identified, all formal iterations of \(a,b,c,d,e,f\) and their inverses are added, but since the diagram of \(\mathbb{I} \times \mathbb{I}\) is commutative, we have certain relations between these generators. For example, \(f=a\circ b^{-1}=d^{-1} \circ c\) and \(e=d \circ a=c \circ b\). Then, one can check that the resulting abelian group is freely generated by \(a,b,c\), and thus, \(\ab(\mathbb{I} \times \mathbb{I})\) is isomorphic to \(\Z^3\). So, \(\ab(\mathbb{I} \times \mathbb{I})\) is not isomorphic to \(\ab(\mathbb{I})\times\ab(\mathbb{I})\). Thus, \(\ab^\flat\) does not preserve products, and is therefore a reflector which is not linear.
\end{example}

\section{Abelianization for Unital Regular Categories}\label{sec abelianization}

The main motivation for this paper is the observation that the category of groups is a tangent category via the abelianization functor. In this section, we generalize this to the setting of unital regular categories. This allows us to provide a bountiful list of novel examples of tangent categories. For an in-depth introduction to unital and regular categories, we invite the reader to see \cite{borceuxbourn04,Bourn1996}.

Let \(\mathbb{X}\) be a category with finite products and zero morphisms (in other words, a pointed category with binary products). For every pair of objects \((X,Y)\), we can define morphisms called the \textbf{quasi-injections}, \(\iota_1\colon X \to X \times Y\) and \(\iota_2\colon Y \to X \times Y\), as follows:
\begin{equation}
    \iota_1 \coloneq \langle 1_X, 0_{X,Y} \rangle, \qquad \iota_2 \coloneq \langle 0_{Y,X}, 1_Y \rangle.
\end{equation}
A \textbf{unital category} \cite[Definition 1.2.5]{borceuxbourn04} is a finitely complete category which admits zero morphisms, and such that the quasi-injections are \textbf{jointly strongly epic} (or equivalently, since we have pullbacks, are \textbf{jointly extremally epic}). This means that whenever we have a monomorphism \(m\) and morphisms \(f\) and \(g\) making the diagram
\[
    \xymatrix{& M \ar[d]^-m\\
    X \ar[ru]^-f \ar[r]_-{\iota_1} & X\times X & X \ar[l]^-{\iota_2} \ar[lu]_-{g}}
\]
commute, \(m\) is an isomorphism. Many examples of unital categories can be found below.

We concluded the previous section by showing that, even if a linear subcategory was reflective, the reflector need not be linear. We show that this situation cannot occur in a unital category:

\begin{proposition}\label{Prop:unital-reflector} Let \(\mathbb{X}\) be a unital category and let \(\mathbb{L}\) be a linear (resp.\ additive) subcategory of \(\mathbb{X}\) which is also reflective. Then any reflector  \(\mathfrak{L}\colon \mathbb{X} \to \mathbb{L}\) is a linear (resp.\ additive) reflector.
\end{proposition}
\begin{proof} We need to show that \(\mathfrak{L}\) preserves finite products. To do so, first note that since we have zero morphisms, the terminal object \(\ast\) of \(\mathbb{X}\) is a zero object both in \(\mathbb{X}\) and \(\mathbb{L}\). Since left adjoints preserve zero objects and zero morphisms, \(t_{\mathfrak{L}(0)}\colon \mathfrak{L}(\ast) \to \ast\) is an isomorphism and \(\mathfrak{L}(0_{X,Y}) = 0_{\mathfrak{L}(X), \mathfrak{L}(Y)}\). We now show that \(\omega_{X,Y}\colon \mathfrak{L}(X \times Y) \to  \mathfrak{L}(X) \times  \mathfrak{L}(Y)\) is an isomorphism. To do so, we will first show that \(\mathfrak{L}(X \times Y)\) is a coproduct of \(\mathfrak{L}(X)\) and \(\mathfrak{L}(Y)\) in \(\mathbb{L}\), with injections \(\mathfrak{L}(\iota_1)\colon \mathfrak{L}(X) \to \mathfrak{L}(X \times Y)\) and \(\mathfrak{L}(\iota_2)\colon \mathfrak{L}(Y) \to \mathfrak{L}(X \times Y)\).

    Consider two morphisms \(f\colon \mathfrak{L}(X) \to Z\) and \(g\colon \mathfrak{L}(Y) \to Z\) in \(\mathbb{L}\). Observe that, since \(\mathbb{L}\) is a full linear subcategory of \(\mathbb X\), the product of two objects of \(\L\) is also a coproduct, and the quasi-injections into this product are the coproduct injections. By the universal property of the coproduct, there is a unique morphism \([f,g]\colon \mathfrak{L}(X) \times \mathfrak{L}(Y) \to Z\) such that \([f,g] \circ \iota_1 = f\) and \([f,g] \circ \iota_2 = g\). Define the morphism \(\overline{[f,g]}\colon \mathfrak{L}(X \times Y) \to Z\) as the composite \(\overline{[f,g]} \coloneq [f,g] \circ \omega_{X,Y}\). One can easily show that \(\omega_{X,Y} \circ \mathfrak{L}(\iota_1) = \iota_1\) and \(\omega_{X,Y} \circ \mathfrak{L}(\iota_2) = \iota_2\). It follows that \(\overline{[f,g]} \circ \mathfrak{L}(\iota_1) = f\) and \(\overline{[f,g]} \circ \mathfrak{L}(\iota_2) = g\). Suppose that there is another morphism \(h\colon \mathfrak{L}(X \times Y) \to Z\) in \(\mathbb{L}\) such that \(h \circ \mathfrak{L}(\iota_1) = f\) and \(h \circ \mathfrak{L}(\iota_2) = g\). Let \((-)^\flat\) be the transpose operation of the adjunction, which takes a morphism of type \(k\colon \mathfrak{L}(E) \to F\) to a morphism of type \(k^\flat\colon E \to F\). We have \(\overline{[f,g]}^\flat \circ \iota_1 = f^\flat =  h^\flat \circ \iota_1\) and \(\overline{[f,g]}^\flat \circ \iota_2 = g^\flat =  h^\flat \circ \iota_2\). However, since \(\iota_1\) and \(\iota_2\) are jointly epic in \(\mathbb{X}\), it follows that \(\overline{[f,g]}^\flat= h^\flat\), and therefore, that \(\overline{[f,g]} = h\). So, we conclude that \(\mathfrak{L}(X \times Y)\) is a coproduct of \(\mathfrak{L}(X)\) and \(\mathfrak{L}(Y)\) in \(\mathbb{L}\).

    Now since \(\mathbb{L}\) is linear, this implies that \(\mathfrak{L}(X \times Y)\) is a product with projections given by the morphisms \(\overline{[1_{\mathfrak{L}(X)}, 0_{\mathfrak{L}(X), \mathfrak{L}(Y)}]}\) and \(\overline{[0_{\mathfrak{L}(Y), \mathfrak{L}(X)}, 1_{\mathfrak{L}(Y)}]}\). However, since \([1_{\mathfrak{L}(X)}, 0_{\mathfrak{L}(X), \mathfrak{L}(Y)}] = \pi_1\) and \([0_{\mathfrak{L}(Y), \mathfrak{L}(X)}, 1_{\mathfrak{L}(Y)}] = \pi_2\), we then get \(\overline{[1_{\mathfrak{L}(X)}, 0_{\mathfrak{L}(X), \mathfrak{L}(Y)}]} = \mathfrak{L}(\pi_1)\) and \(\overline{[0_{\mathfrak{L}(Y), \mathfrak{L}(X)}, 1_{\mathfrak{L}(Y)}]}= \mathfrak{L}(\pi_2)\). Thus, \(\mathfrak{L}(X \times Y)\) is a product in \(\mathbb{L}\), with projections \(\mathfrak{L}(\pi_1)\) and \(\mathfrak{L}(\pi_2)\), and this is equivalent to the fact that \(\omega_{X,Y}\colon \mathfrak{L}(X \times Y) \to  \mathfrak{L}(X) \times  \mathfrak{L}(Y)\) is an isomorphism. We conclude that \(\mathfrak{L}\) preserves finite products, which implies that \(\mathcal{L}_{\mathfrak{L}}\) also preserves finite products, as desired.
\end{proof}

In a unital category, every object admits at most one magma structure \cite[Theorem 1.4.5]{borceuxbourn04}. In particular, every object admits at most one commutative monoid structure. As such, in a unital category, being a commutative monoid is a property of an object rather than an additional structure. Objects in a unital category having this property are referred to as \textbf{commutative objects} \cite[Definition 1.4.1]{borceuxbourn04}. Moreover, every morphism between commutative objects is automatically a monoid morphism. Thus, the category of commutative monoid objects of a unital category is, equivalently, the full subcategory of its commutative objects \cite[Proposition 1.4.11]{borceuxbourn04}. By a slight abuse of notation, for a unital category \(\mathbb{X}\), we will denote by \(\CMON[\X]\) the full subcategory of commutative objects in~\(\mathbb X\), and use the forgetful functor \(\mathcal{U}\colon \CMON[\X]\to \X\) as the inclusion functor. Furthermore, \(\CMON[\X]\) is a linear subcategory of \(\mathbb{X}\).

For a unital category \(\mathbb{X}\) which is also regular and finitely cocomplete, \(\CMON[\X]\) is in fact a linear reflective subcategory. Let us first briefly review the definition of a regular category. Let \(\mathbb X\) be a finitely complete category. Then, for all morphisms \(f\colon Z\to X\), the pullback of \(f\) with itself is called the \textbf{kernel pair} of \(f\). Also, a morphism \(\rho\colon X\to Y\) in~\(\X\) is called a \textbf{regular epimorphism} if there is a pair of morphisms \(\diag{Z\ar@<.5ex>[r]\ar@<-.5ex>[r]&X}\) of which \(\rho\) is the coequalizer. (Note that a regular epimorphism is indeed always an epimorphism.) We then say that \(\X\) is a \textbf{regular category} \cite[Definition A.5.1]{borceuxbourn04} if:
\begin{itemize}
    \item The kernel pair of any morphism admits a coequalizer,
    \item The pullback of any regular epimorphism along any morphism is again a regular epimorphism.
\end{itemize}
If \(\X\) is not only finitely complete, but also finitely cocomplete, then the first condition above is automatically verified, so that \(\X\) is regular if and only if regular epimorphisms are preserved by pullbacks.

In the case where \(\X\) is unital, finitely cocomplete, and regular, then the inclusion of full subcategory \(\mathcal{U}\colon \CMON[\X]\to \X\) discussed above admits a left adjoint \cite[Proposition 1.7.5]{borceuxbourn04} \(\com\colon \X\to \CMON[\X]\), where for an object \(X\), \(\com(X)\) is defined as the coequalizer of the quasi-injections \(\iota_1\colon X \to X \times X\) and \(\iota_2\colon X \to X \times X\). Thus, applying Proposition \ref{Prop:unital-reflector} immediately gives us that:

\begin{proposition}\label{prop unital-linear-reflector} For a finitely cocomplete regular unital category \(\mathbb{X}\), the functor
    \[
        \com\colon \X\to \CMON[\X]
    \]
    is a linear reflector, which in turn induces a monadic linear assignment \(\mathcal{L}_\com\colon \mathbb{X} \to \mathbb{X}\). Moreover, the \(\mathcal{L}_\com\)-algebras correspond precisely to the commutative objects, and so, we have an isomorphism of categories \(\CMON[\X] \cong \mathcal{L}\-\mathsf{ALG}^\sharp\).
\end{proposition}

Applying Theorem \ref{theo:linearprojtotangent}, we may define a cartesian tangent structure on any finitely cocomplete regular unital category. Moreover, since unital categories have zero morphisms and kernels, following Theorem \ref{thm ker-zero-morph}, the differential bundles and differential objects in the resulting tangent category correspond precisely to commutative objects. This gives us the following result:

\begin{theorem}\label{thm unital-tan} Let \(\mathbb{X}\) be a finitely cocomplete regular unital category. Then \((\mathbb{X}, \mathbb{T}_{\mathcal{L}_\com})\) is a cartesian tangent category, where the tangent bundle functor is given by:
\[\mathcal{T}(X) = X \times \com(X)\] 
Moreover, \(\mathsf{DBUN}[(\mathbb{X}, \mathbb{T}_{\mathcal{L}_\com})] \simeq \mathbb{X} \times \CMON[\X]\), and for every object \(X\), \(\mathsf{DBUN}[(\mathbb{X}, \mathbb{T}_{\mathcal{L}_\com})]_X \simeq \CMON[\X]\).
\end{theorem}

Let us now review some new examples of tangent categories built this way. Our first example is a non-Rosický generalization of our main example on the category of groups:

\begin{example}\label{ex monoids} Let \(\mathsf{MON}\) be the category of monoids and let \(\mathsf{CMON}\) be the category of commutative monoids. \(\mathsf{MON}\) is a finitely cocomplete regular unital category whose commutative objects are precisely the commutative monoids, so \(\CMON[\mathsf{MON}] = \mathsf{CMON}\). For a monoid \(M\), \(\com(M)\) is the quotient of \(M \times M\) by the smallest congruence containing \((x,y) \sim (y,x)\) for all \(x,y \in M\). Then, \(\mathsf{MON}\) is a cartesian tangent category with tangent bundle \(\mathcal{T}(M) = M \times \com(M)\), and whose differential bundles (and differential objects) correspond precisely to commutative monoids.
\end{example}

A good source of examples comes from the notion of a variety of universal algebras. Indeed, a variety of universal algebras is always a finitely cocomplete regular category. Then, a variety of universal algebras is also a unital category precisely when it is a \textbf{Jónsson--Tarski variety} \cite[Theorem 1.2.15]{borceuxbourn04}, which essentially means that its signature admits a unique constant \(0\) and a binary operation \(+\) satisfying the equations \(x+0=x=0+x\) \cite[Definition 1.2.14]{borceuxbourn04}. Thus, every Jónsson--Tarski variety admits a cartesian tangent structure given by abelianization:

\begin{example}\label{ex pmag}
    The free Jónsson--Tarski variety is precisely the category of pointed magmas \(\PMAG\). So, \(\PMAG\) is a finitely cocomplete unital regular category, and moreover, the commutative objects in \(\PMAG\) are precisely the commutative monoids. For a pointed magma \(M\) (with binary operation \(\bullet\) and chosen point \(e\)), the commutative monoid \(\com(M)\) is the quotient of \(M\times M\) by the smallest equivalence relation containing \((x,e)\sim (e,x)\) for all \(x\in M\), and which is compatible with the magma structure, in the sense that for all \(a\), \(b\), \(x\), \(y\), \(t\), \(z\in M\), if \((x,y)\sim (t,z)\), then \((a\bullet x,b\bullet y)\sim(a\bullet t,b\bullet z)\) and \((x\bullet a,y\bullet b)\sim(t\bullet a,z\bullet b)\). It turns out that \(\com(M)\) is generated, as a monoid, by classes of elements of the form \((x,e)\) for \(x\in M\). Then, \(\PMAG\) is a cartesian tangent category whose tangent bundle functor satisfies \(\mathcal{T}(M) = M \times \com(M)\), and whose differential bundles (and differential objects) correspond precisely to commutative monoids. Note that \(\MON\) is a full sub-cartesian tangent category of \(\PMAG\). More generally, any pointed variety whose algebras have an underlying pointed magma structure is a finitely cocomplete regular unital category, and thus, admits a cartesian tangent structure.
\end{example}

In order to obtain a Rosický tangent structure induced by abelianization, we need our base category to be not only unital, but strongly unital. In a category \(\mathbb{X}\) with finite products, for every object \(X\), let \(\Delta_X\colon X \to X \times X\) be the canonical diagonal morphism, that is, the morphism defined as follows:
\begin{align}
    \Delta_X = \langle 1_X, 1_X \rangle.
\end{align}
A \textbf{strongly unital category} \cite[Definition 1.8.3]{borceuxbourn04} is a category \(\mathbb{X}\) with finite limits and zero morphisms, such that the quasi-injection \(\iota_1\) (or equivalently \(\iota_2\)) and the diagonal morphism \(\Delta_X\) are jointly strongly epic (or equivalently, since we again have pullbacks, jointly extremally epic). This means that whenever we have a monomorphism \(m\) and morphisms \(f\) and \(g\) making the diagram
\[
    \xymatrix{& M \ar[d]^-m\\
    X \ar[ru]^-f \ar[r]_-{\iota_1} & X\times X & X \ar[l]^-{\Delta_X} \ar[lu]_-{g}}
\]
commute, \(m\) is an isomorphism. This is not the original definition, but one of the equivalent characterizations in~\cite[Theorem~1.8.15]{borceuxbourn04}. Every strongly unital category is in particular unital \cite[Proposition 1.8.4]{borceuxbourn04}. In a strongly unital category \(\mathbb{X}\), every commutative monoid object is an abelian group~\cite[Corollary 1.8.20]{borceuxbourn04}, so commutative objects are called \textbf{abelian objects} \cite[Definition 1.5.4]{borceuxbourn04}, and hence, by a slight abuse of notation, \(\CMON[\X]=\AbG[\mathbb{X}]\). Thus, for a finitely cocomplete regular and strongly unital category \(\mathbb{X}\), we get a left adjoint to the forgetful functor \(\mathcal{U}\colon \mathsf{AbG}[\mathbb{X}]\to \X\), which we will denote by \(\ab\colon \mathbb{X} \to \mathsf{AbG}[\mathbb{X}]\). For all object \(X\), \(\ab(X)\) is defined just like \(\com(X)\) as above. As such, \(\ab\) is an additive reflector, which we refer to as the \textbf{abelianization functor}, which, in turn, induces an additive assignment \(\mathcal{L}_\ab\).

\begin{theorem}\label{thm unital-Ros-tan} Let \(\mathbb{X}\) be a finitely cocomplete, regular, strongly unital category. Then, \((\mathbb{X}, \mathbb{T}_{\mathcal{L}_\ab})\) admits a cartesian Rosický tangent structure, where the tangent bundle functor is \(\mathcal{T}(X) = X \times \ab(X)\). Moreover, \(\mathsf{DBUN}[(\mathbb{X}, \mathbb{T}_{\mathcal{L}_\ab})] \simeq \mathbb{X} \times \AbG[\X]\), and for every object~\(X\), \(\mathsf{DBUN}[(\mathbb{X}, \mathbb{T}_{\mathcal{L}_\ab})]_X \simeq \AbG[\X]\).
\end{theorem}

By \cite[Theorem~1.8.16]{borceuxbourn04}, a variety of algebras is a strongly unital category precisely when it admits a unique constant \(0\) and a ternary operation \(p\) satisfying the equations \(p(x,0,0)=x\) and \(p(x,x,z)=z\). A good source of examples then comes from looking at Mal'tsev varieties, which were originally introduced by Smith in \cite{smith2006mal}. Recall that a \textbf{pointed Mal'tsev variety} \cite[Definition 2.2.1]{borceuxbourn04} is a variety \(\mathsf V\) with a unique constant \(0\) and a ternary operation \(p\) satisfying the equations \(p(x,z,z)=x\) and \(p(x,x,z)=z\).

\begin{example}\label{ex Mal} Every pointed Mal'tsev variety \(\mathsf V\) is a finitely cocomplete, regular, strongly unital category \cite[Corollary 2.2.10]{borceuxbourn04}. By~\cite[Proposition~2.3.8]{borceuxbourn04}, a \(\mathsf V\)-algebra \((A,0,p)\) is abelian exactly when \(p\) is autonomous in the sense of~\cite{Johnstone1989}, that is, for all \(a\), \(b\), \(c\), \(x\), \(y\), \(z\), \(u\), \(v\), \(w\in A\),
    \[
        p(p(a, b, c), p(x, y, z), p(u, v, w)) = p(p(a, x, u), p(b, y, v), p(c, z, w)).
    \]
    This is equivalent to \(p\) being itself a morphism of Mal'tsev algebras. For any \(\mathsf V\)-algebra \(\mathsf A=(A,0,p)\), the abelian object \(\ab(\mathsf A)\) is then an autonomous Mal'tsev algebra whose underlying set is the quotient of \(A\times A\) by the congruence generated by the elements of type
    \[
        \bigl(p(p(a, b, c), p(x, y, z), p(u, v, w)),\; p(p(a, x, u), p(b, y, v), p(c, z, w))\bigr).
    \]
    Thus, a pointed Mal'tsev variety \(\mathsf V\) admits a cartesian Rosický tangent structure with tangent bundle \(\mathcal{T}(\mathsf A) = \mathsf A \times \ab(\mathsf A)\), and whose differential bundles and differential objects correspond precisely to the autonomous \(\mathsf V\)-algebras.
\end{example}

A convenient class of regular strongly unital categories are the semiabelian categories, which were introduced by Janelidze, M{\'a}rki, and Tholen in \cite{Janelidze-Marki-Tholen}. Briefly, a \textbf{semiabelian category} \cite[Definition 5.1.1]{borceuxbourn04} is a category which admits zero morphisms, binary coproducts, and is Barr exact and Bourn protomodular. For an in-depth introduction to semiabelian categories, we refer the reader to \cite{borceuxbourn04,Janelidze-Marki-Tholen}. Every semiabelian category is a finitely cocomplete, regular, strongly unital category (and, in fact, a Mal'tsev category) \cite[Proposition 5.1.2 and 5.1.3]{borceuxbourn04}. For a semiabelian category \(\X\), we then get an abelianization functor \(\ab\colon \mathbb{X} \to \mathsf{AbG}[\mathbb{X}]\) as above, which is left adjoint to the forgetful functor. Moreover, the unit \(\eta_X\) of this adjunction is a normal epimorphism, and its kernel is called the \textbf{commutator} of \(X\) \cite[Definition 2.8.15]{borceuxbourn04}, which we denote by \([X,X]\). Therefore, since a normal epimorphism is always a cokernel of its kernel, we may express the abelianization of an object \(X\) in a semiabelian category as:
\begin{equation*}
    \ab(X)=X/[X,X].
\end{equation*}
One can prove that the reflective subcategory \(\AbG[\mathbb{X}]\) is closed under both subobjects and quotients in \(\X\), which makes it a Birkhoff subcategory~\cite{janelidzekelly96}. This means that, when \(\X\) is a variety of algebras, its subcategory of abelian objects is a subvariety (determined by the equations that characterize abelianness in the given variety). As such, Jónsson--Tarski varieties whose objects have groups as their underlying magma structure, or equivalently, varieties of \(\Omega\)-groups in the sense of Higgins~\cite{Higgins}, are semiabelian~\cite{Janelidze-Marki-Tholen,borceuxbourn04}, which provides us with several examples. We then conclude this paper with some interesting new examples of Rosický tangent categories obtained from semiabelian categories:

\begin{example} \(\Grp\) is semiabelian, and thus, it is finitely cocomplete, regular, and strongly unital. Applying Theorem \ref{thm unital-Ros-tan} to \(\Grp\) results precisely in the Rosický tangent structure introduced in Example \ref{ex tanL-group}.
\end{example}

\begin{example}\label{ex rings} Let \(\mathsf{RNG}\) be the category of non-unital associative rings. Then \(\mathsf{RNG}\) is a semiabelian category, and thus, it is finitely cocomplete, regular, and strongly unital. The abelian objects in \(\mathsf{RNG}\) correspond to \textbf{abelian ring}, which are non-unital rings \(R\) whose multiplication is trivial, that is, \(xy=0\) for all \(x,y \in R\). In other words, abelian rings are essentially abelian group with trivial multiplication, and therefore \(\AbG[\mathsf{RNG}] = \Ab\). For a non-unital associative ring \(R\), its commutator \([R,R]\) is \(R^2\), the set of products of two elements of \(R\): \([R,R]= R^2 = \lbrace xy \,\vert\, x,y \in R \rbrace\). Therefore, the abelianization of \(R\) is \(\ab(R) = R/R^2\). Then, \(\mathsf{RNG}\) admits a cartesian Rosický tangent structure with tangent bundle \(\mathcal{T}(R) = R \times (R/R^2)\), and whose differential bundles and differential objects correspond precisely to abelian groups with trivial multiplication. This example easily generalizes to the commutative case, but also to the non-associative case.
\end{example}

\begin{example} The previous example can also be generalized to the category of algebras over a reduced operad. An algebraic operad \(\P\)\cite{LV} is said to be reduced when \(\P(0)=0\). In this case, the resulting category of \(\P\)-algebras is known to be semiabelian (since it is a variety of \(\Omega\)-groups). The abelianization of a \(\P\)-algebra \(A\) is a \(\P\)-algebra obtained by quotienting \(A\) by the ideal \(A^{\ge2}\) containing all elements of the form \(\mu(a_1,\dots,a_n)\), for \(\mu\) an operation of \(\P\) of arity \(n\ge 2\). The abelian objects are then the \(\P\)-algebras for which all such operations are trivial, which we call \textbf{abelian \(\P\)-algebras}. Therefore, the category of \(\P\)-algebras admits a cartesian Rosický tangent structure with tangent bundle \(\mathcal{T}(A) = A \times (A/A^{\ge2})\) and whose differential bundles and differential objects correspond precisely to the abelian \(\P\)-algebras. We recapture the previous example by taking the operad of non-unital associative rings. For another specific example, consider the operad \(\Lie\) of Lie algebras. Abelian algebras over the operad \(\Lie\) also coincides with the usual notion of abelian Lie algebras. Moreover, for a Lie algebra \(\mathfrak g\), we have that \(\mathfrak g^{\ge 2}=[\mathfrak g,\mathfrak g]\), where the usual bracket notation coincides with the commutator notation here. Therefore, the category of Lie algebras admits a cartesian Rosický tangent structure with tangent bundle \(\mathcal{T}(\mathfrak g) = \mathfrak g \times (\mathfrak g/[\mathfrak g,\mathfrak g])\), and whose differential bundles and differential objects correspond precisely to the abelian Lie algebras.
\end{example}

\begin{example}\label{ex crossed} Crossed modules form a semiabelian category, in fact a variety of \(\Omega\)-groups. Actually, for any semiabelian category \(\mathbb X\), one can define internal crossed modules in~\(\mathbb X\) \cite{janelidze03}, and those still form a semiabelian category. The classical crossed modules then correspond to the internal crossed modules in the category of groups. Abelian objects in the category of crossed modules and their commutator were described in \cite{CCG}. A crossed module \(\partial\colon T\to G\) is abelian precisely when \(T\) and \(G\) are abelian groups and \(G\) acts trivially on \(T\). For a given crossed module \({\partial\colon T\to G}\), its abelianization is the quotient \({\overline{\partial}\colon T/{[G,T]}\to \ab(G)}\) where the commutator \([G,T]\) is generated by all \({}^xtt^{-1}\) for \(x\in G\), \(t\in T\). Then crossed modules for a cartesian Rosický tangent category, where the differential bundles and differential objects correspond to the abelian crossed modules.
\end{example}

\begin{example}\label{ex loop} Let \(\mathsf{LOOP}\) be the category of loops. Then \(\mathsf{LOOP}\) is semiabelian, and the abelian objects are precisely the abelian groups, so \(\AbG[\mathsf{LOOP}] = \Ab\). For a loop \(L\) (with division operator \(\backslash\) and \(\slash\)), the commutator \([L,L]\) is the normal subloop of \(L\) generated by the commutator elements \((xy)\slash(yx)\) and associator elements \(((xy)z)\slash(x(yz))\) for \(x\), \(y\), \(z\in L\) (see~\cite[Section~5.1]{EverVdL4} for further details). Then, \(\mathsf{LOOP}\) admits a cartesian Rosický tangent structure with tangent bundle \(\mathcal{T}(L) = L \times (L/[L,L])\) and whose differential bundles and differential objects correspond precisely to abelian groups.
\end{example}

\begin{example}\label{ex hopf} The category of cocommutative Hopf algebras over a field of characteristic 0 is semiabelian~\cite{GSV-Takeuchi}. In this category, a cocommutative Hopf algebra \(H\) is abelian if and only if it is (bi)commutative \cite{vespaW18}. For a cocommutative Hopf algebra~\(H\), the commutator \([H,H]\) is the normal Hopf subalgebra generated by the elements of the form \(xy-yx\). Thus, recalling that the product of cocommutative Hopf algebras is \(\otimes\), we see that the category of cocommutative Hopf algebras over a field of characteristic 0 admits a cartesian Rosický tangent structure with tangent bundle \(\mathcal{T}(H) = H \otimes (H/[H,H])\) and whose differential bundles and differential objects correspond precisely to the commutative Hopf algebras.
\end{example}

\end{document}